\documentclass[11pt]{article}

\usepackage{amsmath, amsthm, amsfonts, amssymb,sansmath}
\usepackage[letterpaper]{geometry}
\usepackage[colorlinks=true,linkcolor=black,anchorcolor=black,citecolor=black,urlcolor=black,plainpages=false,pdfpagelabels]{hyperref}
\usepackage[capitalize]{cleveref}
\usepackage{subcaption}
\usepackage{graphicx}
\usepackage{algorithm}
\usepackage[noend]{algpseudocode}
\usepackage{enumitem}
\usepackage{tikz}
\usepackage{tikz-cd}
\usepackage[textsize=scriptsize]{todonotes}
\usepackage[margin=40pt,font=small,labelfont=bf]{caption}
\usepackage{xspace}
\usepackage{array, booktabs, multirow, ragged2e}
\usepackage{mathtools} 
\usepackage{comment}
\numberwithin{equation}{section}
\usepackage[sort]{natbib}

\let\originalparagraph\paragraph
\renewcommand{\paragraph}[2][.]{\originalparagraph{#2#1}}
\setitemize{itemsep=0pt}

\newtheorem{theorem}{Theorem}[section]
\usepackage{thmtools}
\newtheorem{proposition}[theorem]{Proposition}
\newtheorem{lemma}[theorem]{Lemma}

\theoremstyle{definition}
\newtheorem{definition}[theorem]{Definition}

\theoremstyle{definition}
\newtheorem{remark}[theorem]{Remark}

\usetikzlibrary{calc}
\usetikzlibrary{decorations.pathreplacing}
\usetikzlibrary{positioning,chains,fit,shapes}
\tikzset{>=stealth}

\newcolumntype{M}[1]{>{\RaggedRight\hspace{0pt}}m{#1}}

\newcommand\category[1]{\ensuremath{\mathbf{#1}}}

\DeclareMathOperator{\im}{im}

\DeclareMathOperator{\coker}{coker}

\DeclareMathOperator{\rank}{rank}

\DeclareMathOperator{\id}{Id}

\DeclareMathOperator{\colim}{colim}
\DeclareMathOperator{\push}{push}

\DeclareMathOperator{\gr}{gr}
\DeclareMathOperator{\Id}{Id}

\DeclareMathOperator{\supp}{supp}

\DeclareMathOperator{\lub}{\vee}

\DeclareMathOperator{\Lan}{Lan}

\newcommand{\cell}{\mathcal{A}}

\newcommand{\kvect}[0]{\category{Vec}}

\newcommand{\fb}{\mathcal{F}(M)}

\newcommand{\B}[1]{\mathcal B\ifthenelse{\equal{#1}{}}{}{(#1)}}
\newcommand{\Bi}[1]{{\mathcal B}_i\ifthenelse{\equal{#1}{}}{}{(#1)}}
\newcommand{\Bopen}[1]{U\ifthenelse{\equal{#1}{}}{}{(#1)}}

\newcommand{\barc}[1]{\mathcal B\ifthenelse{\equal{#1}{}}{}{(#1)}}

\newcommand{\G}[0]{\mathcal E}

\renewcommand{\H}[0]{\mathcal H}
\newcommand{\I}[0]{R}

\renewcommand{\L}[0]{\mathcal L}
\newcommand{\Lplus}[0]{\mathcal{L_{[0,\infty)}}}
\newcommand{\M}[0]{\mathcal M}

\newcommand{\pth}[0]{\Gamma}
\renewcommand{\P}[0]{\mathcal B}

\newcommand{\PCal}[0]{\mathcal P}

\newcommand{\R}[0]{\mathbb R}

\renewcommand{\S}[0]{\mathcal A^\bullet}

\newcommand{\T}[0]{\mathbb T}
\newcommand{\U}[0]{\mathbb U}

\renewcommand{\u}[2]{\mathsf{lifts}_{#1}\ifthenelse{\equal{#2}{}}{}{[#2]}}
\MakeRobust{\u}
\newcommand{\g}[2]{\mathsf{grades}_{#1}\ifthenelse{\equal{#2}{}}{}{[#2]}}
\newcommand{\sig}[2]{\mathsf{sig}_{#1}\ifthenelse{\equal{#2}{}}{}{[#2]}}
\newcommand{\sigInv}[2]{\mathsf{sigInv}_{#1}\ifthenelse{\equal{#2}{}}{}{[#2]}}
\MakeRobust{\sigInv}

\newcommand{\K}[0]{\mathbf{k}}
\newcommand{\Nbb}[0]{\mathbb N}

\newcommand{\Y}[0]{\mathcal K}
\newcommand{\Z}[0]{\mathbb Z}

\newcommand{\tpcolor}[0]{blue!20!white!80!green}

\newcommand{\mapUpdate}[2]{\mathsf{mapUpdate}_{#1}\ifthenelse{\equal{#2}{}}{}{(#2)}}

\newcommand{\low}{\rho}

\newcommand{\XiSp}[0]{S}

\newcommand{\length}{\mathrm{weight}}

\newcommand{\lift}[0]{\mathrm{lift}}

\newcommand{\pfd}{p.f.d.\@\xspace}

\newcommand{\rl}{\mathrm{rl}}
\newcommand{\cl}{\mathrm{cl}}

\newcommand{\ce}{\Delta}  

\title{Fast Queries of Fibered Barcodes}
\author{Michael Lesnick\thanks{University at Albany, Albany, NY, USA; \texttt{mlesnick@albany.edu}} \and Matthew Wright\thanks{St.\ Olaf College, Northfield, MN, USA; \texttt{wright5@stolaf.edu}}}

\begin{document}

\date{}
\maketitle

\begin{abstract}
The \emph{fibered barcode} $\fb$ of a bipersistence module $M$ is the map sending each non-negatively sloped affine line $\ell \subset \R^2$ to the barcode of the restriction of $M$ along $\ell$.  The simplicity, computability, and stability of $\fb$ make it a natural choice of invariant for data analysis applications.  In an earlier preprint \cite{lesnick2015interactive}, we introduced a framework for real-time interactive visualization of  $\fb$, which allows the user to select a single line $\ell$ via a GUI and then plots the associated barcode.  This visualization is a key feature of our software RIVET for the visualization and analysis of bipersistent homology.  Such interactive visualization requires a framework for efficient queries of  $\fb$, i.e., for quickly obtaining the barcode along a given line $\ell$.  To enable such queries, we introduced a novel data structure based on planar line arrangements, called an \emph{augmented arrangement}.  The aim of the present paper is to give an updated and improved exposition of the parts of \cite{lesnick2015interactive} concerning the mathematics of the augmented arrangement and its computation.  Notably, by taking the input to be a minimal presentation rather than a chain complex, we are able to substantially simplify our algorithm for computing the augmented arrangement and its complexity analysis.  
\end{abstract}


\medskip


\section{Introduction}

\subsection{Overview}\label{Sec:Overview}
Topological data analysis (TDA) studies the shape of data using tools from algebraic topology.  Persistent homology, one of the central tools of TDA, provides invariants of data called barcodes by constructing a filtration of topological spaces from the data, taking homology with coefficients in a field, and applying a standard structure theorem to the resulting diagram of vector spaces.  In the last 25 years, persistent homology has been used in hundreds of applications \cite{DONUT} and has been the subject of extensive theoretical work.   Introductions can be found in many textbooks and surveys, e.g., \cite{boissonnat2018geometric,carlsson2009topology,carlsson2021topological,chazal2021introduction,dey2022computational,edelsbrunner2010computational,oudot2015persistence,rabadan2019topological,wasserman2018topological}. 

For many data types, such as noisy point cloud data, a single filtered space does not adequately encode the structure of interest in our data \cite{carlsson2010multiparameter,carlsson2009theory,chazal2011geometric}. 
This motivates the consideration of multiparameter persistent homology \cite{carlsson2009theory,frosini1999size},  which begins with the construction of a \emph{multiparameter filtration}.  For algebraic reasons, defining barcodes in the multiparameter setting is problematic.  This creates significant challenges for the theoretical and practical development of multiparameter persistent homology.  In recent years, there has been a surge in effort on these challenges, leading to substantial progress in theory, computation, and software; see \cref{Sec:Related_Work} for a brief overview of some of this work, and \cite{botnan2022introduction} for a detailed introduction.

Visualization plays a major role in applications of persistent homology, but is more
challenging in the multiparameter setting than in the 1-parameter setting.  In 2015, we posted a preprint  \cite{lesnick2015interactive} which introduced a novel framework for the interactive visualization of bipersistent homology (i.e., 2-parameter persistent homology).  This framework is implemented in the software RIVET \cite{rivet}, which was the first publicly available software for analysis and visualization of bipersistence.  To elaborate, RIVET implements algorithms for several computational problems in the bipersistence pipeline. Specifically, RIVET computes
\begin{itemize}
\item function-Rips and degree-Rips bifiltrations of metric data,
\item a minimal presentation of a homology module of a bifiltration, and
\item  three invariants of such a homology module, namely, the Hilbert function, bigraded Betti numbers, and fibered barcode. 
\end{itemize}
 
The fibered barcode $\fb$ of a bipersistence module $M$ is the map sending each non-negatively sloped affine line $\ell \subset \R^2$ to the barcode of the restriction of $M$ along $\ell$; see \cref{Sec:Fibered_Barcode_Intro}.  RIVET provides a novel visualization of the three invariants, whose key feature is an interactive scheme for exploring the fibered barcode.  

The preprint \cite{lesnick2015interactive} explained in detail the mathematical and algorithmic foundations of our scheme for visualizing the fibered barcode.  The scheme hinges on a framework for efficient queries of the fibered barcode, i.e., for quickly obtaining the barcode along a given line $\ell\subset \R^2$.  In \cite{lesnick2015interactive}, we introduced a data structure based on planar line arrangements, called an \emph{augmented arrangement}, on which we can perform such queries; see \cref{Sec:Augmented_Arrangement_Intro} for an overview of this data structure.  Once the augmented arrangement has been computed, querying it for the barcode along a given line is far more efficient than computing the barcode from scratch.  On typical input,  queries of our data structure in RIVET are fast enough to enable smooth animations of the changing barcode as the query line $\ell$ is varied by the user.  Moreover, the cost of computing and storing the augmented arrangement is small enough for practical use.  See \cref{SimpleAugArrComplexity,Thm:Query_Cost} below for complexity bounds.

The isomorphism type of a bipersistence module $M$ can be efficiently encoded in a \emph{minimal presentation}.  Concretely, this is a matrix with an $\R^2$-valued label for each row and each column, where the labels of the rows and columns give the 0th and 1st bigraded Betti numbers of $M$, respectively; see \cref{Sec:Free_Modules}.  Notably, early versions of RIVET did not compute minimal presentations, and instead computed augmented arrangements directly from the chain complexes of bifiltrations.  The augmented arrangement computation detailed in \cite{lesnick2015interactive} closely parallels what was implemented in the early versions of RIVET, and thus also does not exploit minimal presentations. 

In subsequent work \cite{lesnick2022minimal}, we introduced a novel algorithm for computing minimal presentations of bipersistence modules and implemented this in RIVET.  This algorithm runs in cubic time and quadratic space, and is asymptotically more memory efficient than a classical approach based on Schreyer's algorithm.  Computational experiments reported in \cite{lesnick2022minimal} showed that RIVET's minimal presentation computations are much faster and more scalable than comparable computations in standard open-source commutative algebra software.  Several subsequent works have revisited the problem of computing minimal presentations for bipersistent homology, leading to substantial further advances; see \cref{Sec:Related_Work} for an overview.

It turns out that our  algorithm from \cite{lesnick2015interactive} for computing the augmented arrangement admits major simplifications if the input is assumed to be a minimal presentation rather than a chain complex.   These simplifications have not previously appeared in print.   Moreover, first computing a minimal presentation typically leads to drastic speedups to RIVET's augmented arrangement computation, partly because the chain complexes are typically far larger than the minimal presentations of their homology.  

With this in mind, our aim in this paper is to present a more concise and polished exposition of the scheme from \cite{lesnick2015interactive} for fast queries of fibered barcodes, assuming that the input is a minimal presentation.  
This paper mostly supersedes the mathematical core of \cite{lesnick2015interactive}, but does not replace the material from \cite{lesnick2015interactive} on visualization and software.  Some updates of such material from \cite{lesnick2015interactive}  can be found in the RIVET documentation \cite{rivet}.   

In the remainder of this introduction, we formally define the fibered barcode, give a high-level introduction to the augmented arrangement, state the main results of this paper (which are variants of \cite[Theorems 1.3 and 1.4]{lesnick2015interactive}), and discuss related work on multiparameter persistence.

\subsection{The Fibered Barcode}\label{Sec:Fibered_Barcode_Intro}
We prepare for the definition of fibered barcodes with some background; additional background will be covered in \cref{sec:prelim}.  Let $\kvect$ denote the category of vector spaces over a fixed field $\K$.    
We regard a poset $\mathbb{P}$ as a category with object set $\mathbb{P}$ and  morphisms the relations $p\leq q$.  Define a ($\mathbb{P}$-indexed)  \emph{persistence module} to be a functor $M\colon \mathbb{P}\to \kvect$ for some poset $\mathbb{P}$.  Thus, $M$ consists of a vector space $M_p$ for each $p\in \mathbb{P}$, together with a linear map  $M_{p,q}$ for each $p\leq q\in \mathbb{P}$, such that $M_{p,r}=M_{q,r}\circ M_{p,q}$ for all $p\leq q\leq r$.  In addition, each map $M_{p,p}$ is the identity.  We say $M$ is \emph{pointwise finite dimensional (\pfd)} if $\dim M_p < \infty$ for all $p\in \mathbb{P}$.  
If $\mathbb{P}$ is totally ordered, then we call $M$ a \emph{1-parameter persistence module}, and if $\mathbb{P}=\R^2$ with the usual product partial order (i.e., $(a_1,a_2)\leq (b_1,b_2)$ if and only if $a_1\leq b_1$ and $a_2\leq b_2$), then we call $M$ a \emph{bipersistence module}.

An \emph{interval} $I$ in a totally ordered set $\T$ is a non-empty subset $I$ such that whenever $r<s<t\in \T$ and $r,t\in I$, we also have $s\in I$.  We regard $\T\sqcup\{\infty\}$ as poset with $t<\infty$ for all $t\in \T$.  Given $s<t\in \T\sqcup\{\infty\}$, let $[s,t)\subset \T$ be the interval 
\[ [s,t)=\{x\in \T\mid s\leq x<t\}.\] 
According to a standard structure theorem for persistence modules, for any totally ordered set $\T$ and \pfd 1-parameter persistence module $M\colon \T\to \kvect$, $M$ decomposes in an essentially unique way into indecomposables, where the isomorphism classes of the indecomposables are parameterized by intervals in $\T$ \cite{crawley2012decomposition,zomorodian2005computing}.  
This yields a multiset $\B{M}$ of intervals in $\T$ called the \emph{barcode} of $M$, which determines $M$ up to isomorphism.  If $M$ is finitely presented (see  \cref{Sec:Free_Modules}) then it is \pfd, and each interval of $\B{M}$ is of the form $[s,t)$.

We now formally define the fibered barcode of a \pfd  bipersistence module $M$.   For $\ell\subset\R^2$ an affine line with non-negative slope, we let $M^{\ell}:\ell \to \kvect$ denote the restriction of $M$ to $\ell$.  As the restriction of the partial order on $\R^2$ to $\ell$ is a total order, $M^{\ell}$ has a well-defined \emph{barcode} $\B{M^\ell}$, which is a multiset of intervals in $\ell$.  We call the map $\ell \mapsto \B{M^\ell}$ the \emph{fibered barcode} of $M$, and denote it as $\fb$.  Fibered barcodes were introduced in \cite{cerri2013betti}, and this name for them was introduced in \cite{lesnick2015interactive}.  As observed in \cite{cerri2013betti}, $\fb$ is equivalent to the \emph{rank invariant of $M$} \cite{carlsson2009theory}, i.e., the map $(a\leq b)\mapsto \rank M_{a,b}$.

\begin{remark}
The fibered barcode $\fb$ satisfies two natural stability properties, which we call \emph{internal stability} and \emph{external stability}.   In brief, internal stability  \cite{cerri2011new,landi2014rank} says that at lines $\ell\subset \R^2$ of positive, finite slope, the barcode $\B{M^\ell}$ is continuous with respect to perturbations of $\ell$; see \cite{landi2014rank} for a precise quantitative statement.  External stability \cite{cerri2013betti,landi2014rank} is the inequality \[d_m(\fb,\mathcal{F}(N))\leq d_I(M,N),\] where  $d_I$ is the \emph{interleaving distance}  on bipersistence modules \cite{lesnick2015theory} and $d_m$ is the  \emph{matching distance} on fibered barcodes \cite{cerri2013betti}.  For \pfd bipersistence modules $M$ and $N$, $d_m$ is defined by 
\[d_m(\fb,\mathcal{F}(N))=\sup_{\ell} w_\ell\cdot  d_b(\B{M^{\ell}}, \B{N^{\ell}}),\]
where $\ell$ ranges over affine lines in $\R^2$ of positive, finite slope, $d_b$ is the \emph{bottleneck distance} on barcodes \cite{cohen2007stability} (whose usual definition over $\R$ extends to $\ell$ because $\R$ and $\ell$ are isometric), and for $s_\ell$ the slope of $\ell$,
\[w_\ell=\begin{cases}
\frac{1}{\sqrt{1+s_{\ell}^2}} &\textup{ if }s_{\ell}\geq 1,\\
\frac{1}{\sqrt{1+\frac{1}{s_{\ell}^2}}} &\textup{ if } s_{\ell}<1.\
\end{cases}
\]
\end{remark}

\subsection{The Augmented Arrangement}\label{Sec:Augmented_Arrangement_Intro}
Now let $M$ be a finitely presented bipersistence module.  We let  $\S(M)$ denote the \emph{augmented arrangement} of $M$,  the data structure underlying  our scheme for fast queries of $\fb$.  It consists of three parts:  
\begin{enumerate}
  \item a line arrangement $\cell(M)$ in $[0,\infty)\times \R$ (i.e., a cell decomposition induced by a set of intersecting lines), 
  \item for each face (i.e., 2-cell) $\ce$ of $\cell(M)$, a collection $\P^\ce$ of pairs $(a,b)\in \R^2\times (\R^2\cup\infty)$, called the \emph{barcode template} at $\ce$, and 
  \item a search data structure for performing point location queries in $\cell(M)$.
\end{enumerate}
We give the definition  of $\S(M)$ in \cref{Sec:Augmented_Arrangement}.  Here, to convey the flavor of the definition, we briefly mention a few key ideas underlying the construction.  The line arrangement $\cell(M)$ is formed by dual lines of certain joins of elements in the union of the support of the 0th and 1st Betti numbers of $M$.  Each barcode template $\P^\ce$ is (up to a canonical bijection) the barcode of the restriction of $M$ to a finite, totally ordered subset of $\R^2$.  As discussed in \cref{sec:lineArrangements}, search data structures for log-time point location queries are standard objects in computational geometry.  A choice of such a data structure is required for our asymptotic runtime  bound on queries of $\fb$ (\cref{Thm:Query_Cost}), though in practice our query scheme can be implemented without it.

We next briefly describe how we query $\S(M)$ for the barcodes $\B{M^\ell}$.  
A standard point-line duality, described in \cref{sec:arr_def}, parameterizes the set of affine lines in $\R^2$ with non-negative, finite slope by $[0,\infty)\times \R$.  For simplicity, here we restrict attention to the generic case where the dual of $\ell$ lies in a face $\ce$ of $\cell(M)$; the general case is similar and is treated in  \cref{Sec:QueryingMath}.  

\begin{figure}[ht]
  \begin{center}
    \begin{tikzpicture}
      \draw[<->,ultra thick,color=green!60!black] (0,0) -- (7,4.2);
       \node[right,color=green!60!black] at (7.05,4.2) {$\ell$};   
      
      \draw[line width=3pt,color=violet] (1.667,1) -- (3.167,1.9);
      \draw[line width=3pt,color=violet] (4.5,2.7)-- (6,3.6);
      \node[above left,color=violet] at (2.55,1.45) {$I_1$};   
      \node[above left,color=violet] at (5.35,3.13) {$I_2$};   

      \coordinate (a1) at (1.667,0.3);
      \coordinate (pa1) at (1.667,1);
      \coordinate (b1) at (3.167,0.3);
      \coordinate (pb1) at (3.167,1.9);

      \draw[->,thick,blue!70!white] (a1) -- (1.667,0.9);
       \draw[->,thick,blue!70!white] (b1) -- (3.167,1.8);         
      \fill (a1) circle (0.08);
      \fill (b1) circle (0.08);
      \draw[fill=violet] (pa1) circle (0.08);
      \draw[fill=white] (pb1) circle (0.08);
      \node[below] at (1.667,0.3) {$a_1$};
      \node[below] at (3.167,0.3) {$b_1$};
      
      \coordinate (a2) at (3.167,2.7);
      \coordinate(pa2) at (4.5,2.7);
      \coordinate(b2) at (6,2.7);
      \coordinate(pb2) at (6,3.6);
      \draw[->,thick,blue!70!white] (a2) -- (4.4,2.7);  
      \draw[->,thick,blue!70!white] (b2) -- (6,3.5);       
      \draw[fill=black] (a2) circle (0.08);
      \draw[fill=black] (b2) circle (0.08);
      \draw[fill=violet] (pa2) circle (0.08);
      \draw[fill=white] (pb2) circle (0.08);
      \node[above] at (3.167,2.7) {$a_2$};
      \node[below] at (6.05,2.7) {$b_2$};
      
    \end{tikzpicture}
  \end{center}
  \caption{The barcode $\protect\B{M^\ell}$ is obtained   from the barcode template $\protect\P^\ce$ by ``pushing'' the points of each pair onto $\ell$.  In this example, $\protect\P^\ce=\{(a_1,b_1),(a_2,b_2)\}$, and $\protect\B{M^\ell}$ consists of two disjoint intervals $I_1$, $I_2$.}  
  \label{fig:Barcode_Templates}
\end{figure}
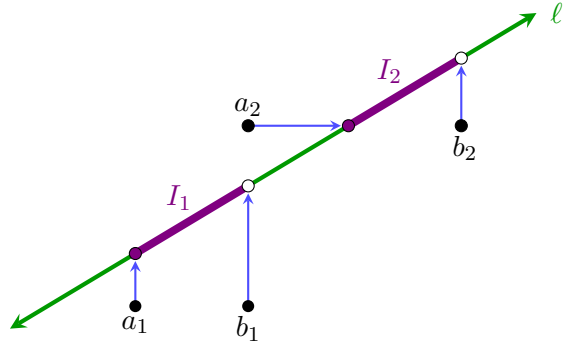

To obtain $\B{M^\ell}$, for each pair $(a,b)\in \P^\ce$, we ``push'' both $a$ and $b$ either upward or rightward onto the line $\ell$, as illustrated in \cref{fig:Barcode_Templates}.  This gives a pair of points $\push_\ell(a), \push_\ell(b)\in \ell$.  We take $\push_\ell(\infty)=\infty$.  See \cref{Sec:CriticalLines} for the formal definition of the map $\push_\ell$.  
\cref{Thm:QueriesMain}  says that
\[ \B{M^\ell} = \{ \,[\push_\ell(a),\,\push_\ell(b)) \mid (a,b) \in \P^\ce\, \}, \]
where $[\push_\ell(a),\,\push_\ell(b))$ is understood to be an interval in $\ell$.  
Thus, to obtain the barcode of $\B{M^\ell}$ it suffices to identify the face $\ce$ via a point location query and then compute $\push_\ell(a)$ and $\push_\ell(b)$ for each $(a,b)\in \P^\ce$.

\subsection{Main Results}\label{Sec:Main_Results_Intro}
We now state our two main theorems, which bound the cost of computing,  storing, and querying $\S(M)$.  
As is standard in computational geometry, we adopt the \emph{real RAM} model of computation, where arithmetic operations over the reals require constant time \cite{preparata2012computational}.  We also assume that all elementary arithmetic operations in the field $\K$ require constant time.  In TDA, most homology computations are done over finite fields, where the latter assumption holds.

Let $\kappa =\kappa_x \cdot \kappa_y$, where $\kappa_x$ and $\kappa_y$ are the number of unique $x$ and $y$ coordinates, respectively, of points in the union of the supports of the 0th and 1st bigraded Betti numbers of $M$; see \cref{Sec:Free_Modules}. 

\begin{theorem}\label{SimpleAugArrComplexity}
    Given a minimal presentation $\eta$ of $M$ with a total of $m$ rows and columns, 
    \begin{enumerate}[label=(\roman*)]
        \item $\S(M)$ has size $O(m\kappa^2)$,
        \item We can compute $\S(M)$ from $\eta$ in $O(m^3 \kappa + \kappa^2 (m + \log\kappa) )$ time and $O(m^2+m\kappa^2)$ memory.
    \end{enumerate}
\end{theorem}

\begin{theorem}\label{Thm:Query_Cost}\mbox{}
\begin{enumerate}[label=(\roman*)]
    \item For a line $\ell$ whose dual lies in a face of $\cell(M)$, we can query $\S(M)$ for $\B{M^\ell}$ in time $O(\log \kappa+|\B{M^\ell}|)$, where $|\B{M^\ell}|$ denotes the number of intervals in $\B{M^\ell}$. 
    \item For all other non-negatively sloped lines $\ell$, we can query $\S(M)$ for $\B{M^\ell}$ in  time $O(\log \kappa+|\B{M^{\ell'}}|)$, for $\ell'$ some arbitrarily small perturbation of $\ell$. 
\end{enumerate}
\end{theorem}

Theorems \ref{SimpleAugArrComplexity}\,(i), \ref{SimpleAugArrComplexity}\,(ii), and \ref{Thm:Query_Cost} are proven in Sections \ref{Sec:Discrete_Barcodes}, \ref{sec:barcodeTemplateCost}, and \ref{Sec:Queries}, respectively.  The term $O(\log \kappa)$ in the bound of \cref{Thm:Query_Cost}\,(i) arises from the point location query used to find the face of the arrangement containing the dual of $\ell$, while the term $|\B{M^\ell}|$ arises from the step of pushing the points of the barcode template onto the line $\ell$.  

\begin{remark}
All practical algorithms for computing persistent homology hinge on variants of Gaussian elimination and have runtime $O(n^3)$, where $n$ is the size of the input.  Moreover, no better bound is known for any such algorithm.  Theoretically, persistent homology can be computed in time  $O(n^{\omega})$, where $\omega$ is the exponent of matrix multiplication \cite{morozov_et_al:LIPIcs.SoCG.2025.68,milosavljevic2011zigzag}, but to our knowledge, no algorithm with sub-cubic complexity has been implemented.  Thus, for computing $\B{M^\ell}$ directly from a minimal presentation of $M$ with a total of $m$ rows and columns, the best available runtime bound is $O(m^{\omega})$ (or $O(m^3)$, among practical algorithms).  This is far larger than the bounds of \cref{Thm:Query_Cost} on the cost of querying $\S(M)$ for $\B{M^\ell}$.  
\end{remark}

\begin{remark}\label{Rem:Quantize}
In the worst case, $\kappa=O(m^2)$.  Thus, the size and time to compute $\S(M)$ are both $O(m^5)$.  
To control the cost of computing and storing $\S(M)$, we can
approximate $M$ by a module $M'$ for which $\kappa$ is a small constant, by snapping the grades of generators and relations onto a grid.  If this snapping changes each coordinate of each generator or relation by at most $\delta$, then we have \[d_m(\mathcal F(M),\mathcal F(M')) \leq d_I(M,M')\leq \delta.\]
Thus, $\mathcal F(M')$ is a $\delta$-approximation of $\mathcal F(M)$ in the matching distance.
\end{remark}

\subsection{Related Work on Multiparameter Persistence}\label{Sec:Related_Work}

The computational study of (multi)graded modules over polynomial rings via Gröbner bases has a long history, e.g., see \cite{cox1998using,kreuzer2000computational,kreuzer2005computational,la1998strategies}.  
To our knowledge, the three earliest TDA works to study the algorithmic aspects of multiparameter persistent homology were: a 2009 paper by Carlsson et al.~\cite{carlsson2009computing}, which applied the classical Gröbner basis algorithms of Buchberger and Schreyer to multiparameter filtrations; a 2011 paper by Cerri et al.~\cite{cerri2011new}, which gave an algorithm for the approximate computation of the matching distance on fibered barcodes; and a 2014 paper by Chacholski et al.~\cite{chacholski2017combinatorial} which studied the problem of expressing the homology of a multiparameter filtration as the homology of a chain complex of free persistence modules.  Our previous preprint \cite{lesnick2015interactive}, on RIVET and fast queries of fibered barcodes, seems to have been the fourth TDA paper on algorithmic aspects of multiparameter persistent homology.  Since the release of these early papers, research activity in the computational aspects of multiparameter persistence has accelerated significantly, and this is now one of the main areas of research in TDA; see \cite{botnan2022introduction} for a recent overview.  In what follows, we mention some of the work in this area most closely related to this paper.

\subsubsection*{Applications of the Ideas Underlying the RIVET Query Scheme}
Our approach in \cite{lesnick2015interactive} of using a line arrangement to decompose the space of affine lines in $\R^2$ has subsequently been used to give the first polynomial-time algorithms to exactly compute the matching distance between fibered barcodes \cite{bjerkevik2021asymptotic,kerber2019matching}.  In addition, Kerber and Wang \cite{kerber_et_al:LIPIcs.SoCG.2026.62} have used  similar ideas to obtain state-of-the-art bounds for the related problem of computing distances between persistent homology transforms.  Fernandes et al.~\cite{fernandes2025computation} have adapted RIVET's query scheme to a scheme for querying \emph{$\gamma$-linear projected barcodes}, a different family of barcodes constructed from a multiparameter persistence module.  
Hickok \cite{hickok2022computing} has also adapted our query scheme to a scheme for querying piecewise-linear families of barcodes indexed by a simplicial complex.

\subsubsection*{Software Based on RIVET}
Aspects of RIVET's visualization scheme have been used in Persistable \cite{scoccola2023persistable}, a software for 2-parameter cluster analysis and visualization.  Persistable does not implement our fast query scheme or provide real-time interactive updates of barcodes, but instead computes barcodes of 1-parameter slices from scratch.  (Persistable does not directly use any code from RIVET.)  RIVET's backend for computing fibered barcodes is used in Oliver Vipond's software for computing multiparameter persistence landscapes (MPLs) \cite{vipond2020multiparameter}; MPLs are vectorizations of multiparameter persistence modules that are equivalent to the restriction of the fibered barcode to lines of slope one.  MPLs were the first of several vectorizations of multiparameter persistence to be introduced and explored in computations \cite{carriere2020multiparameter,cheng2023gril,corbet2019kernel,loiseaux2023framework,loiseaux2023stable}.   A software called Multipers \cite{loiseaux2024multipers} implements computational tools for multiparameter persistence introduced in several recent papers \cite{loiseaux2023framework,loiseaux2022fast,loiseaux2023stable,pmlr-v235-scoccola24a}, and also incorporates code and functionality from many other 
software packages for multiparameter persistence, including both RIVET and Persistable.

\subsubsection*{Applications of Multiparameter Persistence}
Compared to the theoretical and algorithmic foundations of multiparameter persistence, practical applications of multiparameter persistence are much less developed, though we expect that this will change as software improves.  That said, several papers have explored applications of multiparameter persistence to real data.  Two recent publications   \cite{benjamin2024multiscale,vipond2021multi} have applied Vipond's MPL software to problems in immuno-oncology and spatial transcriptomics.  RIVET has also been used to analyze high-dimensional data arising from Wikipedia articles \cite{Wright_Zheng_2020}, to address the virtual screening problem in computational chemistry \cite{keller2018persistent}, and to study molecular generative modules in computational chemistry \cite{schiff2020characterizing}.  Other applications of multiparameter persistence include further work on virtual screening \cite{demir2022todd}, and applications to image analysis \cite{chung2022multi,chung2024morphological}, time series data \cite{chen2022tamp,coskunuzer2024time}, polymer prediction \cite{zhang2024multi}, structural analysis of biomolecules  \cite{xia2015multidimensional}, collective motion of animals \cite{xian2022capturing}, and graph representation learning \cite{chen2024emp}.  These applications  work with multiparameter persistent homology via very simple algebraic invariants, such as the Hilbert function or the barcodes of coordinate-wise 1-parameter slices.  

\subsubsection*{Computing Minimal Presentations}

As noted in \cref{Sec:Overview}, following our initial work on minimal presentation computation in the 2-parameter setting \cite{lesnick2022minimal}, several other works have made significant progress on this problem.  We now give a brief overview.  

Building on an efficient algorithm of Fugacci and Kerber \cite{fugacci_et_al:LIPIcs.SoCG.2019.37} for minimizing a chain complex of free modules, Kerber and Rolle \cite{fugacci2023compression, kerber2021fast} have introduced \emph{mpfree}, an optimized version of our algorithm from \cite{lesnick2022minimal}  that performs far better in practice.  Recently, Luo, Morozov, and Scoccola \cite{morozov2024computing,luo2026computing} have given a specialized algorithm for the case of $0$th homology that runs in nearly linear time, together with an efficient implementation.   Alonso, Kerber, and Pritam \cite{alonso2023filtration} have introduced an algorithm  that replaces a bifiltered clique complex with a smaller one, such that the two bifiltrations have isomorphic homology.  Their experiments demonstrate that for \emph{function-Rips bifiltrations}, this typically yields far smaller bifiltrations and hence leads to much faster minimal presentation computation.  

Bauer, Lenzen, and Lesnick~\cite{bauer2023efficient,bauer2026dualities} have developed a cohomological algorithm for computing a minimal presentation (or resolution) of 2-parameter persistent (co)homology.  This extends the cohomological approach to 1-parameter persistence computation \cite{de2011dualities,bauer2017phat} that underlies Bauer's software Ripser \cite{bauer2021ripser}, a highly efficient code for the special case of Vietoris--Rips filtrations.  The algorithm of \cite{bauer2023efficient} is optimized for function-Rips bifiltrations, and can be combined with the preprocessing algorithm of \cite{alonso2023filtration}.
While this algorithm is competitive with  mpfree on function-Rips bifiltrations and achieves state-of-the-art performance on some examples, the cohomological approach turns out to be more complicated for two parameters than for one, and the benefits are far more modest.  Notably, \cite{bauer2023efficient} does not incorporate two of the main optimizations underlying Ripser, namely \emph{implicit matrix reduction} and \emph{apparent/emergent pairs}; incorporating such optimizations might yield further performance improvements.

\subsubsection*{Barcode-like Representations of (Refinements  of) the Rank Invariant}
As noted above, the fibered barcode is equivalent to the rank invariant.  Since RIVET was introduced, other works have introduced different representations of (refinements of) the rank invariant.  For one example,    
in the special case of function-Rips bifiltrations of $\R$-valued functions on finite metric spaces, 
Cai et al.~\cite{cai2021elder} have defined the \emph{elder-rule-staircase code}, a computable barcode of 0th persistent homology that determines the rank invariant.  
Besides this, a large body of recent work on multiparameter persistence has studied \emph{signed barcodes} of multiparameter persistence modules.  
In the 2-parameter case, a signed barcode is a function from a collection of intervals in $\R^2$ to $\Z$; the word \emph{signed} indicates that the function can take negative values.  
Typically, a signed barcode is equivalent to some other simple invariant, such as the rank invariant or the \emph{generalized rank invariant} \cite{kim2021generalized}; the idea is to encode that invariant in a barcode-like object.  

Two main approaches are used to define signed barcodes: Möbius inversion \cite{botnan2024signed,kim2021generalized} (following ideas of Patel \cite{patel2018generalized}), and relative homological algebra \cite{blanchette2023exact,botnan2024bottleneck,oudot2024stability} (following ideas of Botnan et al.~\cite{botnan2024signed} and Blanchette et al.~\cite{blanchette_brustle_hanson_2022}). 
In particular, both approaches provide invariants equivalent to the rank invariant \cite{botnan2024signed}.  Notably, the \emph{hook signed barcode} obtained via relative homological algebra is Lipschitz stable with respect to appropriate notions of distance \cite{botnan2024bottleneck,oudot2024stability}.  The computation of signed barcodes has been studied in several works, and is viable in some cases \cite{botnan2024signed,chacholski2024koszul,clause2025meta,morozov2021output}.  However, for some proposed constructions of signed barcodes, size, computability, and instability are significant  limitations \cite{botnan2024bottleneck,kim_et_al:LIPIcs.SoCG.2025.64}.  

In the 2-parameter case, signed barcodes provide a mathematically elegant way to fully visualize the rank invariant in a single 2-D plot  \cite{botnan2024signed}.   However, even for very simple examples (e.g., interval modules whose support has a staircase shape), it can be challenging to interpret such a plot; in contrast, RIVET's interactive visualization is straightforward to interpret, but does not fully visualize the rank invariant in a single plot.  Apart from their use in visualization, signed barcodes have been used to construct features of data for machine learning \cite{cheng2023gril,loiseaux2023stable}.

\subsubsection*{Decompositions of Multiparameter Persistence Modules}

Multiparameter persistence modules with finite-dimensional vector spaces satisfy a Krull-Schmidt-Azumaya theorem, i.e., they decompose uniquely into indecomposable summands \cite{botnan2020decomposition}.  Recent work has introduced algorithms and software for computing such decompositions that is optimized for TDA input \cite{dey2025decomposing,dey2019generalized}.  Neither the definition of fibered barcode nor our scheme for fast queries of fibered barcodes depends on the decomposition of a multiparameter persistence module.  However, when a decomposition of the input module $M$ is available, this can potentially be exploited in the setting of our work, either to expedite the computation of the augmented arrangement (as proposed in \cite{dey2025decomposing}), or to provide more refined invariants.  For the latter, some caution is warranted, as decompositions of multiparameter persistence modules are known to be unstable \cite{bauer2025multi}.  Recent work of Bjerkevik has introduced ideas for systematically handling this instability \cite{bjerkevik2025stabilizing}.  

\subsubsection*{Bifiltrations}

Data analysis using 2-parameter persistence, whether using the fibered barcode or other invariants, begins with the construction of  
a bifiltration from data.  Following the initial work of Carlsson and Zomorodian \cite{carlsson2009theory}, which proposed several constructions of bifiltrations, many papers have introduced new bifiltrations of data and studied the theoretical and computational properties of bifiltrations, e.g., 
\cite{alonso2024decomposition,alonso2024delaunay,alonso2023filtration,chung2024morphological,de2022valueoffset}.  Much of this work has focused on density-sensitive bifiltrations of point cloud and metric data \cite{alonso2024sparse,blaser2024core,blumberg2024stability,buchet2024sparse,corbet2023computing,edelsbrunner2021multi,edelsbrunner2020simple,hellmer2024density,lesnick2024nerve,lesnick2024sparse,lesnick2015interactive,sheehy2012multicover}, leading to substantial advances in our understanding of such bifiltrations, as well as new computational tools.

\subsection{Outline}
\cref{sec:prelim} covers preliminaries needed for the rest of the paper, including presentations, resolutions, Betti numbers, barcode computation, updates to $RU$-decompositions, and line arrangements.  \cref{Sec:Augmented_Arrangement} defines the augmented arrangement $\S(M)$ and proves  \cref{SimpleAugArrComplexity}\,(i), which bounds the size of $\S(M)$.  \cref{Sec:QueryingMath} gives our algorithm for querying $\S(M)$, establishes its correctness, and proves \cref{Thm:Query_Cost}, which bounds the time complexity of a query.  \cref{Sec:Computing} gives our  algorithm for computing $\S(M)$ and  \cref{sec:barcodeTemplateCost} bounds its time and memory cost, thereby proving \cref{SimpleAugArrComplexity}\,(ii).

\section{Mathematical Preliminaries}\label{sec:prelim}

\subsection{Free Persistence Modules and Presentations}\label{Sec:Free_Modules}
Given a poset $\mathbb P$, a persistence module $M\colon \mathbb{P}\to \kvect$, an element $p\in \mathbb{P}$, and a vector $v\in M_p$, we write $\gr(v)=p$.  A \emph{set of generators} for $M$ is a set $G\subset \bigsqcup_{p\in \mathbb{P}} M_p$ such that for any $p\in \mathbb{P}$ and $v\in M_p$, we have \[v=\sum_{i=1}^n c_i M_{\gr(g_i),p}(g_i)\] for some $g_1,\ldots,g_n\in G$ and $c_1,\ldots,c_n\in \K$.  

For $p\in \mathbb{P}$, let $\K^{\langle p\rangle}$ denote the persistence module given by 
\begin{align*}
    \K^{\langle p\rangle}_q &=
        \begin{cases}
            \K &\text{if } p\leq q, \\
            0 &\text{otherwise,}
        \end{cases}
    & \K^{\langle p\rangle}_{q,r}=
        \begin{cases}
            \id_\K &{\text{if } p\leq q},\\
            0 &{\text{otherwise.}}
        \end{cases}
\end{align*}
We say a persistence module $F\colon \mathbb{P}\to \kvect$ is \emph{free} if there exists a multiset $\mathcal G$ of elements in $\mathbb{P}$ such that $F\cong \oplus_{p\in \mathcal G}\, \K^{\langle p\rangle}$.  For $F$ free, we define a \emph{basis} of $F$ to be  a minimal set of generators for $F$.  For $B$ a basis of a finitely generated free module $F$, define $\beta^F\colon \mathbb{P}\to \Nbb$ by 
\[\beta^F(p)=|\{b\in B\mid \gr(b) = p\}|.\]  
Though bases of free modules are usually not unique, a straightforward linear algebra argument shows that $\beta^F$ is independent of the choice of $B$.  

A \emph{morphism} of $\mathbb{P}$-indexed persistence modules $\eta\colon M\to N$ is a natural transformation, i.e., a collection of linear maps $\{\eta_p\colon M_p\to N_p\}_{p\in \mathbb{P}}$ such that $\eta_q\circ M_{p,q}=N_{p,q}\circ \eta_p$ for all $p\leq q\in \mathbb{P}$.  With this definition of morphism, the persistence modules over a fixed poset form an abelian category, denoted $\kvect^{\mathbb{P}}$.

Given a morphism $\eta\colon F\to F'$ of finitely generated free $\mathbb{P}$-indexed persistence modules and ordered bases $B$ and $B'$ for $F$ and $F'$, respectively, we represent $\eta$ via a matrix $[\eta]$ with coefficients in the field $\K$, with each row and column labeled by an element of $\mathbb{P}$.  To do so, note that $\eta$ induces a map $\colim \eta\colon \colim F\to \colim F'$, and $B,B'$ induce bases for $\colim F$ and $\colim F'$, respectively.  We take $[\eta]$ to be the usual matrix representation of $\colim \eta$ with 
respect to these induced bases, with the $i^{\mathrm{th}}$ row given the label $\gr{B'_i}$, and the $j^{\mathrm{th}}$ column given the label $\gr{B_j}$.  Regarding the morphism $\eta$ as a functor $\{0<1\}\to \kvect^{\mathbb P}$, it is easily checked that $[\eta]$ determines $\eta$ up to natural isomorphism.

A \emph{presentation} of a persistence module $M\colon \mathbb{P}\to \kvect$ is a morphism $\eta\colon F\to F'$ of free persistence modules such that $M\cong \coker \eta$, where $\coker \eta \coloneqq F'/\im \eta$ is the cokernel of $\eta$.
$M$ is \emph{finitely presented} if there exists a presentation $\eta\colon F\to F'$ of $M$ with $F$ and $F'$ finitely generated.  The above discussion makes clear that we can represent such $\eta$ with respect to a choice of bases for $F$ and $F'$ as a matrix with $\mathbb P$-valued labels for the rows and columns; we call this labeled matrix a \emph{presentation matrix} for $M$, or simply (abusing terminology slightly) a presentation.

We say a presentation $\eta\colon F\to F'$ of $M$ is \emph{minimal} if any presentation $\mu \colon E
\to E'$ of $M$ is naturally isomorphic to $\eta\oplus \nu$ for some $\nu$ of form
\[\nu \colon Y\oplus Z\xrightarrow{\Id_{Y}\oplus\, 0} Y, \]
where $Y$ and $Z$ are free.  
By Azumaya's theorem \cite{azumaya1950corrections}, a minimal presentation is unique up to natural isomorphism.  If $M\colon \R^d\to \kvect$ is a finitely presented persistence module, where $\R^d$ is given the product partial order, then $M$ has a minimal presentation $\eta\colon F\to F'$ with $F$ and $F'$ finitely generated; the analogous result for finitely generated $\Z$-graded $\K[x_1,x_2,\dots,x_d]$-modules is proven in \cite{peeva2011syzygies}, and this proof adapts readily.  We write $\beta_0^M=\beta^{F'}$ and $\beta_1^M=\beta^F$.  More generally, one can define $\beta_i^M$ for all $i\in \mathbb N$ in an analogous way using \emph{minimal resolutions}.
For $a\in \mathbb{\R^2}$, $\beta_i^M(a)$ is called the \emph{$i^{\mathrm{th}}$ Betti number} of $M$ at $a$.
By \emph{Hilbert's syzygy theorem} \cite{peeva2011syzygies}, we have $\beta_i^M=0$ for $i>d$.

\subsection{Computation of Persistence Barcodes}\label{sec:ComputationPersistenceBarcodes}
We next discuss computation of the barcode $\B{M}$ of a persistence module $M$ indexed by a finite totally ordered set, given as input a presentation matrix $Q$ for $M$.  We  require that $Q$ is \emph{ordered}, meaning that the row and columns labels are both in increasing order. Typically, $\B{M}$ is computed by applying the \emph{standard reduction} \cite{edelsbrunner2010computational,edelsbrunner2002topological,zomorodian2005computing} to $Q$.  This a variant of Gaussian elimination which performs left-to-right column additions to construct a certain factorization of $Q$, from which $\B{M}$ can be read off directly.  

To elaborate, for any matrix $R$, let $R_{*j}$ denote the $j^{\mathrm{th}}$ column of $R$.  If $R_{*j}\ne 0$, let $\low_j$ denote the maximum index of a non-zero entry in $R_{*j}$.  We say $R$ is \emph{reduced} if $\low_j\ne \low_{j'}$ whenever $j\ne j'$ are the indices of non-zero columns in $R$.  An \emph{$RU$-decomposition} of an $r \times c$ matrix $A$ is a matrix factorization $A=RU$ such that $R$ is reduced and $U$ is upper-triangular.  The standard reduction computes an $RU$-decomposition of $A$ in time $O(rc \cdot \min(r,c))$.  

Given an $RU$-decomposition $Q=RU$ of our presentation matrix $Q$, we can read the barcode $\B{M}$ from $R$ via the following result, where we denote the label of row $i$ in $Q$ as $\rl(i)$ and the label of column $j$ as $\cl(j)$. 

\begin{proposition}[\cite{cohen2006vines}]\label{Prop:Barcode_from_Red_Mat}
    For any $RU$-decomposition of $Q$, $\B{M}$ is determined from $R$ as follows:
    \[ \B{M} = \{ [\rl(\low_j),\cl(j))  \mid R_{*j} \ne 0,\ \rl(\low_j)<\cl(j) \} \sqcup \{[\rl(i),\infty) \mid i \ne \low_j \textup{ for any $j$}\}. \]
\end{proposition}

Although $\B{M}$ can be obtained from $R$, without considering $U$, maintaining $U$ enables efficient updates to $\B{M}$ when the presentation matrix $Q$ is modified \cite{cohen2006vines}, as we explain next.  

\subsection{Updating an \texorpdfstring{$RU$}{RU}-Decomposition}\label{Sec:Updating_RU_Prelim}
Motivated by the problem of computing the persistent homology of each member of a parameterized family of data sets, Cohen-Steiner et al.~\cite{cohen2006vines} studied the problem of updating an $RU$-decomposition $Q=RU$ when the rows and columns of $Q$ are permuted.  We now briefly discuss the approach of \cite{cohen2006vines}, as well as closely related approaches to this problem appearing in subsequent work \cite{busaryev2010tracking,luo2024warmstarts,piekenbrock2024}.  In \cref{updateRU}, we apply $RU$-updates to compute barcode templates.

Given the $RU$-decomposition $Q=RU$ of an $r\times c$ matrix $Q$ and permutation matrices $P$ and $P'$ of respective dimensions $r\times r$ and $c\times c$, the aim is to compute an $RU$-decomposition of $\hat Q\coloneqq PQP'$ by updating $R$ and $U$.  Writing $\bar R\coloneqq P R P'$ and $\bar U \coloneqq (P')^{-1} U P'$, note that $\hat Q=\bar R\bar U$, but this is generally not an $RU$-decomposition of $\hat Q$, since $\bar R$ might not be reduced and $ \bar U$ might not be upper triangular.  We therefore perform matrix operations on $\bar R$ and $\bar U$ to obtain an $RU$-factorization $\hat Q=\hat R \hat U$.  These updates can be done using the following \emph{global strategy}, where we use the same notation for a matrix and the result after performing elementary operations on it:
\begin{enumerate}
    \item Do downward row operations on $\bar U$ to make it upper triangular, also doing the corresponding leftward column operations on $\bar R$.
    \item Then do rightward column operations on $\bar R$ to make it reduced, also doing the corresponding upward row operations on $\bar U$.
\end{enumerate}
Note that the row operations performed in step 2 preserve the upper-triangular structure of $\bar U$.  Hence, the resulting matrices $\hat R$ and $\hat U$ are reduced and upper triangular, respectively.  To see that $\hat Q=\hat R\hat U$, note that 
\begin{align*}
\hat R&= \bar R A^{-1} B\\
\hat U&= B^{-1} A \bar U,
\end{align*}
where left multiplication by $A$ performs the row operations on $\bar U$ in step 1, and right multiplication by $B$ performs the column operations on $\bar RA^{-1}$ in step 2.  Thus,  
\[\hat R\hat U=
(\bar R A^{-1} B)(B^{-1} A \bar U)
=\bar R\bar U=\hat Q,\]
as desired.  This update scheme adapts readily to the setting where $R$ and $U^{-1}$ are maintained instead of $R$ and $U$ \cite{busaryev2010tracking,luo2024warmstarts,piekenbrock2024}. 

Step 1 can be done by ordinary Gaussian elimination, while step 2 can be done using the standard reduction mentioned in \cref{sec:ComputationPersistenceBarcodes}.  However, if the permutations are known in advance to have special structure, then it may be advantageous to tailor the reduction algorithms to that structure.  Cohen-Steiner et al.~\cite{cohen2006vines} considered the case where the permutations are transpositions of adjacent indices, and showed that a specialized version of the above update scheme, called a \emph{vineyard update}, runs in time $O(r+c)$.  This scheme can be extended to handle arbitrary permutations, simply by decomposing each permutation into a sequence of transpositions of adjacent rows or columns, and then applying a separate vineyard update for each transposition in the sequence.  However, for general permutations, this can be quite expensive in practice, in part because quadratically many transpositions are required in the worst case.  Thus, it  can sometimes be more efficient in practice to perform a single global update to the $RU$-decomposition, as above, or to recompute the $RU$-decomposition from scratch.  Luo and Nelson \cite{luo2024warmstarts} consider such global updates in the setting where $R$ and $U^{-1}$ are maintained.   As an alternative approach, Piekenbrock and Perea \cite{piekenbrock2024} introduce a variant of the vineyard update scheme that decomposes an arbitrary permutation into a sequence of \emph{non-adjacent} transpositions.  They then perform a separate update for each transposition in the sequence by applying an algorithm of Busaryev et al.~\cite{busaryev2010tracking}.

Typically, the matrices arising in TDA applications are large and sparse, so they must be stored sparsely.  To implement vineyard updates, \cite{cohen2006vines} proposes to store $R$ in a column-sparse way and $U$ in a row-sparse way.  With these choices, it is inexpensive to permute columns of $R$ or rows of $U$, but it is expensive to permute rows of $R$ or columns of $U$.  Therefore, \cite{cohen2006vines}  permutes rows of $R$ implicitly, by maintaining a separate array representation of the row permutation.  Column permutations of $U$ are handled analogously.

\subsection{Line Arrangements}\label{sec:lineArrangements}
Line arrangements in (subregions of) $\R^2$ are classical objects in computational geometry.  In this paper, we will only consider line arrangements in the right half plane $\H \coloneqq [0,\infty)\times \R$.  Here, we define such arrangements and review standard results about computation, storage, and queries of line arrangements.

A \emph{cell} is a topological space homeomorphic to $\R^d$ for some $d\geq 0$.    
We define a \emph{cell complex} on $\H$ to be a decomposition of $\H$ into  finitely many cells, so that the topological boundary of each cell lies in the union of cells of lower dimension.  In such a cell complex, some cells are unbounded, and each cell has dimension at most two.

A \emph{line arrangement} in $\H$ is the cell complex on $\H$ induced by a finite set $W$ of lines in $\R^2$.  For each $i\in \{0,1,2\}$, the maximum number of $i$-cells in such an arrangement is $\Theta(|W|^2)$ \cite[Chapter 28]{toth2017handbook}.  We call the 0-cells, 1-cell, and 2-cells of the arrangement \emph{vertices}, \emph{edges}, and \emph{faces}, respectively.  The 1-skeleton of the arrangement (i.e., the union of vertices and edges) is \[\partial \H\cup \left(\H \cap \bigcup_{\ell\in W} \ell \right).\]   

The standard data structure for storing a line arrangement is the \emph{doubly connected edge list (DCEL)}.  This consists of a collection of vertices, edges, and faces, together with pointers specifying the facet relations.  The size of the DCEL is proportional to the number of cells in the arrangement.  The problem of computing a DCEL representation of a line arrangement, given the set of lines, has been extensively studied.  
The most popular algorithm for this, the \emph{Bentley-Ottmann algorithm} \cite{deberg2008computation,bentley1979algorithms},  constructs the DCEL representation of a line arrangement with $l$ lines and $v$ vertices in time $O((l + v) \log l)$ and memory $O(l+v)$.  
There exist asymptotically faster algorithms for constructing line arrangements, e.g., an algorithm of Chazelle and Edelsbrunner \cite{chazelle1992algorithm} runs in time $O(l \log l + v)$.  However, Bentley-Ottmann is a standard choice in software because it is simple and performs well in practice.

\subsubsection*{Point Location}
As mentioned in \cref{Sec:Augmented_Arrangement_Intro}, our scheme for efficient queries of fibered barcodes involves \emph{point location} in a line arrangement, i.e., finding the cell in the arrangement containing a given query point \cite[Section 38.3]{toth2017handbook}.  To perform many point location queries, it is usually most efficient to precompute a search data structure.  Let $e$ denote the number of edges in the arrangement.  On the theoretical side, there are several deterministic strategies which, in time $O(e\log e)$ and memory $O(e)$, compute a data structure on which we can perform a point location query in time $O(\log e)$ \cite{sarnak1986planar,edelsbrunner1986optimal,kirkpatrick1983optimal}.
Our main results, \cref{SimpleAugArrComplexity,Thm:Query_Cost}, depend on such a strategy.  

However, these deterministic strategies are rarely used in practice, and none are implemented in the popular  computational geometry library CGAL.  Instead, CGAL implements several other point location methods  \cite{cgal-arrangement-2}, any of which could be used to implement our query scheme for fibered barcodes; see \cite{haran2008experimental} for an experimental comparison of these methods on  planar arrangements.  One of these methods is a randomized incremental construction, based on a trapezoidal decomposition of the arrangement \cite{mulmuley1990fast,seidel1991simple,hemmer2012improved,deberg2008computation}. This method performs sufficiently well for practical use, with fast query times (see \cite{haran2008experimental} and \cite[Appendix B, full version]{hemmer2012improved}), and comes with the same guarantees as the deterministic constructions mentioned above, except that the bound $O(e\log e)$ on the time to construct the data structure is only in expectation \cite{hemmer2016optimal}.  

RIVET currently uses a simple point location strategy that requires no auxiliary search data structure.  (A similar strategy is implemented separately in CGAL.)  This strategy selects an initial point $p$ whose containing cell is known, and then walks linearly through the faces of the arrangement along the line segment from $p$ to the query point $q$.   By the \emph{zone theorem} for line arrangements \cite[Theorem  8.5]{deberg2008computation}, the worst-case runtime of a query is $\Theta(l)$ for this strategy.  Nevertheless, the strategy can be a practical choice for  a sequence $q_1,\ldots,q_k$ of query points in which consecutive points are close to each other, since we can choose $q_i$ as the initial point for the query $q_{i+1}$.  The sequence of queries arising in RIVET's interactive visualization has this structure.

\section{The Augmented Arrangement}\label{Sec:Augmented_Arrangement}
In this section, we define $\S(M)$, the \emph{augmented arrangement} of a finitely presented bipersistence module $M$.  
First, in \cref{sec:arr_def} we define the line arrangement $\cell(M)$ underlying $\S(M)$.   In \cref{Sec:CriticalLines}, we formally define the map $\push_\ell$ from the plane onto an affine line $\ell$ of non-negative slope, first mentioned in \cref{Sec:Augmented_Arrangement_Intro}.  Using this, in \cref{Sec:Discrete_Barcodes} we define the barcode template $\P^\ce$ stored at each face $\ce$ of $\S(M)$.

\subsection{Definition of the Line Arrangement}\label{sec:arr_def}
Let $\L$ denote the set of all affine lines in $\R^2$.  For $\mathcal I\subset \R\cup \{\infty\}$, we let $\L_{\mathcal I}$ denote the subset of $\L$ consisting of lines whose slopes lie in $\mathcal I$.  For example, $\L_{[0,\infty]}$ is the set of lines with non-negative, possibly infinite slope.  Note that each $\ell\in \L_{[0,\infty]}$ is a totally ordered subposet of $\R^2$.

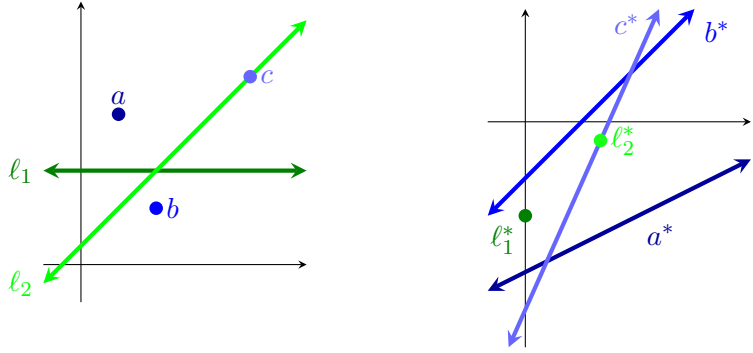
\begin{figure}[ht]
  \begin{center}
    \begin{tikzpicture}
        \draw[->] (-0.5,0) -- (3,0);
        \draw[->] (0,-0.5) -- (0,3.5);
        
        \draw[<->,ultra thick,color=green!50!black] (-0.5,1.25) -- (3,1.25);
        \node[left,color=green!50!black] at (-0.5,1.25) {$\ell_1$};
        \draw[<->,ultra thick,color=green] (-0.5,-0.25) -- (3,3.25);
        \node[left,color=green] at (-0.5,-0.25) {$\ell_2$};
          
        \coordinate (p1) at (0.5,2);
        \fill[blue!60!black] (p1) circle (0.09);
        \node[above,color=blue!60!black] at (p1) {$a$};
      
        \coordinate (p2) at (1,0.75);
        \fill[blue] (p2) circle (0.09);
        \node[right,color=blue] at (p2) {$b$};
        
        \coordinate (p3) at (2.25,2.5);
        \fill[blue!60!white] (p3) circle (0.09);
        \node[right,color=blue!60!white] at (p3) {$c$};
        
        \node at (0,-1) {};
      \end{tikzpicture}
      \hspace{.8in}
      \begin{tikzpicture}
        \draw[->] (-0.5,0) -- (3,0);
        \draw[->] (0,-3) -- (0,1.5);
        
        \draw[<->,ultra thick,color=blue!60!black] (-0.5,-2.25) -- (3,-0.5);
        \node[color=blue!60!black] at (1.8,-1.5) {$a^*$};
        \draw[<->,ultra thick,color=blue] (-0.5,-1.25) -- (2.25,1.5);
        \node[color=blue] at (2.55,1.2) {$b^*$};
        \draw[<->,ultra thick,color=blue!60!white] (-0.222,-3) -- (1.778,1.5);
        \node[color=blue!60!white] at (1.35,1.3) {$c^*$};
        
        \coordinate (q1) at (0,-1.25);
        \fill[green!50!black] (q1) circle (0.09);
        \node[below left,color=green!50!black] at (q1) {$\ell_1^*$};
        
        \coordinate (q2) at (1,-0.25);
        \fill[green!90] (q2) circle (0.09);
        \node[right,color=green] at (q2) {$\ell_2^*$};
    \end{tikzpicture}
  \end{center}
  \caption{Illustration of point-line duality}
  \label{fig:duality}
\end{figure}

The following standard point-line duality gives a parameterization of $\L_{\R}$ by $\R^2$: 
\begin{align*}
	&(-)^* \colon \L_{\R} \to \R^2 & 	&(-)^* \colon \R^2 \to \L_{\R} \\
	&y=qx+r\ \mapsto\ (q,-r) &			&(q,r)\ \mapsto\ y=qx-r.
\end{align*}
These maps are  homeomorphisms with respect to the standard topology on $\R^2$ and the restriction to $\L_{\R}$ of the usual Grassmannian topology.  It is easily checked that this duality preserves incidence, in the sense that for $a\in \R^2$ and $\ell\in \L_{\R}$, $a\in \ell$ if and only if $\ell^*\in a^*$; see \cref{fig:duality}.  This point-line duality does not extend naturally to vertical lines.

Define $\XiSp\subset \R^2$ by $\XiSp=\supp \beta_0^M\cup\supp \beta_1^M$, where $\supp$ denotes the support.  For $X$ a finite subset of $\R^2$, let $\lub X$ denote the \emph{join of $X$}, i.e., the least upper bound of $X$ in $\R^2$.
For example, $\lub\{(2,0),(1,1),(0,2)\}=(2,2)$.  For $a,b \in \R^2$, we often write $\lub\{a,b\}$ as $a\lub b$.

\begin{definition}
\mbox{}
  \begin{itemize}
    \item[(i)] We call a pair of distinct elements $a,b\in \R^2$ \emph{weakly incomparable} if either $a$ and $b$ are incomparable in the partial order on $\R^2$, or $a$ and $b$ share a coordinate (i.e., $a_1=b_1$ or $a_2=b_2$). 
    \item[(ii)] We call $\alpha\in \R^2$ an \emph{anchor} if $\alpha=s \lub t$ for some weakly incomparable pair $s,t\in \XiSp$. 
  \end{itemize}
\end{definition}

\begin{definition}
We define $\cell(M)$ to be the line arrangement in $\H$ induced by the set of lines
\[ \{\alpha^*\mid \alpha\in \R^2 \textup{ is an anchor}\}. \] 
\end{definition}

\begin{remark}
There are two types of vertices in the arrangement $\cell(M)$:
the intersection of the duals of two anchors, and the intersection of the dual of an anchor with $\partial \H$.  These two types are not mutually exclusive, since the duals of two anchors with the same $y$-coordinate intersect on the line $\partial \H$.  To better understand the first type of vertex, note that for two anchors $\alpha\ne \gamma$, $\alpha^*$ and $\gamma^*$ intersect in $\H$ if and only if there exists some $\ell\in \Lplus$ containing both $\alpha$ and $\gamma$; such $\ell$ exists precisely when $\alpha$ and $\gamma$ are comparable and have distinct $x$-coordinates.  
\end{remark}

As in \cref{Sec:Main_Results_Intro}, let $\kappa =\kappa_x \kappa_y$, for $\kappa_x$ and $\kappa_y$ the number of unique $x$ and $y$ coordinates, respectively, of points in $\XiSp$.  Clearly, the number of anchors is at most $\kappa$.  Hence the number of lines in $\cell(M)$, including $\partial \H$, is at most $\kappa+1$.  
The size bounds mentioned in \cref{sec:lineArrangements} thus imply that the numbers of vertices, edges, and faces in $\cell(M)$ are each  $O(\kappa^2)$.  

\subsection{Push Maps}\label{Sec:CriticalLines}
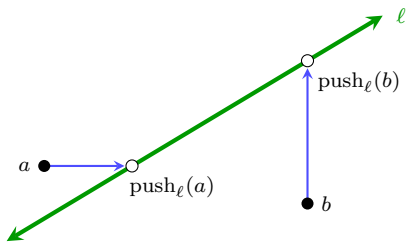
\begin{figure}[ht]
  \begin{center}
    \begin{tikzpicture}
      \tikzstyle{every node}=[font=\scriptsize]
      \draw[<->,ultra thick,color=green!60!black] (0,0) -- (5,3);
      \node[right,color=green!60!black] at (5.05,3) {$\ell$};   
      
      \coordinate (u) at (0.5,1);
      \coordinate(pu) at (1.667,1);
      \draw[->,thick,blue!70!white] (u) -- (1.567,1);      
      \fill (u) circle (0.08);
      \draw[fill=white] (pu) circle (0.08);
      \node[left] at (0.45,1) {$a$};
      \node[below right] at (1.5,1) {$\push_\ell(a)$};
      
      \coordinate (v) at (4,0.5);
      \coordinate(pv) at (4,2.4);
      \draw[->,thick,blue!70!white] (v) -- (4,2.3);      
      \fill (v) circle (0.08);
      \draw[fill=white] (pv) circle (0.08);
      \node[right] at (4.05,0.5) {$b$};
      \node[below right] at (pv) {$\push_\ell(b)$};
    \end{tikzpicture}
  \end{center}
  \caption{Illustration of $\push_{\ell}$ for a line $\ell$ of positive, finite slope.}
  \label{fig:push_map}
\end{figure}

Recall that for each $\ell\in \L_{[0,\infty]}$, the partial order on $\R^2$ restricts to a total order on $\ell$.  
This extends to a total order on $\ell \cup \{\infty\}$ by taking $\infty$ to be the maximum element.  
\begin{definition}
    For $\ell\in \L_{[0,\infty]}$, define the \emph{push map} $\push_\ell:\R^2 \to \ell \cup \{\infty\}$ by taking 
    \[\push_\ell(a)=\min\{b\in \ell\mid a\leq b\},\]
    where it is understood that the minimum over an empty set is $\infty$.
\end{definition}
Note that
\begin{itemize} 
  \item $\infty\in \im \push_\ell$ if and only if $\ell$ is either horizontal or vertical.
  \item For $a\in \R^2$ and $\push_\ell(a)=c\in \ell$, either $a_1=c_1$ or $a_2=c_2$; see \cref{fig:push_map}.
  \item $\push_\ell$ preserves the partial order on $\R^2$, i.e., $\push_\ell(a)\leq \push_\ell(b)$ for all $a \leq b \in \R^2$.
\end{itemize}

For any $a\in \R^2$, the maps $\{\push_\ell\}_{\ell\in \L_{(0,\infty)}}$ induce a map \[\push_{a}:\L_{(0,\infty)} \to \R^2,\] defined by $\push_{a}(\ell)=\push_\ell(a)$.   

\begin{lemma}\label{Lem:ContinuityOfPush}
For any $a\in \R^2$, the map $\push_{a}$ is continuous.  
\end{lemma}
\begin{proof}
  This follows from the observation that for any  $\ell\in \L_{(0,\infty)}$, $\push_{a}(\ell)$ is the unique intersection of $\ell$ with \[\{(a_1,y)\mid y\geq a_2\} \cup \{(x,a_2)\mid x\geq a_1\}.\qedhere\]  
\end{proof}

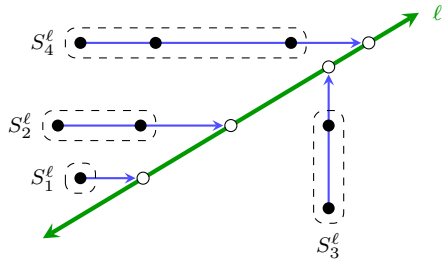
\begin{figure}[ht]
  \begin{center}
    \begin{tikzpicture}
      \tikzstyle{every node}=[font=\scriptsize]

      \draw[<->,ultra thick,color=green!60!black] (0,0) -- (5,3);
      \node[right,color=green!60!black] at (5.05,3) {$\ell$};   
      
      \coordinate (u) at (0.5,0.8);
      \coordinate (pu) at (1.333,0.8);
      \draw[->,thick,blue!70!white] (u) -- (1.233,0.8);
      \fill (u) circle (0.08);
      \draw[fill=white] (pu) circle (0.08);
      \draw[dashed,rounded corners=4pt] (0.3,0.6) rectangle (0.7,1);
      \node[left] at (0.3,0.8) {$\XiSp^\ell_1$};
      
      \coordinate (v1) at (0.2,1.5);
      \coordinate (v2) at (1.3,1.5);
      \coordinate (pv) at (2.5,1.5);
      \draw[->,thick,blue!70!white] (v1) -- (2.4,1.5);
      \fill (v1) circle (0.08);
      \fill (v2) circle (0.08);
      \draw[fill=white] (pv) circle (0.08);
      \draw[dashed,rounded corners=4pt] (0,1.3) rectangle (1.5,1.7);
      \node[left] at (0,1.5) {$\XiSp^\ell_2$};
      
      \coordinate (x1) at (3.8,0.4);
      \coordinate (x2) at (3.8,1.5);
      \coordinate (px) at (3.8,2.28);
      \draw[->,thick,blue!70!white] (x1) -- (3.8,2.18);
      \fill (x1) circle (0.08);
      \fill (x2) circle (0.08);
      \draw[fill=white] (px) circle (0.08);
      \draw[dashed,rounded corners=4pt] (3.6,0.2) rectangle (4,1.7);
      \node[below] at (3.8,0.2) {$\XiSp^\ell_3$};
      
      \coordinate (y1) at (0.5,2.6);
      \coordinate (y2) at (1.5,2.6);
      \coordinate (y3) at (3.3,2.6);
      \coordinate (py) at (4.333,2.6);
      \draw[->,thick,blue!70!white] (y1) -- (4.233,2.6);
      \fill (y1) circle (0.08);
      \fill (y2) circle (0.08);
      \fill (y3) circle (0.08);
      \draw[fill=white] (py) circle (0.08);
      \draw[dashed,rounded corners=4pt] (0.3,2.4) rectangle (3.5,2.8);
      \node[left] at (0.3,2.6) {$\XiSp^\ell_4$};
    \end{tikzpicture}
  \end{center}
  \caption{Illustration of a totally ordered partition $\XiSp^\ell$ of $\XiSp$.  The elements of $S$ are drawn as black dots.}
  \label{fig:partition}
\end{figure}

For $\ell\in \L_{(0,\infty)}$, $\push_\ell$ induces a totally ordered partition $\XiSp^\ell$ of $\XiSp$, where elements of the partition are level sets of the restriction of $\push_\ell$ to $\XiSp$, and the total order on $\XiSp^\ell$ is the pullback of the total order on $\ell$.  We denote the $i^{\mathrm{th}}$ element of $\XiSp^\ell$ as $\XiSp^\ell_i$.
See  \cref{fig:partition} for an illustration.

The following result is the key to defining the augmented arrangement.

\begin{theorem}\label{Thm:2Cells}
If the duals of $\ell, \ell'\in \L_{(0,\infty)}$ are contained in the same (open) face of $\cell(M)$, then $\XiSp^\ell=\XiSp^{\ell'}$.
\end{theorem}
\begin{proof}
Since point-line duality preserves incidence, the 1-skeleton of $\cell(M)$ is 
\begin{equation}\label{Eq:Description_Of_One_Skel}
\{\ell^*\mid \ell \in \L_{(0,\infty)} \textup{ contains an anchor} \}\cup \partial \H, 
\end{equation}
Moreover, $\partial \H=\{\ell^*\mid \ell \in \L_{\{0\}}\}$.  
Thus, since a face of $\cell(M)$ is connected, it suffices to show that for any $\ell\in \L_{(0,\infty)}$ that contains no anchor, there is an open neighborhood $\mathcal{N} \subset \L_{(0,\infty)}$ of $\ell$ such that $\XiSp^{\ell}=\XiSp^{\ell'}$ for all $\ell'\in \mathcal{N}$.
To show this, assume that $\ell \in \L_{(0,\infty)}$ contains no anchor, and consider distinct $s,t\in \XiSp$ with $\push_\ell(s)=\push_\ell(t)$.  
Note that $s,t$ must lie either on the same horizontal line or the same vertical line; otherwise $s$ and $t$ would be incomparable and we would have $\push_\ell(s)=\push_\ell(t)=s\lub t\in \ell$, contradicting the assumption that $\ell$ contains no anchor.  
Assume without loss of generality that $s<t$ and $s,t$ lie on the same horizontal line $h$.  
Then we must also have that $\push_\ell(s)=\push_\ell(t)$ lies on $h$. 
Since $t=s\lub t$, $t$ is an anchor. However, $\ell$ does not contain any anchor, so we must have $t<\push_\ell(t)$.  Any sufficiently small perturbation $\ell'$ of $\ell$ will also intersect $h$ at a point $a$ to the right of $t$, so that we have $\push_{\ell'}(s) = \push_{\ell'}(t) = a$.  
Thus for all $\ell'$ in a neighborhood of $\ell$, $\push_{\ell'}(s)=\push_{\ell'}(t)$.  

In fact, since $\XiSp$ is finite, we can choose a single neighborhood $\mathcal{N}$ of $\ell$ in $\L_{(0,\infty)}$ such that for any $s,t\in \XiSp$ with $\push_\ell(s)=\push_\ell(t)$, we have $\push_{\ell'}(s)=\push_{\ell'}(t)$ for all $\ell'\in \mathcal{N}$.  
In view of \cref{Lem:ContinuityOfPush}, by choosing $\mathcal{N}$ to be smaller if necessary, we may also ensure that if $\push_\ell(s)\ne \push_\ell(t)$, then $\push_{\ell'}(s)\ne \push_{\ell'}(t)$ for all $\ell'\in \mathcal{N}$.  
Thus, the partition $\XiSp^{\ell'}$ is independent of the choice of $\ell'\in \mathcal{N}$.  
Finally, \cref{Lem:ContinuityOfPush} also implies that the total order on $\XiSp^{\ell'}$ is independent of the choice of $\ell' \in\mathcal{N}$. 
\end{proof}

\begin{remark}
  In fact, \cref{Thm:2Cells} can be strengthened to show that for any $\ell,\ell'\in \L_{(0,\infty)}$, the duals of $\ell$ and $\ell'$ lie in the same cell of $\cell(M)$ if and only if $\XiSp^\ell=\XiSp^{\ell'}$.  However, we will not need the full generality of this stronger result.  \end{remark}

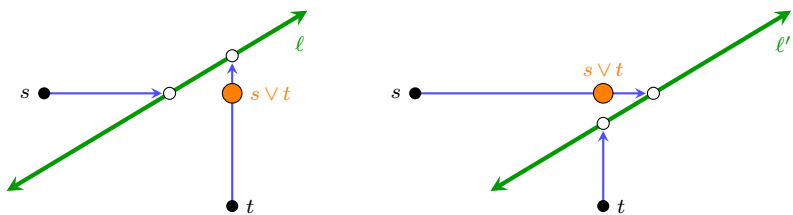
\begin{figure}[ht]
  \begin{center}
    \begin{tikzpicture}[scale=1]
      \tikzstyle{every node}=[font=\scriptsize]
      \draw[<->,ultra thick,color=green!60!black] (1,-.3) -- (5.0,2.1);
      \node[below,color=green!60!black] at (4.9,1.9) {$\ell$};   
      
      \coordinate (u) at (1.5,1);
      \coordinate(pu) at (3.167,1);
      \draw[->,thick,blue!70!white] (u) -- (3.067,1);      
      \fill (u) circle (0.08);
      \draw[fill=white] (pu) circle (0.08);
      \node[left] at (1.45,1) {$s$};
      
      \coordinate (v) at (4,-.5);
      \coordinate(pv) at (4,1.5);
      \draw[->,thick,blue!70!white] (v) -- (4,1.4);      
      \fill (v) circle (0.08);
      \draw[fill=white] (pv) circle (0.08);
      \node[right] at (4.05,-0.5) {$t$};

      \coordinate(LUB) at (4,1);
      \draw[fill=orange] (LUB) circle (0.13);
      \node[right,color=orange] at (4.1,1) {$s \lub t$};
    \end{tikzpicture}
    \hspace{.25in}
    \begin{tikzpicture}[scale=1]
      \tikzstyle{every node}=[font=\scriptsize]
      \draw[<->,ultra thick,color=green!60!black] (2.5,-.3) -- (6.5,2.1);
      \node[below,color=green!60!black] at (6.4,1.9) {$\ell'$};   
      
      \coordinate (u) at (1.5,1);
      \coordinate(pu) at (4.667,1);
      \draw[->,thick,blue!70!white] (u) -- (4.567,1);      
      \fill (u) circle (0.08);
      \draw[fill=white] (pu) circle (0.08);
      \node[left] at (1.45,1) {$s$};

      \coordinate (v) at (4,-.5);
      \coordinate(pv) at (4,0.6);
      \draw[->,thick,blue!70!white] (v) -- (4,0.5);      
      \fill (v) circle (0.08);
      \draw[fill=white] (pv) circle (0.08);
      \node[right] at (4.05,-0.5) {$t$};

      \coordinate(LUB) at (4,1);
      \draw[fill=orange] (LUB) circle (0.13);
      \node[above,color=orange] at (4,1.1) {$s \lub t$};
    \end{tikzpicture}
  \end{center}
  \caption{An illustration of how incomparable points $s,t \in S$ push onto lines $\ell$ and $\ell'$ lying above and below the anchor  $s \lub t$. The dual points $\ell^*$ and $\ell'^*$ lie in neighboring $2$-cells of $\cell(M)$ with shared boundary lying on $(s \lub t)^*$. 
  Note that $\push_{\ell}(s)<\push_{\ell}(t)$ and $\push_{\ell'}(t)<\push_{\ell'}(s)$, i.e.,  the order in which the points push onto the line switches as the line changes.} 
  \label{fig:Criticality_Illustration}
\end{figure}

To derive our main algorithm for computing the augmented arrangement, in \cref{sec:updateLiftMap} we give a detailed analysis of the relationship between  $S^{\ell}$ and $S^{\ell'}$ for lines $\ell$ and $\ell'$ lying in neighboring faces of $\cell(M)$.  A simple such case is illustrated in \cref{fig:Criticality_Illustration}; see also \cref{fig:line_crosses_anchor} for further examples.     

\subsection{Barcode Templates}\label{Sec:Discrete_Barcodes}
Fixing a face $\ce$ of $\cell(M)$, we now define the barcode template $\P^\ce$ stored at the face $\ce$ in the augmented arrangement $\S(M)$. 
    
\begin{definition}
\cref{Thm:2Cells} yields a totally ordered partition \[\XiSp^\ce=\{S^\ce_1,\ldots, S^\ce_k\}\] of $\XiSp$, given by $S^{\ce}=S^{\ell}$ for any line $\ell$ with $\ell^*\in \ce$.  For each $S^\ce_i$, we let 
 \[ \PCal^\ce_i= \bigvee \left(\bigcup_{j\leq i} S^\ce_j \right).\]
 The set of \emph{template points} at $\ce$ is $\{\PCal^\ce_1, \ldots, \PCal^\ce_k\}$.
\end{definition}
   
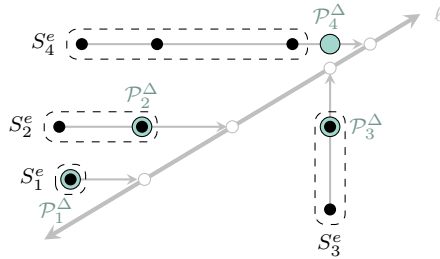
\begin{figure}[ht]
  \begin{center}
    \begin{tikzpicture}[scale=1]
      \tikzstyle{every node}=[font=\scriptsize]

      \draw[<->,ultra thick,color=black!25] (0,0) -- (5,3);
      \node[right,color=black!25] at (5.05,3) {$\ell$};   

      \coordinate (u) at (0.35,0.8);
      \coordinate (pu) at (1.333,0.8);
      \draw[->,thick,black!25] (u) -- (1.24,0.8);
      \draw[fill=\tpcolor] (u) circle (0.13);
      \fill (u) circle (0.08);
      \draw[color=black!25,fill=white] (pu) circle (0.08);
      \draw[dashed,rounded corners=4pt] (0.15,0.6) rectangle (0.55,1);
      \node[left] at (0.15,0.8) {$\XiSp^e_1$};
      \node[\tpcolor!70!black] at (0.16,0.38) {$\PCal^\ce_1$};
      
      \coordinate (v1) at (0.2,1.5);
      \coordinate (v2) at (1.3,1.5);
      \coordinate (pv) at (2.5,1.5);
      \draw[->,thick,black!25] (v1) -- (2.43,1.5);
      \fill (v1) circle (0.08);
      \draw[fill=\tpcolor] (v2) circle (0.13);
      \fill (v2) circle (0.08);
      \draw[color=black!25,fill=white] (pv) circle (0.08);
      \draw[dashed,rounded corners=4pt] (0,1.3) rectangle (1.5,1.7);
      \node[left] at (0,1.5) {$\XiSp^e_2$};
      \node[\tpcolor!70!black] at (1.3,1.92) {$\PCal^\ce_2$};

      \coordinate (x1) at (3.8,0.4);
      \coordinate (x2) at (3.8,1.5);
      \coordinate (px) at (3.8,2.28);
      \draw[->,thick,black!25] (x1) -- (3.8,2.21);
      \fill (x1) circle (0.08);
      \draw[fill=\tpcolor] (x2) circle (0.13);
      \fill (x2) circle (0.08);
      \draw[color=black!25,fill=white] (px) circle (0.08);
      \draw[dashed,rounded corners=4pt] (3.6,0.2) rectangle (4,1.7);
      \node[below] at (3.8,0.2) {$\XiSp^e_3$};
      \node[\tpcolor!70!black] at (4.29,1.5) {$\PCal^\ce_3$};

      \coordinate (y1) at (0.5,2.6);
      \coordinate (y2) at (1.5,2.6);
      \coordinate (y3) at (3.3,2.6);
      \coordinate (py) at (4.333,2.6);
      \draw[->,thick,black!25] (y1) -- (4.24,2.6);
      \fill (y1) circle (0.08);
      \fill (y2) circle (0.08);
      \fill (y3) circle (0.08);
      \draw[fill=\tpcolor] (3.8,2.6) circle (0.13);
      \draw[color=black!25,fill=white] (py) circle (0.08);
      \draw[dashed,rounded corners=4pt] (0.3,2.4) rectangle (3.5,2.8);
      \node[left] at (0.3,2.6) {$\XiSp^e_4$};
      \node[\tpcolor!70!black] at (3.8,2.96) {$\PCal^\ce_4$};
    \end{tikzpicture}
  \end{center}
  \caption{Illustration of template points $\PCal^\ce$ (light blue), for  $S$ and $\ell$ as in \cref{fig:partition} and $\ce$ the face containing $\ell^*$.}
  \label{fig:Template_Points}
\end{figure}

In words, elements of $\PCal^\ce$ are joins of unions of prefixes of $S^\ce$.  See \cref{fig:Template_Points} for an illustration.  
Note that $\PCal^\ce$ is a subset of $\R^2$  and that, by construction, the restriction of the partial order on $\R^2$ to $\PCal^\ce$ is a total order.  We let $M^\ce$ denote the restriction of $M$ to $\PCal^\ce$.  Since $M$ is assumed to be finitely presented, 
 $M^\ce$ is \pfd, so it has a well-defined barcode $\B{M^\ce}$.  

\begin{definition}
    We call $\P^\ce\coloneqq\B{M^\ce}$ the \emph{barcode template} at $\ce$.  
\end{definition}
$\B{M^\ce}$ consists of intervals $[a,b)\subset \PCal^\ce$ with $a<b\in \PCal^\ce\cup\{\infty\}$.  Identifying each interval $[a,b)$  with the pair $(a,b)$, we recover the form of the barcode template specified in \cref{Sec:Augmented_Arrangement_Intro}.  

\begin{definition}\label{Def:Aug_Arrangement}
    The \emph{augmented arrangement} $\S(M)$ consists of the line arrangement $\cell(M)$ and a data structure for log-time point location queries in $\cell(M)$ as in \cref{sec:lineArrangements}, together with the barcode templates $\P^\ce$ at all faces $\ce$ of $\cell(M)$.
\end{definition}

\begin{remark}
The fibered barcode $\fb$ determines neither the Betti numbers of $M$ nor the resulting set of anchors; see \cite[Example 2.2]{lesnick2015interactive}.  Hence $\fb$ does not determine $\S(M)$.  We expect that by defining  the set of anchors differently, one can obtain a variant of $\S(M)$ which depends only on $\mathcal F(M)$.  We leave the exploration of this to future work.
\end{remark}

Having defined the augmented arrangement, it is now easy to see that it satisfies the size bound of \cref{SimpleAugArrComplexity}\,(i): 

\begin{proof}[Proof of \cref{SimpleAugArrComplexity}\,(i)]
As noted in \cref{sec:arr_def}, $\cell(M)$ has size $O(\kappa^2)$.  The point location data structure has size linear in the number of edges in $\cell(M)$, so also has size $O(\kappa^2)$.  Since $\S(M)$ stores a barcode template of size $O(m)$ for each face of $\cell(M)$, the total size of all barcode templates is $O(m\kappa^2)$.  Thus, the total size of $\S(M)$ is $O(m\kappa^2)$.  \end{proof}

\section{Querying the Augmented Arrangement}\label{Sec:QueryingMath}
We now specify the procedure to query $\S(M)$ for the barcode $\B{M^\ell}$ of a given line $\ell\in \L_{[0,\infty]}$.  The main result of this section, \cref{Thm:QueriesMain}, establishes the correctness of this  procedure.  In \cref{Sec:Representing_And_Querying_Aug_Arrangement}, we discuss the computational details of our queries and bound the complexity of a query, proving \cref{Thm:Query_Cost}. 

\subsection{The Query Procedure}\label{Sec:Query_Procedure}
The query of $\S(M)$ for the barcode $\B{M^\ell}$ involves two main steps.  First, we choose a face $\ce$ in $\cell(M)$, as specified below.  Second, we obtain the intervals of $\B{M^\ell}$ from the pairs of $\P^\ce$ by pushing the points in each pair $(a,b)\in \P^\ce$ onto the line $\ell$, via the map $\push_\ell$ of \cref{Sec:CriticalLines}. 

\begin{figure}[ht]
  \begin{center}
    \begin{tikzpicture}[scale=1.5]
      \draw[<->] (0,3) -- (0,0) -- (4,0);

      \draw[<->,ultra thick,color=blue] (-0.3,1) -- (4,1);
      \draw[<->,dashed,ultra thick,color=blue,opacity=.5] (-0.3,1.2) -- (4,1.5);  
      
      \node[below right,color=blue] at (0,1) {$\ell_2$};
      \node[above right,color=blue,opacity=.5] at (0,1.25) {$\ell'_2$};

      \draw[<->,ultra thick,color=purple] (1,-0.3) -- (1,3);
      \draw[<->,dashed,ultra thick,color=purple,opacity=.5] (1.2,-0.3) -- (1.4,3); 
      \node[below left,color=purple] at (1,3) {$\ell_3$};
      \node[below right,color=purple,opacity=.5] at (1.4,3) {$\ell'_3$};

     \foreach \p in {(1,1),(1,2),(2.4,2)}
        { \draw[fill=orange] \p circle (0.07); } 

      \draw[<->,ultra thick,color=green!70!black] (0,-0.3) -- (4,2.2);  
      \node[color=green!70!black,right] at (3.15,2.1) {$\ell_1$};
    \end{tikzpicture}
    \hspace{.8in}
    \begin{tikzpicture}[scale=1.2]
      \draw[<->,ultra thick] (0,-2.5) -- (0,1.5);
      \draw[->,ultra thick] (0,-2) -- (3,-0.5);   
      \draw[->,ultra thick] (0,-2) -- (3,1.6);    
      \draw[->,ultra thick] (0,-1) -- (3,0.5);    

      \fill[blue] (0,-1) circle (0.09);
      \fill[blue,opacity=.5] (0.14,-1.22) circle (0.09);
      \node[left,color=blue] at (-0.05,-.8) {$\ell_2^*$};
      \node[left,color=blue,opacity=.5] at (0.08,-1.4) {${\ell'_2}^*$};
      \node[color=blue] at (0.55,-.95) {$\ce_2$};

      \fill[green!70!black] (1.25,0.3) circle (0.09);
      \node[color=green!70!black] at (1,0.3) {$\ell_1^*$};
      \node[color=green!70!black] at (1.2,1.0) {$\ce_1$};

      \node[color=purple] at (2.5,0.6) {$\ce_3$};
    \end{tikzpicture}
  \end{center}
   \caption{Three anchors are drawn at left in orange; the corresponding line arrangement $\cell(M)$ is shown at right.  For each line $\ell_i$ at left, the corresponding face $\ce_i$, chosen as in \cref{Sec:Query_Procedure}, is labeled at right in the same color.  For the horizontal line $\ell_2$ and vertical line $\ell_3$, choices of the perturbations $\ell'_2$ and $\ell'_3$ are also shown at left as dashed lines.  The duals $\ell_1^*$, $\ell_2^*$, and ${\ell'_2}^*$ are also shown at right.  The dual point ${\ell'_3}^*$ is not shown because this point has very large coordinates, since $\ell'_3$ is nearly vertical.}
  \label{fig:coface}
\end{figure}
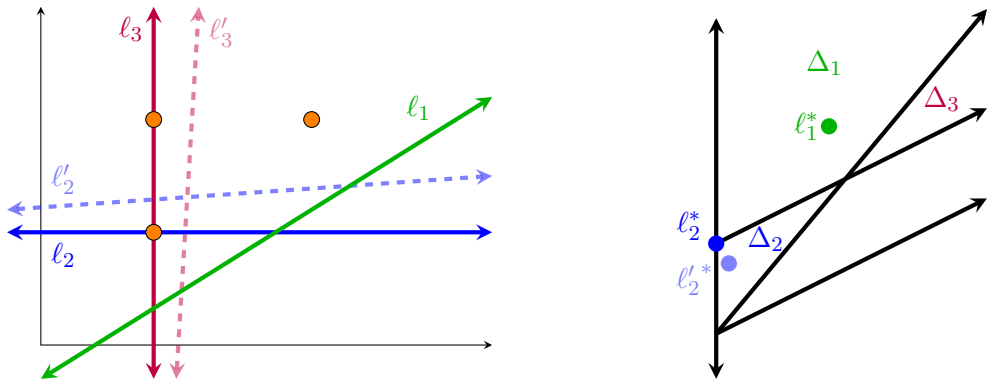

If $\ell$ is not vertical, then we choose $\ce$ to be a 2-dimensional coface of $\ell^*$ (meaning that $\ell^*$ is in the closure of $\ce$); if $\ell$ contains no anchor, then there is only one such coface.  In more detail, we choose $\ce$ as follows (see \cref{fig:coface} for an illustration):
\begin{itemize}
  \item If $\ell$ is neither horizontal nor vertical, then we choose $\ce$ to be any 2-D coface of $\ell^*$. 
  \item If $\ell$ is horizontal, say $\ell$ is the line $y=c$, then we take $\ce$ to be the face whose points are dual to lines $\ell'$ of positive slope lying above all anchors with $y$-coordinate at most $c$ and lying below all anchors with $y$-coordinate greater than $c$.  (It is helpful to think of the line $\ell'$ as arising as a certain small perturbation of $\ell$, as in \cref{fig:coface}.)
      
  \item Symmetrically, if $\ell$ is vertical, say $\ell$ is the line $x=c$, then we take $\ce$ to be the face whose points are dual to lines $\ell'$ of positive slope lying right of all anchors with $x$-coordinate at most $c$ and lying left of all anchors with $x$-coordinate greater than $c$.  
\end{itemize}

\begin{remark}\label{Rem:Concrete_Description_of_e}
    For $\ell$ horizontal or vertical, the face $\ce$ can be described directly in terms of $\cell(M)$, as follows: If $\ell$ is the horizontal line $y=c$, then since $\ell$ has slope 0, $\ell^*\in \partial \H$.  If $\ell^*$ is contained in an edge of $\cell(M)$, then $\ce$ is the unique coface of this edge.  Otherwise, $\ell^*$ is a vertex of $\cell(M)$.  Then the 2-D cofaces of $\ell^*$ are ordered vertically in $\H$ and  $\ce$ is the bottom such coface.
    
    If $\ell$ is the vertical line $x=c$, let $k$ be the line in the arrangement $\cell(M)$ of maximum slope, among those having slope at most  $c$, if a unique such line exists; if there are several such lines, then take $k$ to be the one with the largest $y$-intercept.  Then $k$ contains a unique unbounded edge in $\cell(M)$ and $\ce$ is the face lying directly above this edge.  (To see why this choice of $\ce$ is correct, note that points in $\ce$ are dual to lines of large slope which lie below $k^*$.)
        If a line $k$ as above does not exist, then $\ce$ is the bottom unbounded face of $\cell(M)$; since $\cell(M)$ contains no vertical lines besides $\partial \H$, this face is uniquely defined.
\end{remark}

\subsection{Correctness of the Query Procedure}\label{Sec:Query_Theorem}

\begin{theorem}\label{Thm:QueriesMain}
  For any line $\ell\in \L_{[0,\infty]}$ and a  face $\ce$ chosen as above, we have
  \[\B{M^\ell}=\{\,[\push_\ell(a),\,\push_\ell(b)) \mid [a,b)\in \P^\ce,\ \push_\ell(a)<\push_\ell(b)\,\},\]
  where we take $\push_\ell(\infty)=\infty$.
\end{theorem}

Note that if $\ell^*$ lies in a face of $\cell(M)$, then that face is necessarily $\ce$, and we have  $\push_\ell(a)<\push_\ell(b)$ for all $(a,b)\in \P^\ce.$  Thus, the theorem statement simplifies for such $\ell$.  However, if  $\ell^*$ is contained in the 1-skeleton of $\cell(M)$, then it is possible to have $\push_\ell(a)=\push_\ell(b)$.   

We prepare for the proof of \cref{Thm:QueriesMain} with three lemmas.  For $a\in \R^2$, let 
\[ \I(a)=\{s\in \XiSp \mid s\leq a\}.\]

\begin{lemma}\label{Lem:GradesOfInfluenceAndIsomorphisms}
For $a\leq b\in \R^2$ with $\I(a)=\I(b)$, $M_{a,b}$ is an isomorphism.  
\end{lemma}

\begin{proof}
Let $\eta:F\to F'$ be a minimal presentation for $M$.  It is easy to check that $F_{a,b}$ and $F'_{a,b}$ are both isomorphisms.  As we have \[F'_{a,b}\circ\eta_{a}=\eta_b\circ F_{a,b},\] $F'_{a,b}$ descends to an isomorphism on $\coker(\eta)\cong M$. 
\end{proof}

Given totally ordered sets $\T$ and $\U$ with $\T$ finite, a functor (i.e., order-preserving function) $\G \colon \T \to \U$, and a functor $N \colon \T \to \kvect$, we let $\Lan_{\G}(N) \colon \U \to \kvect$ denote the left Kan extension of $N$ along $\G$ \cite[Section 6.1]{riehl2017category}.  Then $\Lan_{\G}(N)$ admits the following simple, elementary description:  For $u\in \U$, let 
\[\G_{\leq u}=\{t\in \T\mid \G(t)\leq u\}.\]
For $\G_{\leq u}\ne \emptyset$, we define
\[ \G^{-1}(u)=\max\,\G_{\leq u}.\]
Then $\Lan_{\G}(N) \colon \U \to \kvect$ is given by
\begin{align*}
    \Lan_{\G}(N)_{u}&=
        \begin{cases}
        N_{\G^{-1}(u)}&\textup{ if $\G_{\leq u}\ne \emptyset$},\\
        0 &\textup{ otherwise.}
        \end{cases}
    \\
    \Lan_{\G}(N)_{u,v}&=
        \begin{cases}
        N_{\G^{-1}(u),\G^{-1}(v)}
        &\textup{ if $\G_{\leq u}\ne \emptyset$}\\
        0 &\textup{ otherwise.}
        \end{cases}
\end{align*}

\begin{lemma}\label{BarcodesContinuousExtensions}
    Given totally ordered sets $\T$ and $\U$ with $\T$ finite, a functor $\G \colon \T \to \U$, and a \pfd functor $N \colon \T \to \kvect$,
    \[\B{\Lan_{\G}(N)}=\{[\G(s),\G(t))\mid [s,t)\in \B{N},\, \G(s)<\G(t) \},\]
    where we define $\G(\infty)=\infty$. 
\end{lemma}

\begin{proof}
    It is easy to check directly that the map $N\mapsto \Lan_{\G}(N)$ preserves direct sums.  (More generally, left Kan extensions preserve coproducts, because they are left adjoints \cite[Proposition 6.1.5]{riehl2017category}, and left adjoints preserve colimits \cite[Theorem 4.5.3]{riehl2017category}.)  Hence, it suffices to prove the result in the case that the barcode of $N$ contains a single interval.  In this case, it is easy to check the result directly from the concrete description of $\Lan_{\G}(N)$ given above.
\end{proof}

\begin{lemma}\label{Lem:Push_Join}
    Given a line $\ell\in \L_{(0,\infty)}$ and $a^1, \ldots, a^n\in \R^2$, we have \[\push_{\ell}(\vee\{a^1,\ldots,a^n\}) = \max_i(\push_{\ell}(a^i)).\] 
\end{lemma}

\begin{proof}
    Let $b = \vee\{a^1,\ldots,a^n\}$. We prove the result assuming that $b$ is either on or below $\ell$; if $b$ is above $\ell$, a symmetric argument applies.  Given the assumption, $b$ and $\push_\ell(b)$ have the same $x$-coordinate, which is also the $x$-coordinate of some $a^j$.  Since $a^j\leq b$, $a^j$ also lies on or below $\ell$, so we have $\push_\ell(a^j)=\push_\ell(b)$.  Therefore, $\max_i(\push_{\ell}(a^i))\geq \push_{\ell}(b)$.  
    But since each $a^i\leq b$ and $\push_{\ell}$ respects the partial order in $\R^2$, the reverse inequality $\max_i(\push_{\ell}(a^i))\leq \push_{\ell}(b)$ also holds.  
\end{proof}

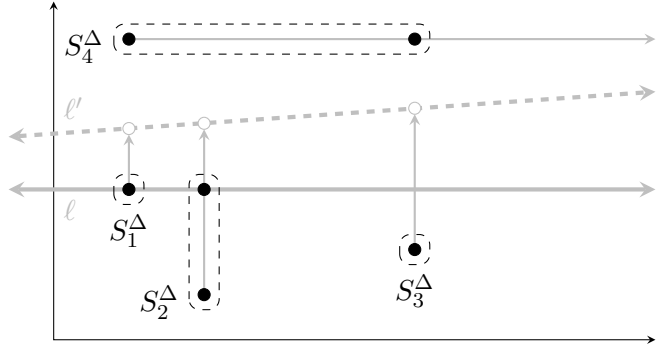
\begin{figure}
    \centering
    \begin{tikzpicture}[scale=1]
      \draw[<->] (0,4.5) -- (0,0) -- (8,0);
      \draw[<->,ultra thick,color=black!25] (-0.6,2) -- (8,2);
      \draw[<->,dashed,ultra thick,color=black!25] (-0.6,2.7) -- (8,3.3);  
      
      \node[below right,color=black!25] at (0,2) {$\ell$};
      \node[above right,color=black!25] at (0,2.8) {$\ell'$};

      \draw[color=black!25,fill=white] (1, 2.81) circle (0.08);      
      \draw[->,thick,black!25] (1,2) -- (1,2.73);
      \draw[color=black!25,fill=white] (2, 2.88) circle (0.08);      
      \draw[->,thick,black!25] (2,0.6) -- (2,2.8);
      \draw[color=black!25,fill=white] (4.8, 3.08) circle (0.08);      
      \draw[->,thick,black!25] (4.8,1.2) -- (4.8,3);

      \draw[->,thick,black!25] (1,4) -- (8,4);
      
      \foreach \p in {(1,2),(2,0.6),(2,2),(1,4),(4.8,4),(4.8,1.2)}
        { \draw[fill=black] \p circle (0.08); }

      \draw[dashed,rounded corners=4pt] (1.8,0.4) rectangle (2.2,2.2);
      \node[left] at (1.8,0.5) {$\XiSp^\ce_2$};
     
      \draw[dashed,rounded corners=4pt] (4.6,1) rectangle (5,1.4);
            \node[below] at (4.8,1.0) {$\XiSp^\ce_3$};
      \draw[dashed,rounded corners=4pt] (0.8,3.8) rectangle (5,4.2);
      \node[left] at (0.8,3.9) {$\XiSp^\ce_4$};

      \draw[dashed,rounded corners=4pt] (0.8,1.8) rectangle (1.2,2.2);
      \node[below] at (1.0,1.8) {$\XiSp^\ce_1$};

    \end{tikzpicture}
    \caption{The sets $S^\ce_i=S^{\ell'}_i$, for a horizontal line $\ell$ and $\ce$ the face containing the dual of the perturbed line $\ell'$, as specified in \cref{Sec:Query_Procedure}. 
    The figure illustrates that the sets $S^\ce_i$ are of the form $V_x$ or $H_y$ specified in the proof of \cref{Thm:QueriesMain}.}
    \label{fig:horizontalLine}
\end{figure}

\begin{proof}[Proof of \cref{Thm:QueriesMain}]
Let $\mathcal Q=\PCal^\ce \cap (\push_{\ell})^{-1}(\ell)$, i.e., $\mathcal Q$ is the subset of $\PCal^\ce$ which pushes onto $\ell$, rather than onto $\infty$.  Note that $\mathcal Q=\PCal^\ce$ unless $\ell$ is horizontal or vertical.  
Let $\G\colon \mathcal Q\to \ell$ denote the restriction of $\push_{\ell}$ to  $\mathcal Q$.  Note that $\G$ is a functor between totally ordered sets.  Let $N$ denote the restriction of $M^\ce$ (hence of $M$) to $\mathcal Q$.
Then $\B{N}$ is given as follows:
\[\B{N}=\{\, [a,b)\in \B{}^{\ce} \mid a,b\in \mathcal Q\cup\{\infty\}\,\} \cup \{\,[a,\infty)\mid [a,b)\in \B{}^{\ce},\, a\in \mathcal Q,\, b\in \PCal^{\ce}\setminus \mathcal Q\,\}.\]
\cref{BarcodesContinuousExtensions} then implies that
\[\B{\Lan_{\G}(N)}=\{\,[\push_\ell(a),\,\push_\ell(b)) \mid [a,b)\in \P^\ce,\ \push_\ell(a)<\push_\ell(b)\,\}.\]
  Thus, to prove the result, it suffices to show that $\Lan_{\G}(N)\cong M^\ell$.  Note that for all $a\in \ell$,
\[ 
    \Lan_{\G}(N)_{a} =
    \begin{cases}
        M_{\G^{-1}(a)} & \textup{ if $\G_{\leq a}\ne \emptyset$}, \\
            0 & \textup{ otherwise,}
    \end{cases}
\]
where we have used the fact that $N_{\G^{-1}(a)}=M_{\G^{-1}(a)}$.

We now define a natural transformation $f:\Lan_{\G}(N)\to M^\ell$ by 
\[ 
    f_{a} =
    \begin{cases}
        M_{\G^{-1}(a),a} & \textup{ if $\G_{\leq a}\ne \emptyset$}, \\
            0 & \textup{ otherwise.}
    \end{cases}
\]
Because the internal maps of both $\Lan_{\G}(N)$ and $M^\ell$ are given in terms of the structure maps of $M$, $f$ is indeed a natural transformation.  We now show that $f$ is a natural isomorphism.  Note that if $\G_{\leq a}=\emptyset$, then $M_a=0$.  Hence, $f_a$ is an isomorphism, because then its domain and codomain are both trivial.  

It remains to check that if $\G_{\leq a}\ne \emptyset$, then $f_a=M_{\G^{-1}(a),a}$ is an isomorphism.  
By \cref{Lem:GradesOfInfluenceAndIsomorphisms}, it suffices to show that 
\[ \I(\G^{-1}(a))= \I(a).\]
Since $\G^{-1}(a)\leq a$, we  have $\I(\G^{-1}(a))\subset \I(a)$.  To show that $\I(a)\subset \I(\G^{-1}(a))$, consider $s\in \I(a)$.  We need to show that $s\leq \G^{-1}(a)$.  By definition, we have 
\begin{equation}\label{Eq:G_Inv(a)}
    \G^{-1}(a) = \max\, \{\PCal^\ce_j\in \PCal^\ce\mid \push_{\ell}(\PCal^\ce_j)\leq a\}.
\end{equation}
Let us write $ \G^{-1}(a)=\PCal^\ce_i$.   
First, consider the case that $\ell$ is neither horizontal nor vertical.  If $s\not \leq \PCal^\ce_i$, then $s\in S^{\ce}_j$ for some $j>i$.  
Suppose $s\in S^{\ell}_{j'}$.  Since $\ell^*$ is in the closure of $\ce$, \cref{Lem:ContinuityOfPush} implies that if $t\in S^{\Delta}_k$ for some $k\leq j$, then $t\in S^{\ell}_{k'}$ for some $k'\leq j'$.  
Therefore, 
\begin{equation}\label{Eq:Levelset_Containment}
    \bigcup_{k=1}^{j} S^{\ce}_{k}\ \subset \ \bigcup_{k=1}^{j'} S^{\ell}_{k}.
\end{equation}
Since $s\leq a$ and $a\in \ell$, we have $\push_{\ell}(s)\leq a$.  Then $\push_{\ell}(t)\leq a$ for all $t\in S^{\ell}_k$ with $k\leq j'$.  Hence, by \cref{Eq:Levelset_Containment}, we have that $\push_{\ell}(t)\leq a$ for all $t\in S^{\ce}_k$ with $k\leq j$.  
\cref{Lem:Push_Join} now implies that $\push_{\ell}(\PCal^\ce_{j})\leq a$,  contradicting the definition of $\G^{-1}(a)$ as a maximum, as in \cref{Eq:G_Inv(a)}.

To finish the proof, we consider the case that $\ell$ is horizontal, say $\ell$ is the line $y=c$; the argument for vertical lines is symmetric.  By the way we have chosen $\ce$, $S^{\ce}$ consists of sets 
\[V_x\coloneqq\{s\in S\mid s_1=x \textup{ and } s_2\leq c\},\] as well as sets \[H_y\coloneqq \{s\in S \mid s_2=y\}\]  for $y>c$, with the sets $V_x$ ordered with respect to $x$ and each $V_x$ ordered before each $H_y$; see \cref{fig:horizontalLine} for an illustration.  It follows that for $a=(a_1,c)\in \ell$, we have 
\[ \G^{-1}(a) = \bigvee \left(\bigcup_{x\leq a_1} V_x\right)=\bigvee \{s\in S\mid s\leq a\},\]
and therefore 
\[ \I(\G^{-1}(a))=\{s\in S\mid s\leq a\}=\I(a). \qedhere \]  
\end{proof}

\subsection{Computational Details and Complexity of Queries}\label{Sec:Queries}
\label{Sec:Representing_And_Querying_Aug_Arrangement}\label{Sec:Data_Structure}

We now discuss the computational details of our queries and bound the complexity of a query.  We assume that the augmented arrangement $\S(M)$ has already been computed as described in \cref{Sec:Computing} below.  
Given a line $\ell\in \L_{[0,\infty]}$, the query of $\S(M)$ for $\B{M^\ell}$ proceeds in two steps.  
The first step performs a search for the face $\ce$ of $\cell(M)$ specified in \cref{Sec:Query_Procedure}.  
Once the face $\ce$ is selected, we obtain $\B{M^\ell}$ from $\P^\ce$ by applying $\push_\ell$ to the endpoints of each pair $(a,b) \in \P^\ce$.  

The search for face $\ce$ works as follows:  If $\ell\in \Lplus$ (i.e., $\ell$ is not vertical), then it suffices to find the cell of $\cell(M)$ containing $\ell^*$.  To do so, we use the point location data structure for $\cell(M)$ that is part of $\S(M)$; see \cref{Def:Aug_Arrangement}.  Since $\cell(M)$ has $O(\kappa^2)$ vertices, the search data structure has size $O(\kappa^2)$, and each point location query takes time $O(\log \kappa)$.  

If $\ell$ is vertical, then $\ell^*$ is not defined, and we require a different approach to find the face $\ce$.  We precompute a separate, simpler search data structure to handle the case of vertical lines:  Let $\Y$ denote the set of lines $k$ in $\cell(M)$ such that there is no other line in $\cell(M)$ with the same slope lying above $k$.  We take $\Y$ to be ordered by slope, and denote the $i^{\mathrm{th}}$ element of $\Y$ as $\Y_i$. We compute a 1-D array whose $i^{\mathrm{th}}$ entry is a pointer to the rightmost (unbounded) edge of $\cell(M)$ contained in $\Y_i$.  Given $\cell(M)$, computing this array takes $O(n)$ time, where $n=O(\kappa)$ is the number of anchor lines.  It follows from the description of the face $\ce$ given in  \cref{Rem:Concrete_Description_of_e} that for any vertical line $\ell$, we can find $\ce$ in $O(\log n)$ time by a binary search of this array.  

We are now ready to prove \cref{Thm:Query_Cost}, which bounds the cost of querying $\S(M)$ for a barcode $\B{M^{\ell}}$.

\begin{proof}[Proof of \cref{Thm:Query_Cost}]
  From the discussion above, it is clear that  once we have precomputed the appropriate data structures, finding the face $\ce$ takes $O(\log \kappa)$ time.  
  Each evaluation of $\push_\ell$ takes constant time, so computing $\B{M^\ell}$ from $\P^\ce$ takes total time $O(|\P^\ce|)$.  
  Thus, the total time to query $\S(M)$ for $\B{M^\ell}$ is $O(|\P^\ce|+\log \kappa)$.  If $\ell^* \in \ce$, then $|\P^\ce|=|\B{M^\ell}|$; this gives \cref{Thm:Query_Cost}~(i).  
  If $\ell^* \not\in \ce$, it may not be that $|\P^\ce|=|\B{M^\ell}|$, but we do have $|\P^\ce|=|\B{M^{\ell'}}|$ for $\ell'$ an arbitrarily small perturbation of $\ell$ with $\ell'^*\in \ce$.  \cref{Thm:Query_Cost}~(ii) follows.
\end{proof}

\section{Computing the Augmented Arrangement}\label{Sec:Computing}
Our algorithm for computing the augmented arrangement $\S(M)$ takes as input a minimal presentation of $M$.  The algorithm computes each of the following:
\begin{enumerate}
    \item the set of anchors of $M$,
    \item the line arrangement $\cell(M)$ and data structure for point location in $\cell(M)$,
    \item the barcode templates $\P^\ce$ at all faces $\ce$.
\end{enumerate}
  
Recall that the row and column labels of a minimal presentation for $M$ encode the 0th and 1st Betti numbers of $M$, respectively, and hence immediately yield $S$, the union of the support of these Betti numbers.  Given $S$, we can readily compute the set of anchors of $M$ in $O(\kappa)$ time and $O(\kappa)$ memory by iterating though the minimal grid containing $S$ in lexicographical order; we omit the easy details.  Since the anchors of $M$ correspond under point-line duality to the lines other than $\partial \H$ in $\cell(M)$, once the anchors have been computed we can use the standard algorithms mentioned in \cref{sec:lineArrangements} to compute $\cell(M)$ and a point location data structure for $\cell(M)$.

\subsection{Computing the Barcode Templates}\label{Sec:Algorithm_High_Level}
In the remainder of this section, we explain the computation of the barcode templates.  
For each face $\ce$ of $\cell(M)$, we define the surjective map 
\[ \lift^\ce : S \to \PCal^\ce \]
where $\lift^\ce(s) = \PCal^\ce_i$ for all $s\in S^\ce_i$.  
Given a presentation matrix $Q$ for $M$, we can obtain an ordered presentation matrix $Q^\ce$ for $M^\ce$ by applying $\lift^\ce$ to every row and column label of $Q$ and then permuting the rows and columns of $Q$ so that both the row and column labels are in increasing order.  (Here and throughout, our convention is that a permutation of the rows or columns of a labeled matrix also permutes the labels.)  We call such $Q^{\ce}$ an \emph{induced presentation} of $M^{\ce}$.  We let $\sigma_{\ce}$ and $\tau_{\ce}$ denote the respective choices of row and column permutations used to form $Q^{\ce}$.

By \cref{Prop:Barcode_from_Red_Mat}, the barcode template $\P^\ce$ can be read from an $RU$-decomposition of $Q^\ce$.  Na\"ively, we could compute $Q^{\ce}$ and its $RU$-decomposition from scratch at each face $\ce$, but this is inefficient.  Instead, our algorithm updates the $RU$-decomposition at one face to obtain the $RU$-decomposition at a neighboring face.

\begin{figure}[ht]
  \begin{center}
    \begin{tikzpicture}[dot/.style 2 args={circle,inner sep=1.5pt,fill,blue,label={#2:\textcolor{blue}{\footnotesize{$#1$}}},name=#1}]
        \draw[<->,ultra thick,black!30!white] (0,-1) -- (0,4.3);
        
        \draw[->,ultra thick,black!30!white] (0,-.5) -- (5.5,4.3);
        \draw[->,ultra thick,black!30!white] (0,1) -- (5.5,1.4);
        \draw[->,ultra thick,black!30!white] (0,1.5) -- (5.5,3.5);
        \draw[->,ultra thick,black!30!white] (0,2.1) -- (5.5,2.3);
        
        \node[dot={a}{above}] at (2,3.2) {};
        \node[dot={b}{right}] at (4.925,3.6) {};
        \node[dot={c}{right}] at (4.5,2.7) {};
        \node[dot={d}{above}] at (3,2.45) {};
        \node[dot={e}{below}] at (1.7,1.65) {};
        \node[dot={f}{left}] at (0.5,1.9) {};
        \node[dot={g}{left}] at (0.5,0.5) {};
        \node[dot={h}{below}] at (2.5,0.3) {};
        \node[dot={i}{right}] at (3.8,1.7) {};
        
        \draw[blue,very thick] (a) -- (b) -- (c) -- (d) -- (e) -- (f) -- (a) -- (d);
        \draw[blue,very thick] (c) -- (i) -- (h) -- (g) -- (e) -- (i);
    \end{tikzpicture}
  \end{center}
  \caption{A line arrangement $\cell(M)$ (in gray), together with its dual graph $G$ (in blue).  The walk $\pth$ through $G$ might visit the vertices in the order $a, b, c, d, e, f, e, g, h, i$.}
  \label{fig:dual_graph}
\end{figure}
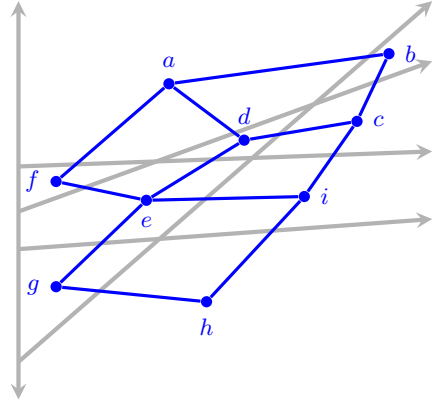

To elaborate, let $G$ denote the \emph{dual graph} of $\cell(M)$, i.e., the undirected graph with vertices the faces of $\cell(M)$ and edges the pairs of adjacent faces; see \cref{fig:dual_graph} for an illustration.  Let $\ce_0$ be the topmost face in $\cell(M)$; $\ce_0$ contains the dual of each line to the right of all points in $\XiSp$.  We compute a walk \[\pth = \ce_0, \ce_1, \ldots, \ce_n\] in $G$ that visits each face at least once, as explained in \cref{sec:path} below.  

For each $i$, we compute an $RU$-decomposition $Q^{\ce_i}=R^{\ce_i} U^{\ce_i}$ of an induced presentation $Q^{\ce_i}$ for $M^{\ce_i}$.  For $i<n$, we compute the $RU$-decomposition of $Q^{\ce_{i+1}}$ by updating the $RU$-decomposition of $Q^{\ce_{i}}$.  To enable efficient $RU$-updates, we will require not only that each $Q^{\ce_i}$ be an ordered presentation, but that it be \emph{strongly ordered}, in the following sense:

\begin{definition}
We say that $Q^{\ce_{i}}$ is strongly ordered if 
\begin{itemize}
\item It is an ordered presentation (see \cref{sec:ComputationPersistenceBarcodes}),
\item Given rows $j$ and $k$ of $Q$ with ($S$-valued) labels $s$ and $t$ such that $\lift^{\ce_i}(s)=\lift^{\ce_i}(t)$ and $s<t$ in the partial order on $\R^2$, we have $\sigma_{\ce_i}(j)<\sigma_{\ce_i}(k)$; and similarly for columns. 
\end{itemize}
\end{definition}

We compute $Q^{\ce_0}$ and its $RU$-decomposition as follows: Since $\ce_0$ contains the duals of all lines lying to the right of $S$, we have that for all $s \in S$, 
\[ \lift^{\ce_0}(s) = \big( \max\,\{t_1 \mid t\in S, t_2\leq s_2\}, s_2 \big). \]
Therefore, to compute $\PCal^{\ce_{0}}=\im \lift^{\ce_0}$ and $\lift^{\ce_0}$, it suffices to sort $S$ in colexicographical order (i.e., where $s$ is ordered before $t$ if and only if $s_2<t_2$ or $s_2=t_2$ and $s_1<t_1$), and then traverse $S$ in increasing order.  Computing $Q^{\ce_0}$ amounts to sorting the rows and columns of $Q$ with respect to the colexicographical order on labels, and then applying $\lift^{\ce_0}$ to each label.  We can then compute an $RU$-decomposition of $Q^{\ce_0}$ via the standard reduction mentioned in \cref{sec:ComputationPersistenceBarcodes}.  

For $0 \le i < n$, write $\ce = \ce_i$ and $\hat \ce = \ce_{i+1}$.  Suppose we have computed the $RU$-decomposition 
$Q^\ce=R^\ce U^\ce$.  We compute an $RU$-decomposition of $Q^{\hat \ce}$ for $M^{\hat \ce}$, in three steps: 
\begin{enumerate}
\item update $\PCal^\ce$ and $\lift^{\ce}$ to obtain  $\PCal^{\hat \ce}$ and $\lift^{\hat \ce}$,
\item permute the rows and columns of $R^\ce$ and $U^\ce$ to obtain a matrix decomposition $Q^{\hat \ce}=\bar R\bar U$ (which is not necessarily an $RU$-decomposition),
\item reduce the matrices $\bar R$ and $\bar U$ to obtain an $RU$-decomposition of $Q^{\hat \ce}$. 
 \end{enumerate}
 The next three subsections discuss these steps in detail.

\subsection{Updating the Template Points and Lift Map}\label{sec:updateLiftMap}
\begin{figure}[p]
  \centering

  \begin{subfigure}{\textwidth}
    \centering
    \begin{tikzpicture}
      \tikzstyle{every node}=[font=\scriptsize]
      \draw[<->,very thick,color=black!25] (-1.85,-2.8) -- (1,0.7); 
      \node[right,color=black!30] at (0.9,0.8) {$\ell$};
      
      \draw[color=black!25,fill=white] (0.43,0) circle (0.08);
      \draw[->,thick,color=black!25] (-4,0) -- (0.35,0);
      \draw[color=black!25,fill=white] (0,-0.528) circle (0.08);
      \draw[->,thick,color=black!25] (0,-2.2) -- (0,-0.62);
      \draw[color=black!25,fill=white] (-0.384,-1) circle (0.08);
      \draw[->,thick,color=black!25] (-3.8,-1) -- (-0.47,-1);
      \draw[color=black!25,fill=white] (-1.036,-1.8) circle (0.08);
      \draw[->,thick,color=black!25] (-2,-1.8) -- (-1.12,-1.8);

      \foreach \p in {(0,0),(-2,-1),(0,-1),(-2,-1.8)} 
        { \draw[fill=\tpcolor] \p circle (0.13); }

      \foreach \p in {(0,-1.8),(0,-2.2),(-3,0),(-4,0),(-3.8,-1),(-3.1,-1),(-2,-1.8)}
        { \fill \p circle (0.08); }

      \node[\tpcolor!70!black] at (0,0.3) {$\alpha$};
      \node[\tpcolor!70!black] at (0.55,-1) {$\PCal^{\ce}_{j-1}$};
      \node[\tpcolor!70!black] at (-2,-0.6) {$\PCal^{\ce}_{j-2}$};

      \draw[dashed,rounded corners=4pt] (-4.2,-0.2) rectangle (-2.8,0.2);
      \node[left] at (-4.2,0) {$\XiSp^\ce_j$};
      \draw[dashed,rounded corners=4pt] (-0.2,-2.4) rectangle (0.2,-1.6);
      \node[right] at (0.2,-2) {$\XiSp^\ce_{j-1}$};
      \draw[dashed,rounded corners=4pt] (-4,-1.2) rectangle (-2.9,-0.8);
      \node[left] at (-4,-1) {$\XiSp^\ce_{j-2}$};
      \draw[dashed,rounded corners=4pt] (-2.23,-2.03) rectangle (-1.77,-1.57);
    \end{tikzpicture}
    \hspace{0.5in}
    \begin{tikzpicture}
      \tikzstyle{every node}=[font=\scriptsize]
      
      \draw[<->,very thick,color=black!25] (-2.15,-2.8) -- (0.2,0.8);
      \node[right,color=black!30] at (0.2,0.8) {$\hat \ell$};

      \draw[color=black!25,fill=white] (-0.322,0) circle (0.08);
      \draw[->,thick,color=black!25] (-4,0) -- (-0.41,0);
      \draw[color=black!25,fill=white] (0,0.494) circle (0.08);
      \draw[->,thick,color=black!25] (0,-2.2) -- (0,0.4);
      \draw[color=black!25,fill=white] (-0.975,-1) circle (0.08);
      \draw[->,thick,color=black!25] (-3.8,-1) -- (-1.06,-1);
      \draw[color=black!25,fill=white] (-1.497,-1.8) circle (0.08);
      \draw[->,thick,color=black!25] (-2,-1.8) -- (-1.59,-1.8);

      \foreach \p in {(0,0),(-2,-1),(-2,0),(-2,-1.8)}
        { \draw[fill=\tpcolor] \p circle (0.13); }
        
      \foreach \p in {(0,-1.8),(-3,0),(-4,0),(-3.8,-1),(-3.1,-1),(-2,-1.8),(0,-2.2)}
        { \fill \p circle (0.08); }

      \node[\tpcolor!70!black] at (0.3,0) {$\alpha$};
      \node[\tpcolor!70!black] at (-2,0.4) {$\PCal^{\hat{\ce}}_{j-1}$};
      \node[\tpcolor!70!black] at (-2,-0.6) {$\PCal^{\hat{\ce}}_{j-2}$};
      
      \draw[dashed,rounded corners=4pt] (-4.2,-0.2) rectangle (-2.8,0.2);
      \node[left] at (-4.2,0) {$\XiSp^{\hat{\ce}}_{j-1}$};
      \draw[dashed,rounded corners=4pt] (-4,-1.2) rectangle (-2.9,-0.8);
      \node[left] at (-4,-1) {$\XiSp^{\hat{\ce}}_{j-2}$};
      \draw[dashed,rounded corners=4pt] (-0.2,-2.4) rectangle (0.2,-1.6);
      \node[right] at (0.2,-2) {$\XiSp^{\hat{\ce}}_{j}$};
      \draw[dashed,rounded corners=4pt] (-2.23,-2.03) rectangle (-1.77,-1.57);
    \end{tikzpicture}
    \caption{Generic case} 
    \label{fig:line_crosses_anchor1}
  \end{subfigure}
  
  \vspace{0.3in}

  \begin{subfigure}{\textwidth}
    \centering
    \begin{tikzpicture}
      \tikzstyle{every node}=[font=\scriptsize]
      \draw[<->,very thick,color=black!25] (-1.85,-2.8) -- (1,0.7); 
      \node[right,color=black!30] at (0.9,0.8) {$\ell$};
      
      \draw[color=black!25,fill=white] (0.43,0) circle (0.08);
      \draw[->,thick,color=black!25] (0,0) -- (0.35,0);
      \draw[color=black!25,fill=white] (0,-0.528) circle (0.08);
      \draw[->,thick,color=black!25] (0,-2.2) -- (0,-0.62);
      \draw[color=black!25,fill=white] (-0.384,-1) circle (0.08);
      \draw[->,thick,color=black!25] (-3.8,-1) -- (-0.47,-1);
      \draw[color=black!25,fill=white] (-1.036,-1.8) circle (0.08);
      \draw[->,thick,color=black!25] (-2,-1.8) -- (-1.12,-1.8);

      \foreach \p in {(0,0),(-2,-1),(0,-1),(-2,-1.8)} 
        { \draw[fill=\tpcolor] \p circle (0.13); }

      \foreach \p in {(0,0),(0,-1.8),(0,-2.2),(-3.8,-1),(-3.1,-1),(-2,-1.8)}
        { \fill \p circle (0.08); }

      \node[\tpcolor!70!black] at (0,0.34) {$\alpha$};
      \node[\tpcolor!70!black] at (0.55,-1) {$\PCal^{\ce}_{j-1}$};
      \node[\tpcolor!70!black] at (-2,-0.6) {$\PCal^{\ce}_{j-2}$};

      \draw[dashed,rounded corners=4pt] (-0.2,-0.2) rectangle (0.2,0.2);
      \node[left] at (-0.2,0) {$\XiSp^\ce_j$};
      \draw[dashed,rounded corners=4pt] (-0.2,-2.4) rectangle (0.2,-1.6);
      \node[right] at (0.2,-2) {$\XiSp^\ce_{j-1}$};
      \draw[dashed,rounded corners=4pt] (-4,-1.2) rectangle (-2.9,-0.8);
      \node[left] at (-4,-1) {$\XiSp^\ce_{j-2}$};
      \draw[dashed,rounded corners=4pt] (-2.23,-2.03) rectangle (-1.77,-1.57);
    \end{tikzpicture}
    \hspace{0.5in}
    \begin{tikzpicture}
      \tikzstyle{every node}=[font=\scriptsize]
      
      \draw[<->,very thick,color=black!25] (-2.15,-2.8) -- (0.2,0.8);
      \node[right,color=black!30] at (0.2,0.8) {$\hat \ell$};

      \draw[color=black!25,fill=white] (0,0.494) circle (0.08);
      \draw[->,thick,color=black!25] (0,-2.2) -- (0,0.4);
      \draw[color=black!25,fill=white] (-0.975,-1) circle (0.08);
      \draw[->,thick,color=black!25] (-3.8,-1) -- (-1.06,-1);
      \draw[color=black!25,fill=white] (-1.497,-1.8) circle (0.08);
      \draw[->,thick,color=black!25] (-2,-1.8) -- (-1.59,-1.8);

      \foreach \p in {(0,0),(-2,-1),(-2,-1.8)}
        { \draw[fill=\tpcolor] \p circle (0.13); }
        
      \foreach \p in {(0,0),(0,-1.8),(-3.8,-1),(-3.1,-1),(-2,-1.8),(0,-2.2)}
        { \fill \p circle (0.08); }

      \node[\tpcolor!70!black] at (0.34,0) {$\alpha$};
      \node[\tpcolor!70!black] at (-2,-0.6) {$\PCal^{\hat\ce}_{j-2}$};
      
      \draw[dashed,rounded corners=4pt] (-4,-1.2) rectangle (-2.9,-0.8);
      \node[left] at (-4,-1) {$\XiSp^{\hat\ce}_{j-2}$};
      \draw[dashed,rounded corners=4pt] (-0.2,-2.4) rectangle (0.2,0.2);
      \node[right] at (0.2,-2) {$\XiSp^{\hat\ce}_{j-1}$};
      \draw[dashed,rounded corners=4pt] (-2.23,-2.03) rectangle (-1.77,-1.57);
    \end{tikzpicture}
  
    \caption{Merge case} 
    \label{fig:line_crosses_anchor2}
  \end{subfigure}

  \vspace{0.3in}
  
  \begin{subfigure}{\textwidth}
    \centering
    \begin{tikzpicture}
      \tikzstyle{every node}=[font=\scriptsize]
      
      \draw[<->,very thick,color=black!25] (-2.15,-2.8) -- (0.2,0.8);
      \node[right,color=black!30] at (0.2,0.8) {$\hat \ell$};

      \draw[color=black!25,fill=white] (0,0.494) circle (0.08);
      \draw[->,thick,color=black!25] (0,-2.2) -- (0,0.4);
      \draw[color=black!25,fill=white] (-0.975,-1) circle (0.08);
      \draw[->,thick,color=black!25] (-3.8,-1) -- (-1.06,-1);
      \draw[color=black!25,fill=white] (-1.497,-1.8) circle (0.08);
      \draw[->,thick,color=black!25] (-2,-1.8) -- (-1.59,-1.8);

      \foreach \p in {(0,0),(-2,-1),(-2,-1.8)}
        { \draw[fill=\tpcolor] \p circle (0.13); }
        
      \foreach \p in {(0,0),(0,-1.8),(-3.8,-1),(-3.1,-1),(-2,-1.8),(0,-2.2)}
        { \fill \p circle (0.08); }

      \node[\tpcolor!70!black] at (0.34,0) {$\alpha$};
      \node[\tpcolor!70!black] at (-2,-0.6) {$\PCal^{\ce}_{j-1}$};
      
      \draw[dashed,rounded corners=4pt] (-4,-1.2) rectangle (-2.9,-0.8);
      \node[left] at (-4,-1) {$\XiSp^{\ce}_{j-1}$};
      \draw[dashed,rounded corners=4pt] (-0.2,-2.4) rectangle (0.2,0.2);
      \node[right] at (0.2,-2) {$\XiSp^{\ce}_{j}$};
      \draw[dashed,rounded corners=4pt] (-2.23,-2.03) rectangle (-1.77,-1.57);
    \end{tikzpicture}
    \hspace{0.5in}    
    \begin{tikzpicture}
      \tikzstyle{every node}=[font=\scriptsize]
      \draw[<->,very thick,color=black!25] (-1.85,-2.8) -- (1,0.7); 
      \node[right,color=black!30] at (0.9,0.8) {$\ell$};
      
      \draw[color=black!25,fill=white] (0.43,0) circle (0.08);
      \draw[->,thick,color=black!25] (0,0) -- (0.35,0);
      \draw[color=black!25,fill=white] (0,-0.528) circle (0.08);
      \draw[->,thick,color=black!25] (0,-2.2) -- (0,-0.62);
      \draw[color=black!25,fill=white] (-0.384,-1) circle (0.08);
      \draw[->,thick,color=black!25] (-3.8,-1) -- (-0.47,-1);
      \draw[color=black!25,fill=white] (-1.036,-1.8) circle (0.08);
      \draw[->,thick,color=black!25] (-2,-1.8) -- (-1.12,-1.8);

      \foreach \p in {(0,0),(-2,-1),(0,-1),(-2,-1.8)} 
        { \draw[fill=\tpcolor] \p circle (0.13); }

      \foreach \p in {(0,0),(0,-1.8),(0,-2.2),(-3.8,-1),(-3.1,-1),(-2,-1.8)}
        { \fill \p circle (0.08); }

      \node[\tpcolor!70!black] at (0,0.34) {$\alpha$};
      \node[\tpcolor!70!black] at (0.42,-1) {$\PCal^{\hat\ce}_{j}$};
      \node[\tpcolor!70!black] at (-2,-0.6) {$\PCal^{\hat\ce}_{j-1}$};

      \draw[dashed,rounded corners=4pt] (-0.2,-0.2) rectangle (0.2,0.2);
      \node[left] at (-0.2,0) {$\XiSp^{\hat\ce}_{j+1}$};
      \draw[dashed,rounded corners=4pt] (-0.2,-2.4) rectangle (0.2,-1.6);
      \node[right] at (0.2,-2) {$\XiSp^{\hat\ce}_{j}$};
      \draw[dashed,rounded corners=4pt] (-4,-1.2) rectangle (-2.9,-0.8);
      \node[left] at (-4,-1) {$\XiSp^{\hat\ce}_{j-1}$};
      \draw[dashed,rounded corners=4pt] (-2.23,-2.03) rectangle (-1.77,-1.57);
    \end{tikzpicture}
  
    \caption{Split case} 
    \label{fig:line_crosses_anchor3}
  \end{subfigure}
  \caption{Illustration of how the sets of template points and lift maps change when passing from a cell $\ce$ to an adjacent cell $\hat\ce$, containing the duals of lines $\ell$ and $\hat\ell$, respectively.}
  \label{fig:line_crosses_anchor}

\end{figure}
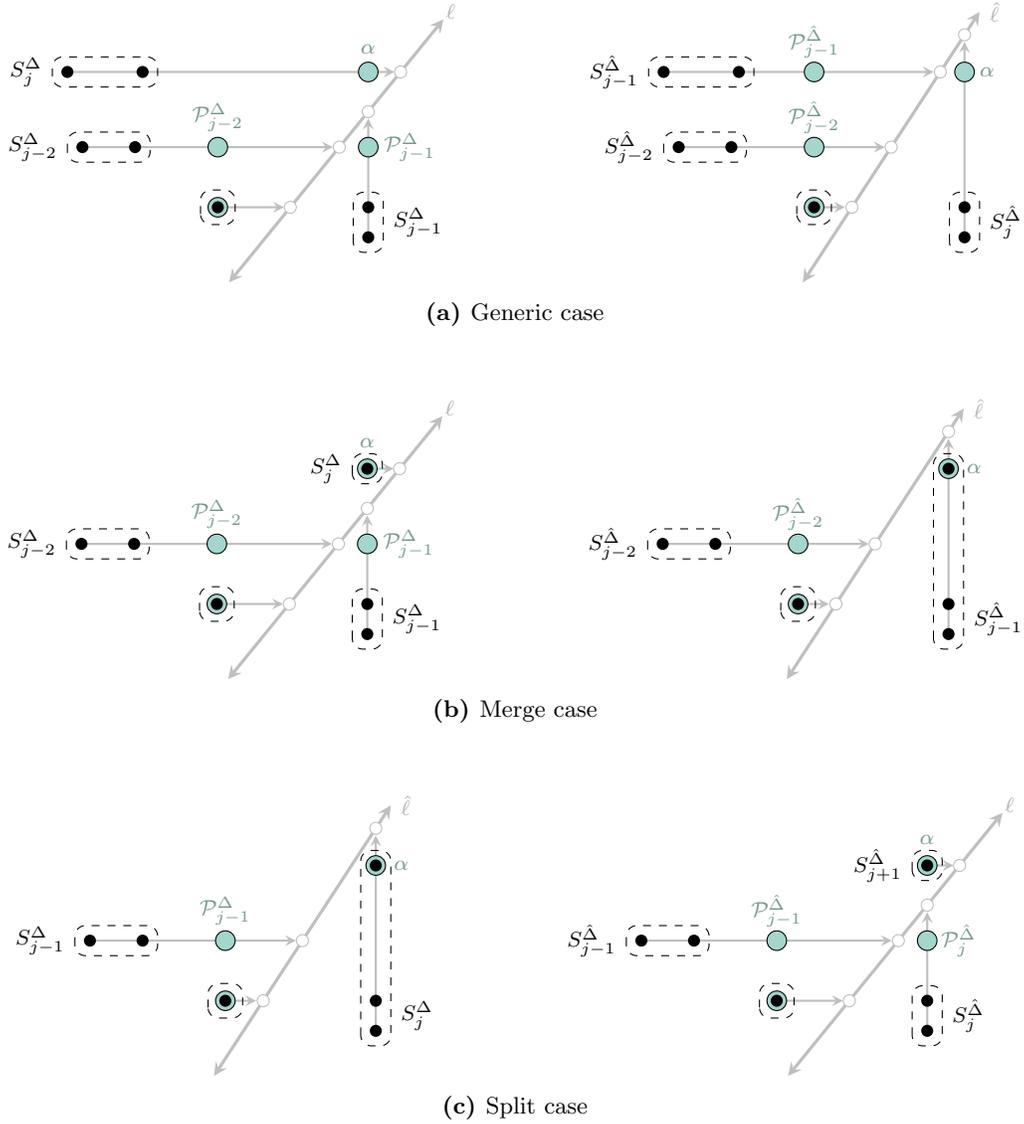

To explain how the template points and lift map change when passing from face $\ce$ to the adjacent face $\hat \ce$,   
suppose that the line containing the shared boundary of faces $\ce$ and $\hat \ce$ is dual to anchor $\alpha$.
Then $\alpha$ is an element of both $\PCal^\ce$ and $\PCal^{\hat \ce}$, as illustrated in  \cref{fig:line_crosses_anchor}. Suppose that $\alpha=\PCal^\ce_j$.  We distinguish between three cases, which are illustrated in \cref{fig:line_crosses_anchor}.  

\begin{enumerate}[leftmargin=1.2in]
  \item[\textbf{Generic case:}] $\alpha=\vee(s,t)$ for some incomparable $s,t\in S$.
  
  \item[\textbf{Merge case:}] $S^\ce_j=\{\alpha\}$.  In this case, $S^\ce_{j-1}$ and $S^\ce_{j}$ merge to form $S^{\hat \ce}_{j-1}$;   
  see \cref{levelsetupdates}\,(ii) below. 
  
  \item[\textbf{Split case:}] $S^{\hat \ce}_{j+1}=\{\alpha\}$.  In this case $S^\ce_j$ splits into $S^{\hat \ce}_j$ and $S^{\hat \ce}_{j+1}$;  
  see \cref{levelsetupdates}\,(iii) below. 
\end{enumerate}
It is easily checked that these cases are mutually exclusive and cover all possible cases.  If $S$ is generic in the sense that no two elements of $S$ share a coordinate, then only the first case arises, which explains our choice of name for this case.  The merge case and split cases are symmetric, in the sense that if passing from $\ce$ to  $\hat \ce$ is a merge, then passing from $\hat \ce$ to $\ce$ is a split, and vice versa; see~\cref{fig:line_crosses_anchor2,fig:line_crosses_anchor3}.    

The following result relates the totally ordered partitions $S^{\ce}$ and $S^{\hat \ce}$.  We omit the easy proof.  

\begin{lemma}\label{levelsetupdates}\mbox{}
\begin{itemize}
  \item[(i)] In the generic case, we have  
    $|S^{\hat\ce}|=|S^{\ce}|$ and 
        
    \[ S^{\hat\ce}_i =
        \begin{cases}
          S^{\ce}_{i} &\textup{ if }i\not\in\{j-1,j\},\\
          S^{\ce}_{j}\setminus \{\alpha\} &\textup{ if }i=j-1,\\
          S^{\ce}_{j-1}\cup\{\alpha\} &\textup{ if }i=j \textup{ and }\alpha\in S,\\
          S^{\ce}_{j-1}&\textup{ if }i=j \textup{ and }\alpha\not \in S.
        \end{cases} \]

  \item[(ii)] In the merge case, we have $|S^{\hat\ce}|=|S^{\ce}|-1$ and
    \[
    S^{\hat\ce}_i=
        \begin{cases}
        S^{\ce}_i &\textup{ if } i< j-1, \\
        S^{\ce}_{j-1}\cup \{\alpha\} &\textup{ if } i=j-1, \\
          S^{\ce}_{i+1} & \textup{ if } i > j-1.
        \end{cases}
        \]
  \item[(iii)] In the split case, we have $|S^{\hat\ce}|=|S^{\ce}|+1$ and
    \[
    S^{\hat\ce}_i=
        \begin{cases}
        S^{\ce}_i &\textup{ if } i< j, \\ 
         S^{\ce}_j\setminus \{\alpha\} &\textup{ if } i=j,\\
          \{\alpha\} &\textup{ if } i=j+1,\\
          S^{\ce}_{i-1} & \textup{ if } i> j+1.
        \end{cases} 
        \]   
\end{itemize}
\end{lemma}

The next result specifies how the template points and lift maps change as we pass from $\ce$ to $\hat\ce$.  To state the result, we let $\PCal^\ce_0=(-\infty,-\infty)$.  
\begin{lemma}\label{Lem:Lift_Update}\mbox{}
\begin{itemize}
  \item[(i)] In the generic case, we have $j>1$ and  
    \begin{align*}
      \PCal^{\hat\ce}_i&=
        \begin{cases}
        \PCal^{\ce}_i &\textup{ if } i\ne  j-1, \\ 
            \max (S^{\ce}_{j}\setminus \{\alpha\})\vee \PCal^{\ce}_{j-2} & \textup{ if } i=j-1,
        \end{cases}       \\ 
      \lift^{\hat\ce}(s) &=
        \begin{cases}
        \lift^\ce(s) &\textup{ if } s\not\in S^{\ce}_{j-1} \cup S^{\ce}_{j},\\ 
            \PCal^{\hat\ce}_j=\alpha  &\textup{ if } s\in S^{\ce}_{j-1} \cup \{\alpha\},\\
             \PCal^{\hat\ce}_{j-1} & \textup{ if } s\in S^{\ce}_{j}\setminus \{\alpha\}.
        \end{cases} 
    \end{align*}
  \item[(ii)] In the merge case, we have 
    \begin{align*}
    \PCal^{\hat\ce}_i&=
        \begin{cases}
        \PCal^{\ce}_i &\textup{ if } i< j-1, \\ 
          \PCal^{\ce}_{i+1} & \textup{ if } i\geq j-1.
        \end{cases}  \\
    \lift^{\hat\ce}(s) &=
        \begin{cases}
         \lift^\ce(s) &\textup{ if } s\not\in S^{\ce}_{j-1},\\ 
            \PCal^{\hat\ce}_{j-1}=\alpha  &\textup{ if } s\in S^{\ce}_{j-1}.
        \end{cases} 
        \end{align*}
  \item[(iii)] In the split case, we have
    \begin{align*}
      \PCal^{\hat\ce}_i&=
        \begin{cases}
        \PCal^{\ce}_i &\textup{ if } i< j, \\ 
          \max (S^{\ce}_{j}\setminus \{\alpha\})\vee \PCal^{\ce}_{j-1} &\textup{ if } i=j,\\
          \PCal^{\ce}_{i-1} & \textup{ if } i> j.
        \end{cases}    \\
    \lift^{\hat\ce}(s) &=
        \begin{cases}
        \lift^\ce(s) &\textup{ if } s\not\in S^{\ce}_{j}\setminus\{\alpha\},\\ 
            \PCal^{\hat\ce}_j  &\textup{ if } s\in S^{\ce}_{j}\setminus\{\alpha\}.
        \end{cases} 
        \end{align*}
\end{itemize}
\end{lemma}

\begin{proof}
The result follows readily from \cref{levelsetupdates} and the definitions of the sets $\PCal^\ce$ and $\PCal^{\hat\ce}$.  We give the details for case (i); the other cases are similar.  

In case (i), if $i\ne j-1$, then we have \[\bigcup_{k=1}^i S^{\ce}_k= \bigcup_{k=1}^i S^{\hat\ce}_k.\]
Therefore,
\[\PCal^{\hat\ce}_i=\bigvee\left( \bigcup_{k=1}^i S^{\hat\ce}_k\right)= \bigvee\left(\bigcup_{k=1}^i S^{\ce}_k\right)=\PCal^{ \ce}_i.\]  
Moreover, letting $\PCal^{\hat \ce}_0=(-\infty,-\infty)$, since $S^{\hat\ce}_{j-1}=S^{\ce}_{j}\setminus \{\alpha\}$ and $\PCal^{\hat \ce}_{j-2}=\PCal^{\ce}_{j-2}$, we have \[\PCal^{\hat\ce}_{j-1}=S^{\hat\ce}_{j-1}\vee \PCal^{\hat\ce}_{j-2}=\max (S^{\ce}_{j}\setminus \{\alpha\})\vee \PCal^{\ce}_{j-2}.\] 
If $i\not\in \{j-1,j\}$, then $S^{\hat\ce}_i=S^{\ce}_i$,  so for $s\in S^\ce_i$, we have 
\[\lift^{\hat\ce}(s)=\PCal^{\hat\ce}_i=\PCal^\ce_i=\lift^{\ce}(s).\]
If $s\in S^{\ce}_{j-1}\cup \{\alpha\}$, then $s\in S^{\hat\ce}_j$, so $\lift^{\hat\ce}(s)=\PCal^{\hat\ce}_j.$  Finally, if $s\in S^{\ce}_j\setminus \{\alpha\}$, then $s\in S^{\hat\ce}_{j-1}$, so $\lift^{\hat\ce}(s)=\PCal^{\hat\ce}_{j-1}.$
\end{proof}

Using \cref{Lem:Lift_Update}, it is straightforward to compute 
 $\PCal^{\hat\ce}$ and $\lift^{\hat\ce}$ by updating $\PCal^{\ce}$ and  $\lift^{\ce}$.  We discuss some low-level details of this update in \cref{sec:implementation}.
 
\cref{Lem:Lift_Update} also immediately implies the following. 

\begin{proposition}\label{lem:liftSwap}
    For $s,t\in S$, the following two conditions are equivalent:
    \begin{enumerate}
    \item $\lift^\ce(s) < \lift^\ce(t)$ and $\lift^{\hat\ce}(s) > \lift^{\hat\ce}(t)$
    \item $s \in S^{\ce}_{j-1}$, $t \in S^{\ce}_{j}\setminus\{ \alpha\}$, and $\alpha=s\vee t$.
    \end{enumerate}
\end{proposition}
Note that the equivalent conditions of \cref{lem:liftSwap} can hold only in the generic case.

\subsection{Permuting Matrices and Updating the \texorpdfstring{$RU$}{RU}-Decomposition}\label{updateRU}
Let $\rl(k)$ and $\cl(k)$ denote the $k^{\mathrm{th}}$ row and $k^{\mathrm{th}}$ column labels of $Q$, respectively.    
Then the $k^{\mathrm{th}}$ row label of $Q^\ce$ is $\lift^\ce(\rl(\sigma_{\ce}^{-1}(k)))$ and the $k^{\mathrm{th}}$ column label of $Q^\ce$ is $\lift^\ce(\cl(\tau_{\ce}^{-1}(k)))$.  Changing the row and column labels of $Q^\ce$  by replacing $\lift^{\ce}$ with $\lift^{\hat \ce}$ in these expressions, without changing the underlying matrix, yields a presentation $\tilde Q^{\hat\ce}$ of $M^{\hat\ce}$.   

Assume that $Q^\ce$ is strongly ordered.  Then in the merge and split cases, \cref{Lem:Lift_Update}\,(ii,iii) imply that $\tilde Q^{\hat\ce}$ is also strongly ordered, so we take $Q^{\hat\ce}=\tilde Q^{\hat\ce}$ and take the $RU$-decomposition of $Q^{\hat\ce}$ to be the same as for $Q^\ce$.  

In the generic case, \cref{lem:liftSwap} implies that $\tilde Q^{\hat\ce}$ is not an ordered presentation.  We identify permutation matrices $P$ and $P'$ such that $Q^{\hat\ce} \coloneqq P\tilde Q^{\hat\ce} P'$ is a strongly ordered presentation of $M^{\hat\ce}$.  Specifically, we take $P$ to be the permutation matrix that acts on rows of $\tilde Q^{\hat\ce}$ by transposing the block of rows labeled by elements of $S^{\ce}_{j-1} $ with the (adjacent) block of rows labeled by elements of $S^{\ce}_{j}\setminus \{\alpha\}$, while fixing the order of rows in each of the two blocks.  We take $P'$ to have the analogous action on columns.

To obtain an $RU$-decomposition of $Q^{\hat \ce}$, we update the $RU$-decomposition $Q^\ce=R^\ce U^\ce$ by reducing the matrices $\bar R\coloneqq P R^\ce P'$ and $\bar U\coloneqq (P')^{-1}U^\ce P'$ as explained in \cref{Sec:Updating_RU_Prelim}.  In principle, any of the update strategies mentioned in \cref{Sec:Updating_RU_Prelim} could be used in our setting.  Our complexity analysis of barcode template computation in \cref{sec:barcodeTemplateCost} considers the standard vineyard update strategy of \cite{cohen2006vines} which, as discussed in \cref{Sec:Updating_RU_Prelim}, decomposes the permutations represented by $P$ and $P'$ into a sequence of elementary transpositions of adjacent indices, performing a separate $RU$-update for each transposition.   In fact, in the generic setting that the coordinates of the row and column labels of $Q$ are distinct, each of $P$ and $P'$ is either an identity matrix or represents a single transposition of adjacent indices; thus, in this setting, the vineyard strategy and the global $RU$-update strategy described in \cref{Sec:Updating_RU_Prelim} are essentially the same.

The current implementation of barcode template computation in RIVET uses the vineyard strategy, together with some engineering tricks to improve practical performance \cite[Sections 8 and 9]{lesnick2015interactive}.  In brief, the code estimates the time required to compute the $RU$-decomposition of $Q^{\hat \ce}$ using vineyard updates, as well as the time required to compute an $RU$-decomposition of $Q^{\hat \ce}$ from scratch; it then chooses the option estimated to be faster.  

\begin{remark}
While the global approach to $RU$-updates described in \cref{Sec:Updating_RU_Prelim} is not implemented in RIVET, this appears to be a promising alternative; it is simpler than what is currently implemented, and it seems likely that it would perform comparably.  In future work, it would be interesting to empirically compare these and other $RU$-update schemes on the problem of barcode template computation. 
   \end{remark}
   
\subsection{Implementation Details}\label{sec:implementation}
To implement the update procedure described above, when passing from face $\ce$ to face $\hat \ce$ we must update the template points, lift map, and the $RU$-decomposition.  In addition, we must read the barcode template $\mathcal B^{\hat \ce}$ from the $RU$-decomposition $Q^{\hat \ce}=R^{\hat \ce}U^{\hat \ce}$.  To do all of this efficiently, we maintain several simple data structures.  We now discuss these.  

We represent $\PCal^{\ce}$ as an ordered list.  To represent $\lift^\ce$ in memory, each element $\PCal^\ce_k\in \PCal^\ce$ stores its level set $S^\ce_k = (\lift^\ce)^{-1}(\PCal^\ce_k)$ again as an ordered list; recall that these level sets form a partition of $S$.  In addition, to enable efficient updates of $\PCal^{\ce}$, we maintain a list $A$ of all anchors, and for each $\alpha\in A$, a pointer to the copy of $\alpha$ in the list representing $\PCal^{\ce}$; if $\alpha\not \in \PCal^{\ce}$, then this pointer is null.  Moreover, each $\alpha\in A$ maintains (possibly null) pointers to the anchors in $A$ immediately below and immediately to the left of $\alpha$.  As we pass from cell $\ce$ to cell $\hat\ce$ we can update these data structures in constant time using the formulae of \cref{Lem:Lift_Update}.  We omit the straightforward details. 

Following \cite{cohen2006vines}, we represent $R^{\ce}$ and $U^{\ce}$ in memory by a column-sparse matrix and a row-sparse matrix, respectively; see \cref{Sec:Updating_RU_Prelim}.  To enable efficient computation of the permuted matrices $\bar R$ and $\bar U$ from $R^{\ce}$ and $U^{\ce}$, 
each $s\in S$ maintains the minimum and maximum indices of rows $k$ of $Q^\ce$ with $\rl(\sigma_\ce^{-1}(k)) = s$.  (Since 
$Q^{\ce}$ is strongly ordered,  such rows are ordered consecutively.)  Similarly, $s$ also maintains the minimum and maximum indices of columns $k$ of $Q^\ce$ with $\cl(\tau_\ce^{-1}(k)) = s$.  This data allows us to determine in constant time the indices of the blocks of rows and columns that must be exchanged to obtain $\bar R$ and $\bar U$ from $R^{\ce}$ and $U^{\ce}$.

To efficiently obtain the barcode template $\mathcal{B}^{\hat \ce}$ from the $RU$-decomposition $Q^{\hat \ce}=R^{\hat \ce}U^{\hat \ce}$, we require that each row and column of $Q^{\hat \ce}$ have access to its label.  Since we maintain the map $\lift^{\hat \ce}\colon S\to \PCal^{\hat \ce}$ by storing its level sets, it suffices for each row $k$ of $Q^{\hat \ce}$ to maintain a pointer to $\rl( \sigma^{-1}_{\hat \ce}(k))\in S$, and similarly for columns.  Once initialized at face $\ce_0$, these pointers never need to be updated as we traverse the path $\pth$.  

\subsection{Computing a Walk Through the Faces}\label{sec:path}
Recall that $\ce_0$ is the topmost face in $\cell(M)$ and that we compute a walk $\pth = \ce_0, \ce_1, \ldots, \ce_n$ in $G$ which visits each face at least once.  We call such a walk \emph{valid}.  We now consider the problem of computing a valid walk.    
Our complexity analysis in \cref{sec:barcodeTemplateCost} below holds for any valid walk $\pth$ that traverses each edge of $G$ at most a constant number of times.
However, in practice we want to choose $\pth$ to minimize the total cost of $RU$-updates performed by our algorithm.  We describe a simple heuristic for this which treats $G$ as a weighted graph and chooses $\pth$ to be a valid walk of small weight.

We take the weight $w_x$ of an edge $x = [\ce, \hat\ce]$ in $G$ to be the total number of row and column transpositions required to update a strongly ordered induced presentation $Q^\ce$ of $M^\ce$ to a strongly ordered induced presentation $Q^{\hat\ce}$ of $M^{\hat\ce}$.  We regard $w_x$ as a proxy for the work required to obtain an $RU$-decomposition of $Q^{\hat\ce}$ from one for $Q^\ce$.   
It is easily checked that $w_x$ 
depends only on the line $l$ of $\cell(M)$ separating $\ce$ and $\hat\ce$.  
Moreover, $w_x\ne 0$ only when the anchor $\alpha\coloneqq l^*$ is the join of incomparable points in $S$, which corresponds to the generic case of \cref{sec:updateLiftMap}.  In this case, to give an explicit formula for $w_x$, suppose $\alpha = \PCal^\ce_j$. 
For $i=0,1$, define 
\[w_{i,j-1} = \sum_{s\in S^{\ce}_{j-1}} \beta_i^M(s), \qquad \qquad w_{i,j} = \sum_{s\in S^{\ce}_{j}\setminus \{\alpha\}} \beta_i^M(s).\]
Then by \cref{lem:liftSwap},  $w_x=w_{0,j-1}w_{0,j}+w_{1,j-1}w_{1,j}$.      

We would like to choose $\pth$ to be a valid walk of minimum total weight, but computing such a walk is an instance of a classic NP-hard problem and likely itself hard.  Therefore, we instead compute a walk $\pth$ whose weight is \emph{approximately} minimum.  We mention two classical approximation schemes for this, yielding approximation ratios of $2$ and $\frac{3}{2}$, respectively: Let $\pth^*$ be a valid walk of minimum total weight.  For the $2$-approximation, we compute a minimum spanning tree $\M$ for $G$, e.g., via Kruskal's Algorithm \cite{cormen2009introduction}. 
Depth-first search of $\M$ starting at $\ce_0$ yields a valid walk $\pth$ in $\M$ traversing  each edge of $\M$ at most twice.  
Since $\length(\M)\leq \length(\pth^*)$, we have that $\length(\pth)\leq 2\,\length(\pth^*)$.  Alternatively, \cite{hoogeveen1991analysis} shows that a variant of the Christofides $\frac{3}{2}$-approximation algorithm for the traveling salesman problem on a metric graph \cite{christofides1976} yields a valid walk $\pth$ with $\length(\pth)\leq \frac{3}{2}\length(\pth^*)$.

\section{Complexity of Computing the Augmented Arrangement}\label{sec:barcodeTemplateCost}

We now turn to the proof of \cref{SimpleAugArrComplexity}\,(ii), which bounds the time and memory cost of our algorithm for computing the augmented arrangement $\S(M)$. 

\begin{proof}[Proof of \cref{SimpleAugArrComplexity}\,(ii)]
Recall that computing $\S(M)$ entails computing the anchors, the line arrangement $\cell(M)$, the point location data structure, and the barcode templates.
To prove \cref{SimpleAugArrComplexity}\,(ii), we analyze the cost of each of these steps.
\Cref{Table:Arrangement_Cost} summarizes the cost of each step.

\begin{table}[ht]\footnotesize
  \begin{center}
    \begin{tabular}{M{6.5cm} m{3.5cm}  l }
      \toprule
      \textbf{Computation} & \textbf{Time} & \textbf{Memory}    \\
      \midrule
      Set of anchors of $M$ & $O(\kappa)$ & $O(\kappa)$ \\
      \midrule[0.3pt]
      Line arrangement and data structure for point location & $O(\kappa^2 \log \kappa)$ & $O(\kappa^2)$ \\
      \midrule[0.3pt]
      Barcode templates  & $O(m^3 \kappa + \kappa^2 (m + \log\kappa) )$ & $O(m^2+m\kappa^2)$  \\
      \bottomrule
    \end{tabular}
  \end{center}
  \caption{Cost of augmented arrangement sub-computations}
  \label{Table:Arrangement_Cost}
\end{table}

As noted at the beginning of \cref{Sec:Computing}, computing the set of anchors requires $O(\kappa)$ time and $O(\kappa)$ memory.  Recall from \cref{sec:lineArrangements} that the Bentley-Ottmann algorithm \cite{deberg2008computation,bentley1979algorithms} computes a line arrangement with $l$ lines and $v$ vertices in $O((l + v) \log l)$ time and $O(l+v)$ memory.  Moreover, for $e$ the number of edges in the arrangement, we can compute a point location data structure for the arrangement in $O(e\log e)$ time and $O(e)$ memory.
As mentioned in \cref{sec:arr_def}, $\cell(M)$ has $O(\kappa)$ lines, $O(\kappa^2)$ vertices, and $O(\kappa^2)$ edges.  Thus,  computing $\cell(M)$ via Bentley-Ottmann requires $O(\kappa^2 \log \kappa)$ time and $O(\kappa^2)$ memory, and a point location data structure for $\cell(M)$  can be computed in $O(\kappa^2 \log \kappa)$ time and $O(\kappa^2)$ memory.\footnote{As noted in \cref{sec:lineArrangements}, there exist algorithms for computing line arrangements that are asymptotically faster than Bentley-Ottmann.  However, since computing the point location data structure takes $O(\kappa^2 \log \kappa)$ time, Bentley-Ottmann is not a  bottleneck in the time complexity of our query scheme.}

The substeps involved in the computation of the barcode templates are listed in \cref{Table:Barcode_Cost}, together with their runtimes.  We now discuss the cost of each substep.
\begin{table}[ht]\footnotesize
  \centering
  \begin{tabular}{M{8cm} m{3.3cm} }
    \toprule
    \textbf{Computation} & \textbf{Time} \\
    \midrule
    Walk $\pth$ (via the 2-approximation algorithm for the optimal walk, using Kruskal's MST algorithm) & $O(\kappa^2 \log \kappa)$  \\
    \midrule[0.3pt]
    $\PCal^{\ce_0}$ and $\lift^{\ce_0}$  & $O(m \log m)$ \\
    \midrule[0.3pt]
    $RU$-decomposition at $\ce_0$ & $O(m^3)$  \\ 
    \midrule 
    updates to $\PCal^\ce$ and $\lift^{\ce}$ & $O(\kappa^2)$  \\
    \midrule[0.3pt]
    updates to permutations and $RU$-decompositions & $ O(m^3 \kappa)$  \\
    \midrule[0.3pt]
reading barcode templates & $ O(m\kappa^2)$  \\
    \bottomrule
  \end{tabular}
  \caption{Cost of barcode template sub-computations}
  \label{Table:Barcode_Cost}
\end{table}

On a graph with $E$ edges and $V$ vertices, Kruskal's algorithm runs in time $O(E \log V)$.  As noted in \cref{sec:arr_def}, the numbers of vertices,  edges, and faces in $\cell(M)$ are each  $O(\kappa^2)$.  By duality, $G$ also has $O(\kappa^2)$ vertices and edges.  Thus, the runtime of Kruskal's algorithm on $G$ is $O(\kappa^2 \log \kappa)$,  and the runtime of depth-first search on the resulting spanning tree is $O(\kappa^2)$.
Hence, the cost of computing $\pth$ using the approach based on Kruskal's algorithm described at end of  \cref{Sec:Algorithm_High_Level} is $O(\kappa^2 \log \kappa)$.

Recall that to compute the $RU$-decomposition and barcode templates at the initial cell $\ce_0$, we first compute the map $\lift^{\ce_0}$ and the ordered presentation $Q^{\ce_0}$.  As explained in \cref{Sec:Algorithm_High_Level}, these computations are simple applications of sorting; they require time $O(m \log m)$. 
Computing the $RU$-decomposition of $Q^{\ce_0}$, and hence the barcode template $\B{}^{\ce_0}$, then requires $O(m^3)$ time.  

We now analyze the cost of performing the updates at each cell $\ce_i$ for $i>0$.  As we move from cell $\ce$ to cell $\hat\ce$, we can update both $\PCal^\ce$ and $\lift^\ce$, as well as the associated pointers described in \cref{sec:implementation}, in constant time.  Since the walk $\pth$ has length $O(\kappa^2)$, the total cost of all such updates is $O(\kappa^2)$.

For our cost analysis of permuting the matrices and updating the RU-decompositions, we consider only the vineyard update strategy of \cite{cohen2006vines}.  
We will show that for this strategy, the total cost of permuting the matrices and updating the $RU$-decompositions is $O(m^3\kappa)$.   
To do so, we count the number of transpositions involving a pair of rows of $Q$ with given labels $a$ and $b$.  These rows can swap only if $a$ and $b$ are incomparable, and they swap only when crossing the line dual to the anchor $a\vee b$.  As we traverse $\pth$, each anchor line is crossed $O(\kappa)$ times.  There are at most $O(m^2)$ label pairs, so at most $O(m^2\kappa)$ row transpositions are performed in total.  The same analysis applies to the column transpositions.  By the results of \cite{cohen2006vines}, if the row transpositions of $R$ and column transpositions of $U$ are performed implicitly, then the $RU$-update corresponding to each row or column transposition requires $O(m)$ time.   Thus, the total cost of all permutations and $RU$-updates is $O(m^3\kappa)$.  
Reading a barcode template from the $RU$-decomposition requires $O(m)$ time. Therefore, the total cost of reading all barcode templates is $O(m\kappa^2)$.

The memory cost of computing all barcode templates is dominated by the cost of storing a presentation matrix, which is $O(m^2)$, and the cost of storing the output, which is $O(m\kappa^2)$.  Thus, computing all barcode templates requires $O(m^3 \kappa + \kappa^2 (m + \log\kappa) )$ time and $O(m^2+m\kappa^2)$ memory.  \cref{SimpleAugArrComplexity}\,(ii) now follows.  
\end{proof}

\section*{Acknowledgments}

This work was funded in part by NSF grant DMS-1606967 and a grant from the Simons Foundation (Award ID 963845).  We thank our collaborators on the RIVET software, especially
 Anway De, Bryn Keller, Simon Segert, Alex Yu, and Roy Zhao.  Our joint work on RIVET has helped shape our understanding of the ideas in this paper.  We also thank the anonymous reviewers of both \cite{lesnick2015interactive} and this paper for valuable feedback.

We, the authors, developed all of the ideas of this paper and wrote all of the text ourselves.  Claude and ChatGPT provided proofreading input, catching many typos and making several helpful suggestions about the exposition.

\bibliographystyle{abbrv}
{\footnotesize \bibliography{fast_queries_references} }

@InProceedings{fugacci_et_al:LIPIcs.SoCG.2019.37,
  author =	{Fugacci, Ulderico and Kerber, Michael},
  title =	{{Chunk Reduction for Multi-Parameter Persistent Homology}},
  booktitle =	{35th International Symposium on Computational Geometry (SoCG 2019)},
  pages =	{37:1--37:14},
  series =	{Leibniz International Proceedings in Informatics (LIPIcs)},
  ISBN =	{978-3-95977-104-7},
  ISSN =	{1868-8969},
  year =	{2019},
  volume =	{129},
  editor =	{Barequet, Gill and Wang, Yusu},
  publisher =	{Schloss Dagstuhl -- Leibniz-Zentrum f{\"u}r Informatik},
  address =	{Dagstuhl, Germany},
  URL =		{https://drops.dagstuhl.de/entities/document/10.4230/LIPIcs.SoCG.2019.37},
  URN =		{urn:nbn:de:0030-drops-104413},
  doi =		{10.4230/LIPIcs.SoCG.2019.37}
}

@InProceedings{morozov_et_al:LIPIcs.SoCG.2025.68,
  author =	{Morozov, Dmitriy and Skraba, Primoz},
  title =	{{Persistent (Co)Homology in Matrix Multiplication Time}},
  booktitle =	{41st International Symposium on Computational Geometry (SoCG 2025)},
  pages =	{68:1--68:16},
  series =	{Leibniz International Proceedings in Informatics (LIPIcs)},
  ISBN =	{978-3-95977-370-6},
  ISSN =	{1868-8969},
  year =	{2025},
  volume =	{332},
  editor =	{Aichholzer, Oswin and Wang, Haitao},
  publisher =	{Schloss Dagstuhl -- Leibniz-Zentrum f{\"u}r Informatik},
  address =	{Dagstuhl, Germany},
  URL =		{https://drops.dagstuhl.de/entities/document/10.4230/LIPIcs.SoCG.2025.68},
  URN =		{urn:nbn:de:0030-drops-232204},
  doi =		{10.4230/LIPIcs.SoCG.2025.68}
}

@InProceedings{milosavljevic2011zigzag,
  author    = {Milosavljevi\'{c}, Nikola and Morozov, Dmitriy and \v{S}kraba, Primo\v{z}},
  title     = {Zigzag Persistent Homology in Matrix Multiplication Time},
  booktitle = {Proceedings of the 27th ACM Symposium on Computational Geometry},
  editor    = {Hurtado, Ferran and van Kreveld, Marc J.},
  year      = {2011},
  pages     = {216--225},
  publisher = {ACM},
  doi       = {10.1145/1998196.1998228}
}

@article{bauer2017phat,
  author    = {Bauer, Ulrich and Kerber, Michael and Reininghaus, Jan and Wagner, Hubert},
  title     = {{PHAT} -- Persistent Homology Algorithms Toolbox},
  journal   = {Journal of Symbolic Computation},
  volume    = {78},
  pages     = {76--90},
  year      = {2017},
  doi       = {10.1016/j.jsc.2016.03.008}
}

@article{mulmuley1990fast,
  title   = {A fast planar partition algorithm, {I}},
  author  = {Mulmuley, Ketan},
  journal = {Journal of Symbolic Computation},
  volume  = {10},
  number  = {3--4},
  pages   = {253--280},
  year    = {1990},
  doi     = {10.1016/S0747-7171(08)80064-8}
}

@article{seidel1991simple,
  title   = {A simple and fast incremental randomized algorithm for computing
             trapezoidal decompositions and for triangulating polygons},
  author  = {Seidel, Raimund},
  journal = {Computational Geometry: Theory and Applications},
  volume  = {1},
  number  = {1},
  pages   = {51--64},
  year    = {1991}
}

@article{haran2008experimental,
  title   = {An experimental study of point location in planar arrangements in {CGAL}},
  author  = {Haran, Idit and Halperin, Dan},
  journal = {ACM Journal of Experimental Algorithmics},
  volume  = {13},
  pages   = {2.3--2.32},
  year    = {2008},
  doi     = {10.1145/1412228.1412237}
}

@inproceedings{hemmer2012improved,
  title         = {Improved implementation of point location in general two-dimensional subdivisions},
  author        = {Hemmer, Michael and Kleinbort, Michal and Halperin, Dan},
  booktitle     = {Algorithms -- {ESA} 2012},
  series        = {Lecture Notes in Computer Science},
  volume        = {7501},
  pages         = {611--623},
  publisher     = {Springer},
  year          = {2012},
  doi           = {10.1007/978-3-642-33090-2_53},
  eprint        = {1205.5434},
  archivePrefix = {arXiv},
  primaryClass  = {cs.CG},
  note          = {Full version: \href{https://arxiv.org/abs/1205.5434}{arXiv:1205.5434}}
}

@article{hemmer2016optimal,
  title   = {Optimal randomized incremental construction for guaranteed logarithmic planar point location},
  author  = {Hemmer, Michael and Kleinbort, Michal and Halperin, Dan},
  journal = {Computational Geometry: Theory and Applications},
  volume  = {58},
  pages   = {110--123},
  year    = {2016}
}

@misc{keller2018persistent,
  author       = {Keller, Bryn and Lesnick, Michael and Willke, Theodore L.},
  title        = {{PHoS}: Persistent Homology for Virtual Screening},
  year         = {2018},
  month        = oct,
  howpublished = {ChemRxiv preprint, \url{https://chemrxiv.org/engage/chemrxiv/article-details/60c73f1cf96a00ccfe28604b}}
}

@incollection{cgal-arrangement-2,
  author    = {Ron Wein and Eric Berberich and Efi Fogel and
               Dan Halperin and Michael Hemmer and Oren Salzman and
               Baruch Zukerman},
  title     = {{2D} Arrangements},
  booktitle = {{CGAL} User and Reference Manual},
  publisher = {{CGAL Editorial Board}},
  edition   = {{6.2}},
  year      = {2026},
  url       = {https://doc.cgal.org/6.2/Arrangement_on_surface_2/index.html}
}

@article{sarnak1986planar,
  title={Planar point location using persistent search trees},
  author={Sarnak, Neil and Tarjan, Robert E},
  journal={Communications of the ACM},
  volume={29},
  number={7},
  pages={669--679},
  year={1986},
  publisher={ACM New York, NY, USA}
}

@article{bentley1979algorithms,
  title={Algorithms for reporting and counting geometric intersections},
  author={Bentley, Jon Louis and Ottmann, Thomas A.},
  journal={IEEE Transactions on computers},
  volume={28},
  number={9},
  pages={643--647},
  year={1979},
  publisher={IEEE}
}

@article{edelsbrunner1986optimal,
  title={Optimal point location in a monotone subdivision},
  author={Edelsbrunner, Herbert and Guibas, Leonidas J and Stolfi, Jorge},
  journal={SIAM Journal on Computing},
  volume={15},
  number={2},
  pages={317--340},
  year={1986},
  publisher={SIAM}
}

@article{kirkpatrick1983optimal,
  title={Optimal search in planar subdivisions},
  author={Kirkpatrick, David},
  journal={SIAM Journal on Computing},
  volume={12},
  number={1},
  pages={28--35},
  year={1983},
  publisher={SIAM}
}

@InProceedings{kerber_et_al:LIPIcs.SoCG.2026.62,
  author =	{Kerber, Michael and Wang, Elena Xinyi},
  title =	{{Computing the Bottleneck Distance Between Persistent Homology Transforms}},
  booktitle =	{42nd International Symposium on Computational Geometry (SoCG 2026)},
  pages =	{62:1--62:15},
  series =	{Leibniz International Proceedings in Informatics (LIPIcs)},
  ISBN =	{978-3-95977-418-5},
  ISSN =	{1868-8969},
  year =	{2026},
  volume =	{367},
  editor =	{Ahn, Hee-Kap and Hoffmann, Michael and Nayyeri, Amir},
  publisher =	{Schloss Dagstuhl -- Leibniz-Zentrum f{\"u}r Informatik},
  address =	{Dagstuhl, Germany},
  URL =		{https://drops.dagstuhl.de/entities/document/10.4230/LIPIcs.SoCG.2026.62},
  URN =		{urn:nbn:de:0030-drops-258693},
  doi =		{10.4230/LIPIcs.SoCG.2026.62}
}

@article{loiseaux2022fast,
  title={Multi-parameter Module Approximation: an efficient and interpretable invariant for multi-parameter persistence modules with guarantees},
  author={Loiseaux, David and Carri{\`e}re, Mathieu and Blumberg, Andrew J},
  journal={Journal of Applied and Computational Topology},
  volume={9},
  number={4},
  pages={1--60},
  year={2025},
  publisher={Springer}
}

@article{chacholski2017combinatorial,
  title={Combinatorial presentation of multidimensional persistent homology},
  author={Chach{\'o}lski, Wojciech and Scolamiero, Martina and Vaccarino, Francesco},
  journal={Journal of Pure and Applied Algebra},
  volume={221},
  number={5},
  pages={1055--1075},
  year={2017},
  publisher={Elsevier}
}

@InProceedings{pmlr-v235-scoccola24a,
  title = 	 {Differentiability and Optimization of Multiparameter Persistent Homology},
  author =       {Scoccola, Luis and Setlur, Siddharth and Loiseaux, David and Carri\`{e}re, Mathieu and Oudot, Steve},
  booktitle = 	 {Proceedings of the 41st International Conference on Machine Learning},
  pages = 	 {43986--44011},
  year = 	 {2024},
  editor = 	 {Salakhutdinov, Ruslan and Kolter, Zico and Heller, Katherine and Weller, Adrian and Oliver, Nuria and Scarlett, Jonathan and Berkenkamp, Felix},
  volume = 	 {235},
  series = 	 {Proceedings of Machine Learning Research},
  month = 	 {21--27 Jul},
  publisher =    {PMLR},
  url = 	 {https://proceedings.mlr.press/v235/scoccola24a.html}
}

@article{busaryev2010tracking,
  title={Tracking a generator by persistence},
  author={Busaryev, Oleksiy and Dey, Tamal K and Wang, Yusu},
  journal={Discrete Mathematics, Algorithms and Applications},
  volume={2},
  number={04},
  pages={539--552},
  year={2010},
  publisher={World Scientific}
}

@article{fernandes2025computation,
  title={Computation of gamma-linear projected barcodes for multiparameter persistence},
  author={Fernandes, Alex and Oudot, Steve and Petit, Fran{\c{c}}ois},
  journal={Journal of Applied and Computational Topology},
  volume={9},
  number={2},
  pages={12},
  year={2025},
  publisher={Springer}
}

@book{preparata2012computational,
  title={Computational geometry: an introduction},
  author={Preparata, Franco P and Shamos, Michael I},
  year={2012},
  publisher={Springer Science \& Business Media}
}

@article{hickok2022computing,
  title={Computing persistence diagram bundles},
  author={Hickok, Abigail},
  journal={arXiv preprint arXiv:2210.06424},
  year={2022}
}

@InProceedings{kim_et_al:LIPIcs.SoCG.2025.64,
  author =	{Kim, Donghan and Kim, Woojin and Lee, Wonjun},
  title =	{{Super-Polynomial Growth of the Generalized Persistence Diagram}},
  booktitle =	{41st International Symposium on Computational Geometry (SoCG 2025)},
  pages =	{64:1--64:20},
  series =	{Leibniz International Proceedings in Informatics (LIPIcs)},
  ISBN =	{978-3-95977-370-6},
  ISSN =	{1868-8969},
  year =	{2025},
  volume =	{332},
  editor =	{Aichholzer, Oswin and Wang, Haitao},
  publisher =	{Schloss Dagstuhl -- Leibniz-Zentrum f{\"u}r Informatik},
  address =	{Dagstuhl, Germany},
  URL =		{https://drops.dagstuhl.de/entities/document/10.4230/LIPIcs.SoCG.2025.64},
  URN =		{urn:nbn:de:0030-drops-232162},
  doi =		{10.4230/LIPIcs.SoCG.2025.64}
}

@article{morozov2021output,
  title={Output-sensitive computation of generalized persistence diagrams for 2-filtrations},
  author={Morozov, Dmitriy and Patel, Amit},
  journal={arXiv preprint arXiv:2112.03980},
  year={2021}
}

@article{clause2025meta,
  title={Meta-diagrams for 2-parameter persistence},
  author={Clause, Nate and Dey, Tamal K and M{\'e}moli, Facundo and Wang, Bei},
  journal={Discrete \& Computational Geometry},
  pages={1--27},
  year={2025},
  publisher={Springer}
}

@article{chacholski2024koszul,
  title={Koszul complexes and relative homological algebra of functors over posets},
  author={Chach{\'o}lski, Wojciech and Guidolin, Andrea and Ren, Isaac and Scolamiero, Martina and Tombari, Francesca},
  journal={Foundations of Computational Mathematics},
  pages={1--45},
  year={2024},
  publisher={Springer}
}

@article{demir2022todd,
  title={ToDD: Topological compound fingerprinting in computer-aided drug discovery},
  author={Demir, Andac and Coskunuzer, Baris and Gel, Yulia and Segovia-Dominguez, Ignacio and Chen, Yuzhou and Kiziltan, Bulent},
  journal={Advances in Neural Information Processing Systems},
  volume={35},
  pages={27978--27993},
  year={2022}
}

@article{blumberg2024stability,
  title={Stability of 2-parameter persistent homology},
  author={Blumberg, Andrew J and Lesnick, Michael},
  journal={Foundations of Computational Mathematics},
  volume={24},
  number={2},
  pages={385--427},
  year={2024},
  publisher={Springer}
}

@article{corbet2023computing,
  title={Computing the multicover bifiltration},
  author={Corbet, Ren{\'e} and Kerber, Michael and Lesnick, Michael and Osang, Georg},
  journal={Discrete \& Computational Geometry},
  volume={70},
  number={2},
  pages={376--405},
  year={2023},
  publisher={Springer}
}

@article{buchet2024sparse,
  title={Sparse Higher Order {\v{C}}ech Filtrations},
  author={Buchet, Micka{\"e}l and B Dornelas, Bianca and Kerber, Michael},
  journal={Journal of the ACM},
  volume={71},
  number={4},
  pages={1--23},
  year={2024},
  publisher={ACM New York, NY}
}

@article{hellmer2024density,
  author = {Hellmer, Niklas and Spaliński, Jan},
  title = {Density sensitive bifiltered Dowker complexes via total weight},
  journal = {Homology, Homotopy and Applications},
  volume = {28},
  number = {3},
  pages = {173--204},
  year = {2026},
  month = aug,
  issn = {1532-0081},
  doi = {10.4310/HHA.2026.v28.n3.a8},
  url = {https://doi.org/10.4310/HHA.2026.v28.n3.a8},
  publisher = {International Press of Boston}
}

@article{edelsbrunner2021multi,
  title={The multi-cover persistence of Euclidean balls},
  author={Edelsbrunner, Herbert and Osang, Georg},
  journal={Discrete \& Computational Geometry},
  volume={65},
  pages={1296--1313},
  year={2021},
  publisher={Springer}
}

@article{edelsbrunner2020simple,
  title={A simple algorithm for higher-order {D}elaunay mosaics and alpha shapes},
  author={Edelsbrunner, Herbert and Osang, Georg},
  journal={Algorithmica},
  volume={85},
  number={1},
  pages={277--295},
  year={2023},
  publisher={Springer}
}

@inproceedings{alonso2023filtration,
  title={Filtration-domination in bifiltered graphs},
  author={Alonso, {\'A}ngel Javier and Kerber, Michael and Pritam, Siddharth},
  booktitle={2023 Proceedings of the Symposium on Algorithm Engineering and Experiments ({ALENEX})},
  pages={27--38},
  year={2023},
  publisher={SIAM}
}

@article{alonso2024sparse,
  title={A Sparse Multicover Bifiltration of Linear Size},
  author={Alonso, {\'A}ngel Javier},
  journal={Discrete \& Computational Geometry},
  pages={1--29},
  year={2026},
  publisher={Springer}
}

@article{alonso2024decomposition,
  title={Decomposition of zero-dimensional persistence modules via rooted subsets},
  author={Alonso, {\'A}ngel Javier and Kerber, Michael},
  journal={Discrete \& Computational Geometry},
  pages={1--21},
  year={2024},
  publisher={Springer}
}

@inproceedings{sheehy2012multicover,
  title     = {A Multicover Nerve for Geometric Inference},
  author    = {Sheehy, Donald R.},
  booktitle = {Proceedings of the 24th Canadian Conference on Computational Geometry (CCCG 2012)},
  year      = {2012},
  pages     = {309--314},
  url       = {https://2012.cccg.ca/e-proceedings.pdf}
}

@inproceedings{alonso2024delaunay,
  title={Delaunay bifiltrations of functions on point clouds},
  author={Alonso, {\'A}ngel Javier and Kerber, Michael and Lam, Tung and Lesnick, Michael},
  booktitle={Proceedings of the 2024 Annual ACM-SIAM Symposium on Discrete Algorithms (SODA)},
  pages={4872--4891},
  year={2024},
  organization={SIAM}
}

@article{lesnick2024nerve,
  title={Nerve Models of Subdivision Bifiltrations},
  author={Lesnick, Michael and McCabe, Kenneth},
  journal={arXiv preprint arXiv:2406.07679},
  year={2024}
}

@article{lesnick2024sparse,
  title={Sparse Approximation of the Subdivision-{R}ips Bifiltration for Doubling Metrics},
  author={Lesnick, Michael and McCabe, Kenneth},
  journal={arXiv preprint arXiv:2408.16716},
  year={2024}
}

@article{blaser2024core,
  title={Core bifiltration},
  author={Blaser, Nello and Brun, Morten and Gardaa, Odin Hoff and Salbu, Lars M},
  journal={Journal of Applied and Computational Topology},
  volume={9},
  number={4},
  pages={30},
  year={2025},
  publisher={Springer}
}

@article{xia2015multidimensional,
  title={Multidimensional persistence in biomolecular data},
  author={Xia, Kelin and Wei, Guo-Wei},
  journal={Journal of computational chemistry},
  volume={36},
  number={20},
  pages={1502--1520},
  year={2015},
  publisher={Wiley Online Library}
}

@article{zhang2024multi,
  title={Multi-Cover Persistence ({MCP})-based machine learning for polymer property prediction},
  author={Zhang, Yipeng and Shen, Cong and Xia, Kelin},
  journal={Briefings in Bioinformatics},
  volume={25},
  number={6},
  year={2024},
  publisher={Oxford University Press}
}

@inproceedings{chen2024emp,
  title={Emp: Effective multidimensional persistence for graph representation learning},
  author={Chen, Yuzhou and Segovia-Dominguez, Ignacio and Akcora, Cuneyt Gurcan and Zhen, Zhiwei and Kantarcioglu, Murat and Gel, Yulia and Coskunuzer, Baris},
  booktitle={Learning on Graphs Conference},
  pages={24--1},
  year={2024},
  organization={PMLR}
}

@inproceedings{chen2022tamp,
  title={TAMP-S2GCNets: coupling time-aware multipersistence knowledge representation with spatio-supra graph convolutional networks for time-series forecasting},
  author={Chen, Yuzhou and Segovia-Dominguez, Ignacio and Coskunuzer, Baris and Gel, Yulia},
  booktitle={International conference on learning representations},
  year={2022}
}

@article{chung2022multi,
  title={A multi-parameter persistence framework for mathematical morphology},
  author={Chung, Yu-Min and Day, Sarah and Hu, Chuan-Shen},
  journal={Scientific reports},
  volume={12},
  number={1},
  pages={6427},
  year={2022},
  publisher={Nature Publishing Group UK London}
}

@article{chung2024morphological,
  title={Morphological multiparameter filtration and persistent homology in mitochondrial image analysis},
  author={Chung, Yu-Min and Hu, Chuan-Shen and Sun, Emily and Tseng, Henry C},
  journal={Plos one},
  volume={19},
  number={9},
  pages={e0310157},
  year={2024},
  publisher={Public Library of Science San Francisco, CA USA}
}

@inproceedings{coskunuzer2024time,
  title={Time-aware knowledge representations of dynamic objects with multidimensional persistence},
  author={Coskunuzer, Baris and Segovia-Dominguez, Ignacio and Chen, Yuzhou and Gel, Yulia R},
  booktitle={Proceedings of the AAAI Conference on Artificial Intelligence},
  volume={38},
  pages={11678--11686},
  year={2024}
}

@inproceedings{schiff2020characterizing,
  title={Characterizing the Latent Space of Molecular Deep Generative Models with Persistent Homology Metrics},
  author={Schiff, Yair and Das, Payel and Vijil, Enara and Ramamurthy, Karthikeyan Natesan},
  booktitle={TDA \& Beyond Workshop, Annual Conference on Neural Information Processing Systems},
  year={2020}
}

@article{azumaya1950corrections,
  title={Corrections and supplementaries to my paper concerning {K}rull-{R}emak-{S}chmidt’s theorem},
  author={Azumaya, Gor{\^o}},
  journal={Nagoya mathematical journal},
  volume={1},
  pages={117--124},
  year={1950},
  publisher={Cambridge University Press}
}

@article{bauer2021ripser,
  author    = {Bauer, Ulrich},
  title     = {Ripser: efficient computation of {Vietoris--Rips} persistence barcodes},
  journal   = {Journal of Applied and Computational Topology},
  volume    = {5},
  number    = {3},
  pages     = {391--423},
  year      = {2021},
  publisher = {Springer},
  doi       = {10.1007/s41468-021-00071-5}
}

@InProceedings{bauer2023efficient,
  author =	{Bauer, Ulrich and Lenzen, Fabian and Lesnick, Michael},
  title =	{{Efficient Two-Parameter Persistence Computation via Cohomology}},
  booktitle =	{39th International Symposium on Computational Geometry (SoCG 2023)},
  pages =	{15:1--15:17},
  series =	{Leibniz International Proceedings in Informatics (LIPIcs)},
  ISBN =	{978-3-95977-273-0},
  ISSN =	{1868-8969},
  year =	{2023},
  volume =	{258},
  editor =	{Chambers, Erin W. and Gudmundsson, Joachim},
  publisher =	{Schloss Dagstuhl -- Leibniz-Zentrum f{\"u}r Informatik},
  address =	{Dagstuhl, Germany},
  URL =		{https://drops.dagstuhl.de/entities/document/10.4230/LIPIcs.SoCG.2023.15},
  doi =		{10.4230/LIPIcs.SoCG.2023.15},
}

@article{bauer2025multi,
  title={Multi-parameter Persistence Modules are Generically Indecomposable},
  author={Bauer, Ulrich and Scoccola, Luis},
  journal={International Mathematics Research Notices},
  volume={2025},
  number={5},
  year={2025},
  publisher={Oxford University Press}
}

@article{bjerkevik2021asymptotic,
  title={Asymptotic improvements on the exact matching distance for $2$-parameter persistence},
  author={Bjerkevik, Havard Bakke and Kerber, Michael},
  journal={Journal of Computational Geometry},
  volume={14},
  number={1},
  pages={309--342},
  year={2023}
}

@article{bjerkevik2025stabilizing,
  title={Stabilizing decomposition of multiparameter persistence modules},
  author={Bjerkevik, H{\aa}vard Bakke},
  journal={Foundations of Computational Mathematics},
  pages={1--60},
  year={2025},
  publisher={Springer}
}

@book{boissonnat2018geometric,
  title={Geometric and topological inference},
  author={Boissonnat, Jean-Daniel and Chazal, Fr{\'e}d{\'e}ric and Yvinec, Mariette},
  volume={57},
  year={2018},
  publisher={Cambridge University Press}
}

@article{botnan2020decomposition,
  title={Decomposition of persistence modules},
  author={Botnan, Magnus and Crawley-Boevey, William},
  journal={Proceedings of the American Mathematical Society},
  volume={148},
  number={11},
  pages={4581--4596},
  year={2020}
}

@inproceedings{botnan2022introduction,
  title={An introduction to multiparameter persistence},
  author={Botnan, Magnus and Lesnick, Michael},
  booktitle={Representations of Algebras and Related Structures: International Conference on Representations of Algebras, ICRA 2020, 9--25 November 2020},
  pages={77--150},
  year={2023},
  organization={EMS Press}
}

@article{botnan2024signed,
  title={Signed barcodes for multi-parameter persistence via rank decompositions and rank-exact resolutions},
  author={Botnan, Magnus Bakke and Oppermann, Steffen and Oudot, Steve},
  journal={Foundations of Computational Mathematics},
  pages={1--60},
  year={2024},
  publisher={Springer}
}

@article{botnan2024bottleneck,
  title={On the bottleneck stability of rank decompositions of multi-parameter persistence modules},
  author={Botnan, Magnus Bakke and Oppermann, Steffen and Oudot, Steve and Scoccola, Luis},
  journal={Advances in Mathematics},
  volume={451},
  pages={109780},
  year={2024},
  publisher={Elsevier}
}

@article{cai2021elder,
  title={Elder-rule-staircodes for augmented metric spaces},
  author={Cai, Chen and Kim, Woojin and M{\'e}moli, Facundo and Wang, Yusu},
  journal={SIAM Journal on Applied Algebra and Geometry},
  volume={5},
  number={3},
  pages={417--454},
  year={2021},
  publisher={SIAM}
}

@article{carriere2020multiparameter,
  title={Multiparameter persistence image for topological machine learning},
  author={Carriere, Mathieu and Blumberg, Andrew},
  journal={Advances in Neural Information Processing Systems},
  volume={33},
  pages={22432--22444},
  year={2020}
}

@inproceedings{cheng2023gril,
  title={Gril: A 2-parameter persistence based vectorization for machine learning},
  author={Cheng, Xin and Mukherjee, Soham and Samaga, Shreyas and Dey, Tamal K},
  booktitle={Proc. ICML 2023 workshop TAGML},
  year={2023}
}

@article{corbet2019kernel,
  title={A kernel for multi-parameter persistent homology},
  author={Corbet, Ren{\'e} and Fugacci, Ulderico and Kerber, Michael and Landi, Claudia and Wang, Bei},
  journal={Computers \& graphics: X},
  volume={2},
  pages={100005},
  year={2019},
  publisher={Elsevier}
}

@Book{deberg2008computation,
  Title                    = {Computational Geometry: Algorithms and Applications},
  Author                   = {de Berg, M. and Cheong, O. and van Kreveld, M. and Overmars, M.},
  Publisher                = {Springer},
  Year                     = {2008}
}

@article{dey2019generalized,
	author = {Dey, Tamal K. and Xin, Cheng},
	da = {2022/02/05},
	doi = {10.1007/s41468-022-00087-5},
	id = {Dey2022},
	isbn = {2367-1734},
	journal = {Journal of Applied and Computational Topology},
	title = {Generalized persistence algorithm for decomposing multiparameter persistence modules},
	ty = {JOUR},
	url = {https://doi.org/10.1007/s41468-022-00087-5},
	year = {2022}}

@inproceedings{dey2025decomposing,
  title={Decomposing Multiparameter Persistence Modules},
  author={Dey, Tamal K and Jendrysiak, Jan and Kerber, Michael},
  booktitle={41st International Symposium on Computational Geometry (SoCG 2025)},
  pages={41--1},
  year={2025},
  organization={Schloss Dagstuhl--Leibniz-Zentrum f{\"u}r Informatik}
}

@article{benjamin2024multiscale,
  title={Multiscale topology classifies cells in subcellular spatial transcriptomics},
  author={Benjamin, Katherine and Bhandari, Aneesha and Kepple, Jessica D and Qi, Rui and Shang, Zhouchun and Xing, Yanan and An, Yanru and Zhang, Nannan and Hou, Yong and Crockford, Tanya L and others},
  journal={Nature},
  volume={630},
  number={8018},
  pages={943--949},
  year={2024},
  publisher={Nature Publishing Group UK London}
}

@article{blanchette_brustle_hanson_2022, 
  title={Homological approximations in persistence theory}, 
  DOI={10.4153/S0008414X22000657}, journal={Canadian Journal of Mathematics}, 
  publisher={Canadian Mathematical Society}, 
  author={Blanchette, Benjamin and Brüstle, Thomas and Hanson, Eric J.}, 
  year={2022}, 
  pages={1–38}
}

@incollection{blanchette2023exact,
  title={Exact structures for persistence modules},
  author={Blanchette, Benjamin and Br{\"u}stle, Thomas and Hanson, Eric J},
  booktitle={Representations of Algebras and Related Topics},
  year={2025},
  pages={121-160},
  publisher={European Mathematical Society}
}

@article{de2011dualities,
  title={Dualities in persistent (co) homology},
  author={De Silva, Vin and Morozov, Dmitriy and Vejdemo-Johansson, Mikael},
  journal={Inverse problems},
  volume={27},
  number={12},
  pages={124003},
  year={2011}
}

@Article{carlsson2009topology,
  Title                    = {Topology and data},
  Author                   = {Carlsson, G.},
  Journal                  = {Bulletin of the American Mathematical Society},
  Year                     = {2009},
  Number                   = {2},
  Pages                    = {255--308},
  Volume                   = {46}
}

@Article{carlsson2010multiparameter,
  Title                    = {Multiparameter hierarchical clustering methods},
  Author                   = {Carlsson, G. and M{\'e}moli, F.},
  Journal                  = {Classification as a Tool for Research},
  Year                     = {2010},
  Pages                    = {63--70},

  Publisher                = {Springer}
}

@Article{carlsson2009computing,
  Title                    = {{Computing multidimensional persistence}},
  Author                   = {Carlsson, G. and Singh, G. and Zomorodian, A.},
  Journal                  = {Algorithms and Computation},
  Year                     = {2009},
  Pages                    = {730--739},

  Publisher                = {Springer}
}

@Article{carlsson2009theory,
  Title                    = {{The theory of multidimensional persistence}},
  Author                   = {Carlsson, G. and Zomorodian, A.},
  Journal                  = {Discrete and Computational Geometry},
  Year                     = {2009},
  Number                   = {1},
  Pages                    = {71--93},
  Volume                   = {42},

  ISSN                     = {0179-5376},
  Publisher                = {Springer}
}

@book{carlsson2021topological,
  title={Topological data analysis with applications},
  author={Carlsson, Gunnar and Vejdemo-Johansson, Mikael},
  year={2021},
  publisher={Cambridge University Press}
}

@Article{cerri2013betti,
  Title                    = {Betti numbers in multidimensional persistent homology are stable functions},
  Author                   = {Cerri, A. and Di Fabio, B. and Ferri, M. and Frosini, P. and Landi, C.},
  Journal                  = {Mathematical Methods in the Applied Sciences},
  Year                     = {2013},
  Number                   = {12},
  Pages                    = {1543--1557},
  Volume                   = {36},

  Publisher                = {Wiley Online Library}
}

@article{cerri2011new,
  title={A new algorithm for computing the 2-dimensional matching distance between size functions},
  author={Biasotti, Silvia and Cerri, Andrea and Frosini, Patrizio and Giorgi, Daniela},
  journal={Pattern Recognition Letters},
  volume={32},
  number={14},
  pages={1735--1746},
  year={2011},
  publisher={Elsevier}
}

@Article{chazal2011geometric,
  Title                    = {Geometric inference for probability measures},
  Author                   = {Chazal, F. and Cohen-Steiner, D. and M{\'e}rigot, Q.},
  Journal                  = {Foundations of Computational Mathematics},
  Year                     = {2011},
  Pages                    = {1--19},

  Publisher                = {Springer}
}

@article{chazal2021introduction,
  title={An introduction to topological data analysis: fundamental and practical aspects for data scientists},
  author={Chazal, Fr{\'e}d{\'e}ric and Michel, Bertrand},
  journal={Frontiers in artificial intelligence},
  volume={4},
  pages={667963},
  year={2021},
  publisher={Frontiers Media SA}
}

@article{chazelle1992algorithm,
    author = {Chazelle, Bernard and Edelsbrunner, Herbert},
    title = {An optimal algorithm for intersecting line segments in the plane},
    year = {1992},
    issue_date = {Jan. 1992},
    publisher = {Association for Computing Machinery},
    address = {New York, NY, USA},
    volume = {39},
    number = {1},
    issn = {0004-5411},
    doi = {10.1145/147508.147511},
    journal = {J. ACM},
    pages = {1–54}
}

@article{bauer2026dualities,
  title={Dualities in Multiparameter Persistence},
  author={Bauer, Ulrich and Lenzen, Fabian and Lesnick, Michael},
  journal={arXiv preprint arXiv:2603.18224},
  year={2026}
}

@Techreport{christofides1976,
  author      = {Christofides, Nicos},
  title       = {Worst-case analysis of a new heuristic for the travelling salesman problem},
  institution = {Graduate School of Industrial Administration, Carnegie-Mellon University},
  year        = {1976},
  number      = {Management Sciences Research Report No. 388},
  address     = {Pittsburgh, PA}
}

@Article{cohen2007stability,
  Title                    = {{Stability of persistence diagrams}},
  Author                   = {Cohen-Steiner, D. and Edelsbrunner, H. and Harer, J.},
  Journal                  = {Discrete and Computational Geometry},
  Year                     = {2007},
  Number                   = {1},
  Pages                    = {103--120},
  Volume                   = {37},

  ISSN                     = {0179-5376},
  Publisher                = {Springer}
}

@InProceedings{cohen2006vines,
  Title                    = {Vines and vineyards by updating persistence in linear time},
  Author                   = {Cohen-Steiner, D. and Edelsbrunner, H. and Morozov, D.},
  Booktitle                = {Proceedings of the twenty-second annual symposium on Computational geometry},
  Year                     = {2006},
  Organization             = {ACM},
  Pages                    = {119--126}
}

@Book{cormen2009introduction,
  Title                    = {Introduction to Algorithms},
  Author                   = {Cormen, T. and Leiserson, C. and Rivest, R. and Stein, C.},
  Publisher                = {MIT Press},
  Year                     = {2009},
}

@Book{cox1998using,
  Title                    = {Using algebraic geometry},
  Author                   = {Cox, D.A. and Little, J.B. and O'Shea, D.},
  Publisher                = {Springer Verlag},
  Year                     = {1998},
  Volume                   = {185}
}

@article{crawley2012decomposition,
  title={Decomposition of pointwise finite-dimensional persistence modules},
  author={Crawley-Boevey, William},
  journal={Journal of Algebra and Its Applications},
  volume={14},
  number={05},
  pages={1550066},
  year={2015},
  publisher={World Scientific}
}

@Article{de2022valueoffset,
  title = {Value-Offset Bifiltrations for Digital Images},
  author = {Anway De and Thong Vo and Matthew Wright},
  journal = {Computational Geometry},
  volume = {109},
  pages = {101939},
  year = {2023},
  issn = {0925-7721},
  doi = {https://doi.org/10.1016/j.comgeo.2022.101939} }

@book{dey2022computational,
  title={Computational topology for data analysis},
  author={Dey, Tamal Krishna and Wang, Yusu},
  year={2022},
  publisher={Cambridge University Press}
}

@misc{DONUT,
  author = {Giunti, Barbara and Lazovskis, J{\=a}nis and Rieck, Bastian},
  title  = {
    {DONUT}: {D}atabase of {O}riginal \& {N}on-{T}heoretical {U}ses of {T}opology
  },
  note   = {\url{https://donut.topology.rocks}},
  year   = {2022},
  key    = {DONUT},
}

@Book{edelsbrunner2010computational,
  Title                    = {Computational topology: an introduction},
  Author                   = {Edelsbrunner, H. and Harer, J.},
  Publisher                = {American Mathematical Society},
  Year                     = {2010}
}

@Article{edelsbrunner2002topological,
  Title                    = {{Topological persistence and simplification}},
  Author                   = {Edelsbrunner, H. and Letscher, D. and Zomorodian, A.},
  Journal                  = {Discrete and Computational Geometry},
  Year                     = {2002},
  Number                   = {4},
  Pages                    = {511--533},
  Volume                   = {28},

  ISSN                     = {0179-5376},
  Publisher                = {Springer}
}

@Article{frosini1999size,
  Title                    = {Size homotopy groups for computation of natural size distances},
  Author                   = {Frosini, P. and Mulazzani, M.},
  Journal                  = {Bulletin of the Belgian Mathematical Society Simon Stevin},
  Year                     = {1999},
  Number                   = {3},
  Pages                    = {455--464},
  Volume                   = {6},

  Publisher                = {Brussels: Societe Mathematique de Belgique, c1994-}
}

@article{la1998strategies,
  title={Strategies for computing minimal free resolutions},
  author={La Scala, Roberto and Stillman, Michael},
  journal={Journal of Symbolic Computation},
  volume={26},
  number={4},
  pages={409--431},
  year={1998},
  publisher={Elsevier}
}

@article{fugacci2023compression,
  author = {Ulderico Fugacci and Michael Kerber and Alexander Rolle},
  title = {Compression for 2-parameter persistent homology},
  journal = {Computational Geometry},
  volume = {109},
  pages = {101940},
  year = {2023},
  issn = {0925-7721},
  doi = {https://doi.org/10.1016/j.comgeo.2022.101940},
  url = {https://www.sciencedirect.com/science/article/pii/S0925772122000839},
}

@article{hoogeveen1991analysis,
  title={Analysis of {C}hristofides' heuristic: Some paths are more difficult than cycles},
  author={Hoogeveen, JA},
  journal={Operations Research Letters},
  volume={10},
  number={5},
  pages={291--295},
  year={1991},
  publisher={Elsevier}
}

@InProceedings{kerber2019matching,
  author =	{Kerber, Michael and Lesnick, Michael and Oudot, Steve},
  title =	{{Exact Computation of the Matching Distance on 2-Parameter Persistence Modules}},
  booktitle =	{35th International Symposium on Computational Geometry (SoCG 2019)},
  pages =	{46:1--46:15},
  series =	{Leibniz International Proceedings in Informatics (LIPIcs)},
  ISBN =	{978-3-95977-104-7},
  ISSN =	{1868-8969},
  year =	{2019},
  volume =	{129},
  editor =	{Barequet, Gill and Wang, Yusu},
  publisher =	{Schloss Dagstuhl -- Leibniz-Zentrum f{\"u}r Informatik},
  address =	{Dagstuhl, Germany},
  URL =		{https://drops.dagstuhl.de/entities/document/10.4230/LIPIcs.SoCG.2019.46},
  URN =		{urn:nbn:de:0030-drops-104505},
  doi =		{10.4230/LIPIcs.SoCG.2019.46}
}

@inproceedings{kerber2021fast,
  author    = {Michael Kerber and Alexander Rolle},
  title     = {Fast Minimal Presentations of Bi-graded Persistence Modules},
  booktitle = {2021 Proceedings of the Symposium on Algorithm Engineering and Experiments ({ALENEX})},
  pages     = {207--220},
  publisher = {SIAM},
  year      = {2021},
  doi       = {10.1137/1.9781611976472.16},
  url       = {https://epubs.siam.org/doi/abs/10.1137/1.9781611976472.16}
}

@article{kim2021generalized,
  title={Generalized persistence diagrams for persistence modules over posets},
  author={Kim, Woojin and M{\'e}moli, Facundo},
  journal={Journal of Applied and Computational Topology},
  volume={5},
  number={4},
  pages={533--581},
  year={2021},
  publisher={Springer}
}

@Book{kreuzer2005computational,
  Title                    = {Computational Commutative Algebra 2},
  Author                   = {Kreuzer, M. and Robbiano, L.},
  Publisher                = {Springer},
  Year                     = {2005}
}

@Book{kreuzer2000computational,
  Title                    = {Computational Commutative Algebra 1},
  Author                   = {Kreuzer, M. and Robbiano, L.},
  Publisher                = {Springer},
  Year                     = {2000}
}

@article{xian2022capturing,
  title={Capturing dynamics of time-varying data via topology},
  author={Xian, Lu and Adams, Henry and Topaz, Chad M and Ziegelmeier, Lori},
  journal={Foundations of Data Science},
  volume={4},
  number={1},
  pages={1--36},
  year={2022},
  publisher={Foundations of Data Science}
}

@article{landi2014rank,
  title={The Rank Invariant Stability via Interleavings},
  author={Landi, Claudia},
  journal={Research in Computational Topology},
  volume={13},
  pages={1},
  year={2018},
  publisher={Springer}
}

@article{lesnick2015theory,
  title={The theory of the interleaving distance on multidimensional persistence modules},
  author={Lesnick, Michael},
  journal={Foundations of Computational Mathematics},
  volume={15},
  number={3},
  pages={613--650},
  year={2015},
  publisher={Springer}
}

@article{lesnick2022minimal,
  title = {Computing Minimal Presentations and Bigraded {B}etti Numbers of 2-Parameter Persistent Homology},
  author = {Michael Lesnick and Matthew Wright},
  journal = {SIAM Journal on Applied Algebra and Geometry},
  volume = 6,
  number = 2,
  pages = {267-298},
  year = {2022},
  doi = {10.1137/20M1388425}
}

@article{lesnick2015interactive,
    title={Interactive Visualization of 2-D Persistence Modules}, 
    author={Michael Lesnick and Matthew Wright},
    journal={arXiv preprint arXiv:1512.00180},  
    year={2015},
    eprint={1512.00180},
    archivePrefix={arXiv},
    primaryClass={math.AT},
    url={https://arxiv.org/abs/1512.00180}, 
}

@article{loiseaux2024multipers,
  title={multipers: Multiparameter persistence for machine learning},
  author={Loiseaux, David and Schreiber, Hannah},
  journal={Journal of Open Source Software},
  volume={9},
  number={103},
  pages={6773},
  year={2024}
}

@article{loiseaux2023stable,
  title={Stable vectorization of multiparameter persistent homology using signed barcodes as measures},
  author={Loiseaux, David and Scoccola, Luis and Carri{\`e}re, Mathieu and Botnan, Magnus Bakke and Oudot, Steve},
  journal={Advances in Neural Information Processing Systems},
  volume={36},
  pages={68316--68342},
  year={2023}
}

@article{loiseaux2023framework,
  title={A framework for fast and stable representations of multiparameter persistent homology decompositions},
  author={Loiseaux, David and Carri{\`e}re, Mathieu and Blumberg, Andrew},
  journal={Advances in Neural Information Processing Systems},
  volume={36},
  pages={35774--35798},
  year={2023}
}

@article{luo2024warmstarts,
  title = {Accelerating iterated persistent homology computations with warm starts},
  journal = {Computational Geometry},
  volume = {120},
  pages = {102089},
  year = {2024},
  issn = {0925-7721},
  doi = {https://doi.org/10.1016/j.comgeo.2024.102089},
  url = {https://www.sciencedirect.com/science/article/pii/S0925772124000117},
  author = {Yuan Luo and Bradley J. Nelson},
}

@inproceedings{morozov2024computing,
  title={Computing {B}etti Tables and Minimal Presentations of Zero-Dimensional Persistent Homology},
  author={Morozov, Dmitriy and Scoccola, Luis},
  booktitle={41st International Symposium on Computational Geometry (SoCG 2025)},
  pages={69:1--69:16},
  year={2025},
  publisher={Schloss Dagstuhl--Leibniz-Zentrum f{\"u}r Informatik}
}

@misc{luo2026computing,
  title={Computing {B}etti Tables and Minimal Presentations of Zero-dimensional Persistent Homology},
  author={Luo, Yuan and Morozov, Dmitriy and Scoccola, Luis},
  year={2026},
  month=feb,
  howpublished={\href{https://arxiv.org/abs/2410.22242v3}{arXiv:2410.22242v3}},
  note={Extended version of the SoCG 2025 paper}
}

@book{oudot2015persistence,
  title={Persistence theory: from quiver representations to data analysis},
  author={Oudot, Steve Y},
  volume={209},
  year={2015},
  publisher={American Mathematical Society Providence}
}

@article{oudot2024stability,
  title={On the stability of multigraded {B}etti numbers and {H}ilbert functions},
  author={Oudot, Steve and Scoccola, Luis},
  journal={SIAM Journal on Applied Algebra and Geometry},
  volume={8},
  number={1},
  pages={54--88},
  year={2024},
  publisher={SIAM}
}

@article{patel2018generalized,
  title={Generalized persistence diagrams},
  author={Patel, Amit},
  journal={Journal of Applied and Computational Topology},
  volume={1},
  number={3},
  pages={397--419},
  year={2018},
  publisher={Springer}
}

@Book{peeva2011syzygies,
  Title                    = {Graded Syzygies},
  Author                   = {Peeva, Irena},
  Publisher                = {Springer},
  Year                     = {2011},
}

@Article{piekenbrock2024,
  Title = {Move schedules: fast persistence computations in coarse dynamic settings},
  Author = {Piekenbrock, Matthew and Perea, Jose},
  Journal = {Journal of Applied and Computational Topology},
  Year = {2024},
  Pages = {301--345},
  Volume = {8},
  doi = {10.1007/s41468-023-00156-3},
  url = {https://doi.org/10.1007/s41468-023-00156-3},
}

@book{rabadan2019topological,
  title={Topological data analysis for genomics and evolution: topology in biology},
  author={Rabad{\'a}n, Ra{\'u}l and Blumberg, Andrew J},
  year={2019},
  publisher={Cambridge University Press}
}

@Book{riehl2017category,
  Title                    = {Category Theory in Context},
  Author                   = {Emily Riehl},
  Publisher                = {Dover Publications},
  Year                     = {2016}
}

@misc{rivet,
  author       = {{The RIVET Developers}},
  title        = {{RIVET}},
  howpublished = {\url{https://github.com/rivetTDA/rivet/}},
  note         = {Version 1.1.0},
  year         = {2020},
}

@article{scoccola2023persistable,
  author = {Luis Scoccola and Alexander Rolle},  
  title = {Persistable: persistent and stable clustering},
  journal = {Journal of Open Source Software},
  doi = {10.21105/joss.05022},
  url = {https://doi.org/10.21105/joss.05022}, 
  year = {2023}, 
  publisher = {The Open Journal}, 
  volume = {8}, 
  number = {83}, 
  pages = {5022}
}

@Book{toth2017handbook,
  Title                    = {Handbook of discrete and computational geometry},
  Editor                   = {T{\'o}th, Csaba D. and O'Rourke, Joseph and Goodman, Jacob E.},
  Publisher                = {CRC press},
  Year                     = {2017},
  Edition                  = {3rd}
}

@article{vipond2020multiparameter,
  title={Multiparameter persistence landscapes},
  author={Vipond, Oliver},
  journal={Journal of Machine Learning Research},
  volume={21},
  number={61},
  pages={1--38},
  year={2020}
}

@article{vipond2021multi,
  author = {Oliver Vipond  and Joshua A. Bull  and Philip S. Macklin  and Ulrike Tillmann  and Christopher W. Pugh  and Helen M. Byrne  and Heather A. Harrington },
  title = {Multiparameter persistent homology landscapes identify immune cell spatial patterns in tumors},
  journal = {Proceedings of the National Academy of Sciences},
  volume = {118},
  number = {41},
  pages = {e2102166118},
  year = {2021},
  doi = {10.1073/pnas.2102166118},
  URL = {https://www.pnas.org/doi/abs/10.1073/pnas.2102166118},
}

@article{wasserman2018topological,
  title={Topological data analysis},
  author={Wasserman, Larry},
  journal={Annual review of statistics and its application},
  volume={5},
  pages={501-532},
  year={2018},
  publisher={Annual Reviews}
}

@article{Wright_Zheng_2020,
  title = {Topological Data Analysis on Simple {E}nglish {W}ikipedia Articles},
  author = {Wright, Matthew and Zheng, Xiaojun},
  journal = {The PUMP Journal of Undergraduate Research},
  volume = 3,
  year = 2020,
  month = {Dec.},
  pages = {308-328},
  url = {https://journals.calstate.edu/pump/article/view/2410}
}

@Article{zomorodian2005computing,
  Title                    = {{Computing persistent homology}},
  Author                   = {Zomorodian, A. and Carlsson, G.},
  Journal                  = {Discrete and Computational Geometry},
  Year                     = {2005},
  Number                   = {2},
  Pages                    = {249--274},
  Volume                   = {33},

  ISSN                     = {0179-5376},
  Publisher                = {Springer}
}

\end{document}